\documentclass[a4paper, 11pt]{article}

\usepackage[english]{babel}
\usepackage[latin2]{inputenc}
\usepackage[T1]{fontenc}

\usepackage{verbatim}
\usepackage[usenames]{xcolor}

\usepackage{amssymb}
\usepackage{amsmath}
\usepackage{amsthm}
	
\usepackage[overload]{empheq}

\usepackage[a4paper]{geometry}
\usepackage{caption}
\usepackage{subcaption}
\usepackage{graphicx}
\usepackage[font=normalsize]{caption}

\usepackage{setspace}
\usepackage{indentfirst}
\usepackage{lmodern}
\numberwithin{equation}{section}

\theoremstyle{plain}
\newtheorem{defi}{Definition}[section]
\newtheorem{theo}{Theorem}[section]
\newtheorem{prop}{Proposition}[section]

\newtheorem{lem}{Lemma}[section]
\newtheorem{cor}{Corollary}[section]
\newtheorem{rem}{Remark}[section]

\newcommand{\prove}{\noindent\textbf{Proof: }}
\newcommand{\quod}{\hfill$\square$}

\newcommand{\brck}[1]{{\left\vert\kern-0.25ex\left\vert\kern-0.25ex\left\vert #1 \right\vert\kern-0.25ex\right\vert\kern-0.25ex\right\vert}}

\usepackage{paralist}

\usepackage{enumerate}
\usepackage{enumitem}
\setenumerate{topsep=-2pt, itemsep=-1pt}

\usepackage[unicode, hidelinks, hyperfootnotes=false, pagebackref=false]{hyperref}

\usepackage{verbatim}

\begin{document}

\begin{center}

\textbf{\Large Numerical analysis of the weakly nonlinear Boussinesq
  system with a freely moving body on the bottom}\\
\vspace{0.3 cm}
\textbf{Krisztian Benyo}\footnote[2]{Institut de Math\'ematiques de Bordeaux, Universit\'e de Bordeaux, France. (\texttt{krisztian.benyo@math.u-bordeaux.fr})}\\

\end{center}

\vspace{1 cm}
\begin{center}\parbox{14cm}{
\noindent
\textsc{Abstract:} In this study, the numerical analysis of a specific fluid-solid
interaction problem is detailed. The weakly nonlinear Boussinesq
system is considered with the addition of a solid object lying
on the flat bottom, allowed to move horizontally under the pressure forces
created by the waves. We present an accurate finite
difference scheme for this physical model, finely tuned to preserve
important features of the original coupled system: nonlinear effects
for the waves, energy dissipation due to the frictional movement of the
solid. The moving bottom case is compared with a system where the
same object is fixed to the bottom in order to observe the qualitative
and quantitative differences in wave transformation. In particular a
loss of wave amplitude is observed. The influence of the friction on the whole system is also
measured, indicating differences for small and large
coefficients of friction. Overall, hydrodynamic damping effects
reminiscent to the dead-water phenomenon can be established.}
\end{center}

\section*{Introduction}

Understanding and predicting surface wave formations and their
evolution has been one of the central elements of coastal engineering
and oceanography for the past few decades. In
$1871$ Boussinesq introduced the first depth averaged model
(\cite{boussinesqpapa}). It
originally described a physical situation with horizontal bottom and was later generalized for
variable depths by Peregrine (\cite{peregrine}). These kind of models
have played a crucial role in water wave modeling, especially in
shallow water regions (such as the shoaling zone for a wave). From a
mathematical point of view, these equations arise as shallow water
asymptotic limits of the full water waves problem (for a thorough
discussion on the subject, please refer to \cite{lannesbible}) and
incorporate (weakly) dispersive and (weakly) nonlinear
effects. Nowadays many reformulations and generalizations of the
original Boussinesq system exist: for example Nwogu's extended equation (\cite{nwoguext}) or the abcd-system introduced in \cite{bonachen}, to cite a few of the most important ones.

Implementing a moving bottom in numerical models has its own challenges, experiments and tentatives to adapt
existing (fixed bottom) shallow water models for a moving bottom
regime appear from time to time in the literature. After observing successively
generated solitary waves due to a disturbance in the bottom topography
advancing at critical speed (\cite{wusolitontrain}) Wu et al. formally
derived a set of generalized channel type Boussinesq systems
(\cite{wutengwaterwaves}), their work was extended later on in a
formal study on more general long wave regimes
(\cite{chengeneralboussinesq}). Tsunami research has also proved to be
a main motivating factor with the consideration of water waves type
problems with a moving bottom, for the full justification of the general
problem in shallow water asymptotic regimes, we refer to \cite{iguchibottommove}.

The scope of the current work revolves around a model in which the
bottom is still moving, but its movement is not prescribed, instead it
is generated by the wave motion, introduced in the author's previous
work (\cite{egoboost}). The main idea is to place a freely moving
object on the bottom of the fluid domain, its movement being governed
by the wave motions. Physical motivations for such a model stem, for
example, from marine energy engineering, most notably submerged wave
energy converters (submerged pressure differential devices,
\cite{LehmannTheWC} and references therein) and oscillating wave surge
converters (WaveRollers and Submerged plate devices,
\cite{waveenergyconverters}). Originating from the same context, one
can also handle the ``complementary'' problem in which the waves are
interacting with a freely moving floating structure. The theoretical
background of such problems have been extensively analyzed in recent
years, see for example \cite{lannesfloating}, \cite{edofloating}.

Although the numerical study of immersed structures is not entirely new, our approach is heavily based on a theoretical analysis of the general system, rather than experimental or numerical considerations usually present in the corresponding bibliography.  An approach with numerous physical and biological applications was
developed by Cottet et al. (\cite{cottetmaitrenumerics},
\cite{cotteteulerian}) based on a level set formulation, adapting an
immersed boundary method and general elasticity theory (see also
\cite{celineengineer}). Modeling underwater landslides provides for an
excellent example of such systems, we refer to \cite{kalischbottom}
and references therein for recent developpments. Numerical models adapted to tsunami generation
due to seabed deformations were presented for example in
\cite{guyennespectrum} or \cite{mitsotakistsunami}. From a control theory point of view,
Zuazua et al. (\cite{zoazoazoo}) performed an analytical and numerical
analysis on underwater wave generator models. Perhaps one of the most
relevant existing studies concerns a submerged spring-block model and the associated experimental and numerical observations (\cite{abadiedeplace}, \cite{abadibadou}).

The structure of the article is as follows. After a brief
introduction, we present the weakly nonlinear Boussinesq system in a
fluid domain with a flat bottom topography and with a solid object
lying on the bottom, capable of moving horizontally under the pressure
forces created by the waves. Following our previous work, this coupled
system admits a unique solution for a long time scale (\cite{egoboost}).

In the third section we detail the finite difference numerical scheme adapted to this system. We elaborate a fourth order accurate staggered grid system for the variables concerning the movement of the fluid, following the footsteps of Lin and Man (\cite{linman}). As for the time discretization, an adapted fourth order accurate Adams--Bashforth predictor-corrector method is implemented, incorporating the discretized ordinary differential equation characterizing the solid displacement via a modified central finite difference scheme. We end this section with some remarks on the boundary conditions implemented and on certain useful properties of our adaptation.

Section $4$ details the numerical experiments concerning the
model. The convergence of the finite difference scheme is measured to
be almost of order $4$ in time and in space for a flat bottom as well
as for large coefficients of friction, greatly improving the
reference staggered grid model (of order $2$ only) in \cite{linman} over a flat bottom. An order $3$ mesh
convergence and an order $2$ convergence in time is observed for small
coefficients of friction. The transformation of a passing wave over the solid is
detailed in various different physical regimes. Wave shoaling effects
are examined and compared to a system with the same parameters
admitting a fixed solid object on the bottom instead of a freely
moving one. The effects of the friction on the motion of the solid are
also measured, revealing that the solid comes to a halt after the wave
has passed over it. Measuring the solid motion also indicates hydrodynamical damping effects reminiscent to the ones
attributed to dead-water phenomena, closely tied to internal wave generation
(for more details, we refer to \cite{ekmandead}, \cite{mercierres}, \cite{vincentintern}). Long term effects by a wave train test are also presented at the end of the section.

\section{The governing equations}\label{sec_equations} 

\subsection{The physical regime}

We are going to work in two spatial dimensions, and we will
reference the horizontal coordinate by $x$ and the vertical coordinate
by $z$. The time parameter shall be $t\in\mathbb{R}^+$. The physical
domain of the fluid is
\begin{equation*}
\Omega_t=\left\{(x,z)\in\mathbb{R}\times\mathbb{R}\,:\,-H_0+b(t,x)<z<\zeta(t,x)\right\},
\end{equation*}
where $H_0$ is a parameter of the system measuring the base water
depth and the functions $\zeta(t,x)$ and $b(t,x)$ stand for the free
surface elevation and the bottom topography variation
respectively (see Figure \ref{fig_domainwithsolid}).

The solid on the bottom is supposed to be moving only in
a horizontal direction, its displacement vector is denoted by
$X(t)$, consequently its velocity is given by $v(t) =
\dot{X}(t)$. Therefore we have that
\begin{equation}\label{form_bottomfunction}
b(t,x)=\mathfrak{b}\left(x-X(t)\right),
\end{equation}
with $\mathfrak{b}$ corresponding to the initial state of the solid,
at $t=0$, thus assuming that $X(0)=0$. This function $\mathfrak{b}$
is of class $\mathcal{C}^\infty(\mathbb{R})$ and compactly supported.

\begin{center}
\captionsetup{type=figure}
\includegraphics[width=0.8\linewidth]{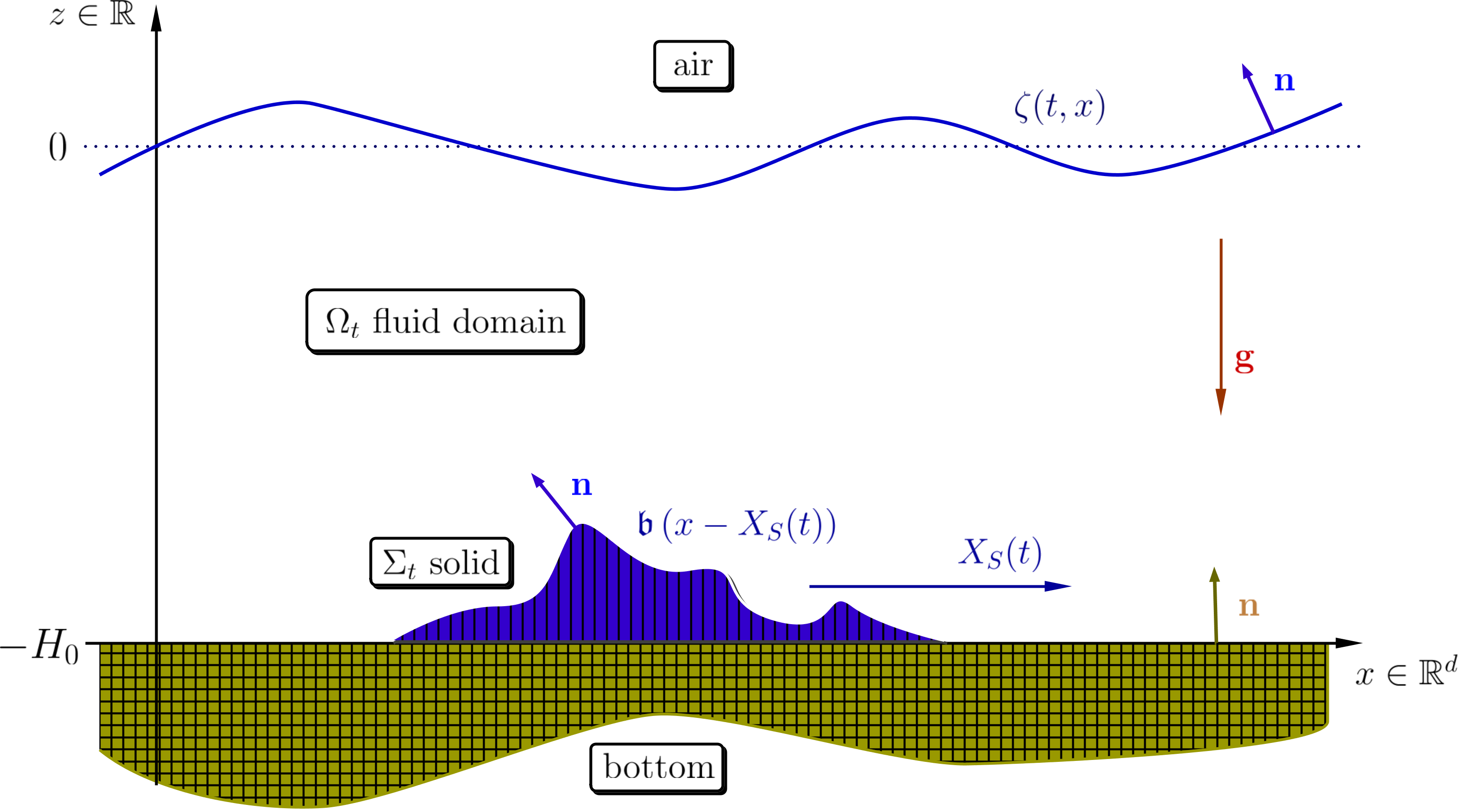}
\captionof{figure}{The coupled water waves setting in the presence of a solid}\label{fig_domainwithsolid}
\end{center}

The fluid is assumed to be homogeneous, inviscid, incompressible, and
irrotational. Its dynamics in general are described by the full
water waves problem, based on the free surface Euler equations (see
for instance \cite{lannesbible}).

The solid is supposed to be rigid and homogeneous with a given mass
$M$, in frictional contact with the flat bottom of the
domain. Its movement is governed by Newton's second law.

Since the weakly nonlinear Boussinesq equation is a shallow water
asymptotic equation, we have to introduce characteristic scales of the
problem. First of all $H_0$ denotes the base water depth, and $L$, the
characteristic horizontal scale of the wave motion as well as the
solid size. Moreover we have $a_{surf}$, the order of the free surface
amplitude, and $a_{bott}$, the characteristic height of the solid.

Using these quantities, we can introduce several dimensionless quantities:
\begin{itemize}
	\item shallowness parameter $\mu=\dfrac{H_0^2}{L^2}$,
	\item nonlinearity (or amplitude) parameter $\varepsilon=\dfrac{a_{surf}}{H_0}$,
	\item bottom topography parameter $\beta=\dfrac{a_{bott}}{H_0}$,
\end{itemize}
which will play an important role in the formulation of the governing equations.

\begin{center}
\captionsetup{type=figure}
\includegraphics[width=0.7\linewidth]{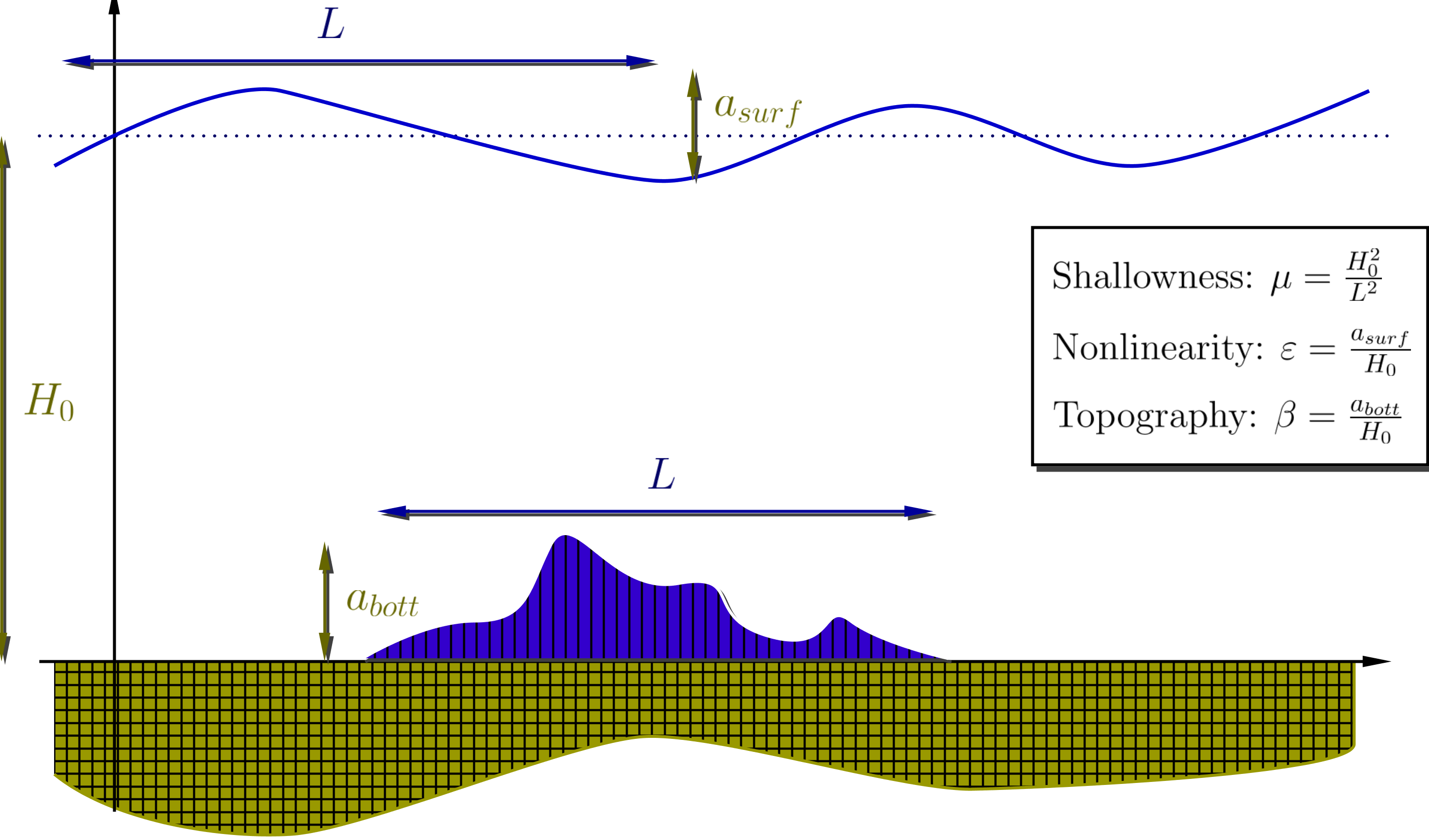}
\captionof{figure}{The characteristic scales of the coupled water waves problem}\label{fig_domainasymptotics}
\end{center}

\subsection{The coupled Boussinesq system}\label{sec_boussrewrite}

We consider here a weakly nonlinear regime, i.e. we assume that
$\varepsilon=\mathcal{O}(\mu)$ and we suppose that the scale of the
solid is relatively small (meaning $\beta=\mathcal{O}(\mu)$). The fluid is then governed by a Boussinesq system. Thus, the asymptotic regime writes as
\begin{equation}\tag{BOUS}\label{hyp_boussinesqregime}
0\leqslant\mu\leqslant\mu_{max}\ll 1,\quad\varepsilon=\mathcal{O}(\mu),\;\beta=\mathcal{O}(\mu).
\end{equation}

Under the Boussinesq regime (\ref{hyp_boussinesqregime}) the governing
fluid equations at order $\mathcal{O}(\mu^2)$ take the form of
\begin{equation}\label{eq_boussinesq}
\begin{cases}
\partial_t\zeta+\partial_x(h\overline{V}) = \dfrac{\beta}{\varepsilon}\partial_t b,\\
\left(1- \displaystyle\frac{\mu}{3}\partial_{xx}\right)\partial_t\overline{V}+\partial_x\zeta+\varepsilon(\overline{V}\cdot\partial_x)\overline{V} = -\displaystyle\frac{\mu\beta}{2\varepsilon}\partial_x\partial_t^2b.
\end{cases}
\end{equation}
with the nondimensionalized fluid height function
$h=1+\varepsilon\zeta-\beta b$ and the vertically averaged horizontal
velocity
\begin{equation*}
\overline{V}(t,x) = \frac{1}{h}\int_{-1+\beta b}^{\varepsilon\zeta}V(t,x,z)\,dz,
\end{equation*}
with $V$ being the horizontal component of the velocity field of the fluid.
In the equation, the variables are the nondimensionalized forms of the corresponding physical quantities described in the previous section. To recover the quantities with the proper dimensions, it is enough to multiply the variable by its corresponding characteristic scale.

\begin{rem}
Without the smallness assumption on $\varepsilon$, and on $\beta$, it
is still possible to perform an asymptotic expansion at
$\mathcal{O}(\mu^2)$. The resulting system is more general than the
Boussinesq system (\ref{eq_boussinesq}) but also more complicated, it
is known as the Serre--Green--Naghdi equations (see for example \cite{iguchigreen}).
\end{rem}

Following the derivation presented in \cite{egoboost}, Newton's second
law for the solid displacement can be written as
\begin{equation}\label{eq_newtonorder1}
\ddot{X}=-\frac{c_{fric}}{\sqrt{\mu}}F_{normal}(X,\zeta)\frac{\dot{X}}{\left|\dot{X}\right|+\delta}+\frac{\varepsilon}{\tilde{M}}\int_{\mathbb{R}}\zeta\partial_xb\,dx,
\end{equation}
with
\begin{equation}
F_{normal}(X,\zeta)=1+\frac{c_{solid}}{\beta}+\frac{\varepsilon}{\beta\tilde{M}}\int_{\operatorname{supp}(\mathfrak{b})+ X}\zeta\,dx,
\end{equation}
where we grouped together the quantities depending only on the
physical parameters of the solid:
\begin{equation}\label{form_solidcoefficient}
c_{solid}=\frac{|\operatorname{supp}(\mathfrak{b})|}{\tilde{M}}\left(\frac{P_{atm}}{\varrho gH_0}+1\right)-\frac{\beta}{\tilde{M}}\int_{\operatorname{supp}(\mathfrak{b})}b\,dx,
\end{equation}
as well as
\begin{equation*}
\tilde{M} = \frac{M}{\varrho L a_{bott}}.
\end{equation*}

Apart from the smallness parameters originating from the system
describing the fluid dynamics (equations (\ref{eq_boussinesq})),
another important parameter in the numerical analysis is the
coefficient of friction $c_{fric}$. The actual measurement of this
coefficient is rather difficult, especially in a complex physical
system such as the current one, mainly because its value depends on
many other physical parameters (like the material structure of the
surfaces, the temperature, the pressure, the velocity of the sliding,
etc.). Generally speaking a coefficient of $10^{-2}\sim 10^{-3}$
corresponds to a relatively frictionless sliding, and values in the
range of $1$ signify an important friction between the contact media.

Moreover, $\delta$ is an artificial parameter to introduce a
uniform dynamical friction law for the movement, its value is taken to
be sufficiently small ($\sim 10^{-10}$) in order to be much smaller
than any physically meaningful solid velocity value. Finally, we have
that $P_{atm}$ stands for
the atmospheric pressure and $g$ the gravitational acceleration constant.

\begin{rem}\label{rem_coefficient}
The first term in $c_{solid}$ corresponds to the effect of the atmospheric and
hydrostatic pressure. It is reasonable to assume that these two
effects should be comparable, which means that for the scale of the
base water depth $H_0$ we can infer that
\begin{equation*}
\frac{P_{atm}}{\varrho g H_0}+1 = \frac{10\textnormal{m}}{H_0}+1,
\end{equation*}
implying that for small depths ($\sim 1$ m) the atmospheric pressure
dominates over even the hydrostatic effects by a factor of $10$, resulting in negligible hydrostatic effects
which is not ideal for our analysis.
\end{rem}

To sum it up, the free surface equations with a solid moving at the bottom in the case of the Boussinesq approximation take the following form
\begin{subequations}\label{eq_systemorder2bous}
\begin{align}[left = {\,\empheqlbrace}]
\begin{split}
&\partial_t\zeta+\partial_x(h\overline{V}) = \frac{\beta}{\varepsilon}\partial_t b,\\
&\left(1-\frac{\mu}{3}\partial_{xx}\right)\partial_t\overline{V}+\partial_x\zeta+\varepsilon(\overline{V}\cdot\partial_x)\overline{V} = -\frac{\mu\beta}{2\varepsilon}\partial_x\partial_{tt}b,
\end{split}  \label{eq_systembous1}\\
&\ddot{X}=-\frac{c_{fric}}{\sqrt{\mu}}F_{normal}(\zeta,X)\frac{\dot{X}}{\left|\dot{X}\right|+\delta}+\frac{\varepsilon}{\tilde{M}}\int_{\mathbb{R}}\zeta\partial_xb\,dx. \label{eq_systembous2}
\end{align}
\end{subequations}

With the numerical scheme in mind, we can rewrite this system in a
more compact form as follows
\begin{subequations}\label{sys_boussstruct}
\begin{align}[left = {\,\empheqlbrace}]
\begin{split}
&\partial_t\zeta = E(\zeta, \overline{V}, X,\dot{X}),\\
&\partial_t\overline{U} = F(\zeta, \overline{V}, X,\dot{X},\ddot{X}),
\end{split} \label{eq_systemstruct1}\\
\begin{split}
&\ddot{X} = G(\zeta, X, \dot{X}),
\end{split}\label{eq_systemstruct2}
\end{align}
\end{subequations}
where
\begin{equation}
\overline{U} = \overline{U}(\overline{V}) = \overline{V}-\frac{\mu}{3}\partial_{xx}\overline{V}.
\end{equation}
The remaining terms are given by
\begin{subequations}
\begin{align}
E(\zeta, \overline{V}, X,\dot{X}) &= -\partial_x(h\overline{V}) -
  \frac{\beta}{\varepsilon}\partial_x\mathfrak{b}(x-X) \dot{X}, \label{eq_bousszeta}\\
F(\zeta, \overline{V}, X,\dot{X},\ddot{X}) &=
                                  -\partial_x\zeta-\frac{\varepsilon}{2}\partial_x(\overline{V}^2)
                                  -\frac{\mu\beta}{2\varepsilon}\partial_{xxx}\mathfrak{b}(x-X)(\dot{X})^2
  +
  \frac{\mu\beta}{2\varepsilon}\partial_{xx}\mathfrak{b}(x-X) \ddot{X}\label{eq_boussvbar}\\
G(\zeta, X,\dot{X}) &=
-\frac{c_{fric}}{\sqrt{\mu}}F_{normal}(\zeta,X)\frac{\dot{X}}{|\dot{X}|+\delta}+\frac{\varepsilon}{\tilde{M}}\int_{\mathbb{R}}\zeta\partial_x\mathfrak{b}(x-X)\,dx
\end{align}
\end{subequations}

\subsection{Relevant properties of the system}

Due to the special structure of our system, let us define the following adapted Banach-space to provide a uniformly formulated energy estimate.
\begin{defi}
The Sobolev type space $\mathcal{X}^s$ is given by
\begin{equation*}
\mathcal{X}^s(\mathbb{R})=\left\{\mathcal{U}=(\zeta,\overline{V})\in L^2(\mathbb{R})\textnormal{ such that }\|\mathcal{U}\|_{\mathcal{X}^s}<\infty\right\},
\end{equation*}
where
\begin{equation*}
\|\mathcal{U}\|_{\mathcal{X}^s}=\|\zeta\|_{H^s}+\|\overline{V}\|_{H^s}+\sqrt{\mu}\|\partial_x\overline{V}\|_{H^s}.
\end{equation*}
\end{defi}
The last term in the $\mathcal{X}^s$ norm appeared due to the necessity to control the dispersive smoothing.

By introducing the wave-structure energy functional
\begin{equation}\label{form_energyvar}
E_B(t)=\frac{1}{2}\int_{\mathbb{R}}\zeta^2\,dx+\frac{1}{2}\int_{\mathbb{R}}h(\overline{V})^2\,dx+\frac{1}{2}\int_{\mathbb{R}}\frac{\mu}{3}h(\partial_x\overline{V})^2\,dx+\frac{1}{2\varepsilon}\left|\dot{X}\right|^2,
\end{equation}
we can establish an $L^2$ type energy estimate for the
coupled system (\ref{eq_systemorder2bous}), from which we are able to deduce a certain control on the velocity of the solid.

We have the following (\cite{egoboost}):
\begin{prop}\label{prop_aprioriestimatebous}
Let $\mu\ll 1$ sufficiently small and $s_0>1$. Then any $\mathcal{U}\in
\mathcal{C}^1([0,T]\times\mathbb{R})\cap\mathcal{C}^1([0,T];H^{s_0})$,
$X\in\mathcal{C}^2([0,T])$ satisfying the coupled system
(\ref{eq_systemorder2bous}), with initial data
$\mathcal{U}(0,\cdot)=\mathcal{U}_{in}\in\mathcal{C}^1(\mathbb{R})\cap
H^{s_0}$ and $(X(0),\dot{X}(0))=(0,v_{S_0})\in\mathbb{R}\times\mathbb{R}$ verifies the energy estimate
\begin{equation}
\sup_{t\in[0,T]}\left\{e^{-\sqrt{\varepsilon}c_0t}E_B(t)\right\}\leqslant
2E_B(0)+2\mu c_0T\|\mathfrak{b}\|_{H^3}^2,
\end{equation}
where
\begin{equation*}
c_0= c(\brck{\mathcal{U}}_{T,W^{1,\infty}},\brck{\mathcal{U}}_{T,H^{s_0}},\|b\|_{W^{3,\infty}}).
\end{equation*}
\end{prop}

\begin{cor}\label{cor_velocityestimate} This energy estimate provides us with a natural control on the solid velocity, namely
\begin{equation}\label{form_veloestim}
\sup_{t\in[0,T]}\left\{e^{-\sqrt{\varepsilon}
    c_0t}\left|\dot{X}(t)\right|^2\right\}\leqslant\varepsilon\|\mathcal{U}_{in}\|_{\mathcal{X}^0}^2+\varepsilon\left|v_{S_0}\right|^2+\varepsilon\mu
c_0T\|\mathfrak{b}\|_{H^3}^2,
\end{equation}
where $c_0$ is as before.
\end{cor}

Since we intend to adapt a similar dissipative property for our
numerical scheme, we present here the brief outline of this proposition.

\noindent\textbf{A sketch of the proof:} It follows a standard energy
estimate argument; we multiply the first equation of
(\ref{eq_systemorder2bous}) by $\zeta$ and the second equation by
$h\overline{V}$, and we integrate over $\mathbb{R}$ with respect to
$x$. 

Handling the right hand side of the first equation is the key to
uniform estimates, since the integral term
\begin{equation*}
-\int_{\mathbb{R}}\zeta\partial_xb\,dx\cdot \dot{X},
\end{equation*}
that appears following the multiplication and the integration is also
present in the right hand side of the solid equation
(\ref{eq_systembous2}) multiplied by $\dot{X}$. Using this equality
one can estimate the integral simply by
$\tilde{M}\varepsilon^{-1}\ddot{X}\dot{X}$ since the new friction term is nonpositive.

The rest of the proof is only straightforward term-by-term estimates,
paying close attention to the higher order derivatives ($(1-\frac{\mu}{3}\partial_{xx})\partial_t\overline{V}$) in the momentum
equation and its source term ($\partial_x\partial_{tt}b$).

The existence and uniqueness theorem for the coupled Boussinesq system then states as
follows (for more details, please refer to \cite{egoboost}):
\begin{theo}
Let us consider the coupled system defined by equations (\ref{eq_systemorder2bous}). Let us suppose that for the initial value $\zeta_{in}$ and $\mathfrak{b}$ the lower bound condition (\ref{genminheightcond}) 
\begin{equation}\label{genminheightcond}
\exists h_{min}>0,\;\forall Y\in\mathbb{R},\;1+\varepsilon\zeta_{in}(Y)-\varepsilon\mathfrak{b}(Y)\geqslant h_{min}
\end{equation}
is satisfied. If the initial values $\zeta_{in}$ and $\overline{V}_{in}$ are in $\mathcal{X}^s(\mathbb{R})$ with $s\in\mathbb{R}$, $s>3/2$, and $ V_{S_0}\in\mathbb{R}$ then there exists a maximal $T_0>0$ independent of $\varepsilon$ such that there is a unique solution 
\begin{equation*}
\begin{aligned}
&(\zeta,\overline{V})\in C\left(\left[0,\frac{T_0}{\sqrt{\varepsilon}}\right];\mathcal{X}^s(\mathbb{R})\right)\cap C^1\left(\left[0,\frac{T_0}{\sqrt{\varepsilon}}\right];\mathcal{X}^{s-1}(\mathbb{R})\right),\\ 
&X\in C^2\left(\left[0,\frac{T_0}{\sqrt{\varepsilon}}\right]\right)
\end{aligned}
\end{equation*}
with uniformly bounded norms for the system (\ref{eq_systemorder2bous}) with initial conditions $(\zeta_{in},\overline{V}_{in})$ and $(0,\sqrt{\varepsilon}V_{S_0})$.
\end{theo}

\section{The discretized model}

In this section we present the elements of the
discretized scheme. We follow the ideas of Lin and Man
(\cite{linman}) for a staggered grid approach. We improve their model
by aiming for a fourth order overall
accuracy for the quantities characterizing the fluid dynamics and the
time stepping in general. In light of the a priori estimate in
Proposition \ref{prop_aprioriestimatebous}. we will examine the
dissipative property of the discretized solid equation (Lemma \ref{lem_discretedissip}).

\subsection{The finite difference scheme on a staggered grid}

Since the solid movement is time dependent only ($X(t)$ does not
depend on the horizontal coordinate $x$), spatial discretization only
concerns the fluid variables, for which we aim for a fourth order
precision.

In their article, Lin and Man obtained a stable
and accurate model for a Boussinesq type
system (the Nwogu equations). They observed good
conservative properties for the fluid system, attributed mainly to the
staggered grid method they implemented, which is an important
factor for long-term measurements, furthermore it is well-adapted to
accurate energy measurements of the system.

We implement their staggered grid method in which the
``scalar'' quantities, such as the surface elevation $\zeta$, the
bottom topography $b$ and the fluid height $h$, are defined on the
grid points, and the ``vectorial''
variable, the averaged velocity $\overline{V}$ (that is still a scalar
since we are only working in one horizontal dimension) is defined on the
mid-points of the mesh. The mesh size will be chosen as $\Delta x$,
numbered by $i=1,2,\ldots N_{space}$. The discrete equation for $\zeta$ (based
on (\ref{eq_bousszeta})) will be defined for mesh points and the
equation for $\overline{V}$ (from equation (\ref{eq_boussvbar})) for
the mid-points. 

In order to be able to do this, we will have to define
the ``scalar'' quantities for the mid-points as well, we shall do so
by a four point centered fourth order interpolation, that is,
\begin{equation*}
\zeta_{i+1/2} = \frac{-\zeta_{i-1}+9\zeta_i + 9\zeta_{i+1}-\zeta_{i+2}}{16}.
\end{equation*}
In the reference article \cite{linman}, only a linear interpolation was
used. Even though it was not mentioned at all, we believe it to be one
of the main reasons for the loss of mesh convergence in their scheme.

For the spatial discretization in general, we chose
fourth order accurate central finite difference schemes for the different
orders of derivatives. In their work, Lin and Man chose only second
order schemes for the higher order derivatives which is another reason
for the resulted loss in mesh convergence in their case. In our implementation even
higher order derivatives are discretized by fourth order schemes. 

Observing the right hand sides of
equations (\ref{eq_systemstruct1}) we may separate four
different types of terms. Once again, we emphasize on the fact that the first equation will act
on mesh points while the second one will be defined on mid-points of
the grid.

\begin{itemize}
\item First order derivative on grid points for a ``scalar'' quantity,
  this concerns the term $\partial_x\mathfrak{b}$, and equivalently
  first order derivative on mid-points for ``vectorial'' variables,
  this concerns the term $\partial_x(\overline{V}^2)$. For this case,
  the classical four point central difference scheme of order $4$ writes as
\begin{equation}
(\partial_x\mathfrak{b})_i = \frac{\mathfrak{b}_{i-2} -
  8\mathfrak{b}_{i-1}+8\mathfrak{b}_{i+1} -
  \mathfrak{b}_{i+2}}{12\Delta x},
\end{equation}
and similarly for the derivative of $\overline{V}^2$ with mid-points.
\item First order derivative on mid-points for a quantity having
  values in grid points,
  this concerns the term $\partial_x(h\overline{V})$, and equivalently
  first order derivative on grid points for quantities having values
  in mid-points,
  this concerns the term $\partial_x\zeta$. For this case,
  the adapted four point central difference scheme of order $4$ writes
  as
\begin{equation}
(\partial_x\zeta)_{i+1/2} = \frac{\zeta_{i-1} -
  27\zeta_{i}+27\zeta_{i+1} -
  \zeta_{i+2}}{24\Delta x},
\end{equation}
and similarly for the derivative of $h\overline{V}$ with mid-points.
\item Second order derivative on mid-points for quantities having
  values at grid points, this concerns both
  $\partial_{xx}\mathfrak{b}$ and $\partial_{xx}\overline{V}$. A classical fourth
  order accurate central finite difference scheme is implemented,
  meaning
\begin{equation}
(\partial_{xx}\mathfrak{b})_{i+1/2} = \frac{-\mathfrak{b}_{i-3/2} +
  16\mathfrak{b}_{i-1/2}-30\mathfrak{b}_{i+1/2}+16\mathfrak{b}_{i+3/2} -
  \mathfrak{b}_{i+5/2}}{12(\Delta x)^2},
\end{equation}
and similarly for $\overline{V}$.
\item Third order derivative on mid-points for a term having
  values on grid points, this concerns
  $\partial_{xxx}\mathfrak{b}$. Once again, a fourth order accurate
  central scheme is applied,
\begin{equation}
(\partial_{xxx}\mathfrak{b})_{i+1/2} = \frac{\mathfrak{b}_{i-5/2}-8\mathfrak{b}_{i-3/2} +
  13\mathfrak{b}_{i-1/2}-13\mathfrak{b}_{i+3/2}+8\mathfrak{b}_{i+5/2} -
  \mathfrak{b}_{i+7/2}}{8(\Delta x)^3}.
\end{equation}
\end{itemize}

The high accuracy guarantees that we can capture
more precisely the nonlinear interaction between the fluid and the
solid without posing problems for the numerical scheme due to the
necessity of information on many grid points, since the solid is
localized to its support (the middle section of the wave tank, as
explained in Section \ref{sec_tank}).

\subsection{Time stepping with Adams--Bashforth}

As elaborated in \cite{linman}, we adapt a fourth order accurate Adams
predictor-corrector method. Starting from the initial condition at time $t=0$, the
first two values of the quantities may be generated by a fourth order
classic Runge--Kutta (RK4) time stepping algorithm. Let us
suppose that currently we are at time step $n\geqslant 2$ and as such,
all information on the main variables ($\zeta$, $\overline{V}$, and
$X$) is known. The method consists of two steps:
\begin{enumerate}
\item First, the predictor step is implemented on the fluid equations
  (equations (\ref{eq_systemstruct1})) by the explicit
  third order Adams--Bashforth scheme
\begin{equation*}
\begin{aligned}
\zeta_i^{n+1^*} &= \zeta_i^n + \frac{\Delta t}{12}\left(23E_i^n-16E_i^{n-1}+5E_i^{n-2}\right),\\
\overline{U}_{i+1/2}^{n+1^*} &= \overline{U}_{i+1/2}^n + \frac{\Delta t}{12}\left(23F _{i+1/2}^n-16F _{i+1/2}^{n-1}+5F _{i+1/2}^{n-2}\right),
\end{aligned}
\end{equation*}
in addition we apply the algorithm for calculating the solid position
(presented in the following section), 
\begin{equation}
X^{n+1} = \overline{G}(\zeta^n,X^n,X^{n-1},X^{n-2}).
\end{equation}
\item With the knowledge of the predicted values, the next step is the
  correction by a fourth order Adams--Moulton method
\begin{equation*}
\begin{aligned}
\zeta_i^{n+1} &= \zeta_i^n + \frac{\Delta t}{24}\left(9E(\zeta_i^{n+1^*}, \overline{V}_{i+1/2}^{n+1^*}, X^{n+1})+ 19E_i^n-5E_i^{n-1}+E_i^{n-2}\right),\\
\overline{U}_{i+1/2}^{n+1} &= \overline{U}_{i+1/2}^n + \frac{\Delta
  t}{24}\left(9F(\zeta_i^{n+1^*}, \overline{V}_{i+1/2}^{n+1^*}, X^{n+1}) +19F _{i+1/2}^n-5F _{i+1/2}^{n-1}+F
  _{i+1/2}^{n-2}\right).
\end{aligned}
\end{equation*}
\end{enumerate}

\begin{rem}
Additionally, the predictor-corrector method can be iterated to guarantee
even more accuracy for the algorithm.
\end{rem}

We remark that in \cite{linman}, the same algorithm was used for the
determination of the fluid quantities, and in the end the method was
presented, without any explications, as a method with a convergence of order $2$ instead of the theoretical
order $4$ in the time variable. One of the reasons for their loss of
convergence can be the negative effects of a second order
discretization in the spatial variable on a fourth order in time
scheme. Another likely scenario would entail a not too accurate
generation of the initial steps, indeed the predictor-corrector
algorithm requires the first two steps to be generated by another
algorithm, implying that if one were to use less accurate methods for
the initial steps, that could negatively influence the convergence in
time. Since in their case, it was not described how they obtained the
first few steps, this can not be excluded. 

\subsection{Time discretization for the solid movement}

Equation (\ref{eq_systemstruct2}) is an ODE
which involves no further spatial discretization. We shall discretize
it in time, considering a time step $\Delta t$ indexed by
$n=1,2,\ldots N_{time}$. We have to be careful though, since this
equation is coupled to the first two equations of the system, having
the source terms depending on $X$, $\dot{X}$, and $\ddot{X}$ as
well, and as such it is incorporated in the Adams scheme presented in
the previous section.

Notice the presence of integrals of the fluid variables in the
expression $G(\zeta,X,\dot{X})$ which implies some
restrictions for calculating the numerical integral since these
variables are only known for grid and mid-points. For further details,
please refer to Section \ref{sec_remarks}.

Another remark
concerns the order of magnitude of the two terms of $G$. We are
working with a regime where the shallowness parameter $\mu$ is supposed to
be small, meaning that the first term is at least of order
$\mathcal{O}(\mu^{-1/2})$ while the second term is small, of
order $\mathcal{O}(\mu)$. As explained in Proposition
\ref{prop_aprioriestimatebous}. we have a strong control on the solid
velocity, moreover good dissipative properties
for the coupled system. We want to achieve a similar property for the discretized
system, that is a similar dissipation of the discrete energy which in
turn ensures that oscillations or other instabilities do not appear
in the simulation.

First of all, we can write that
\begin{equation}\label{eq_solidprediscr}
G(\zeta,X,\dot{X}) =
-C(\zeta,X)\frac{\dot{X}}{|\dot{X}|+\delta} +
\overline{C}(\zeta,X),
\end{equation}
where we introduced the following two quantities
\begin{subequations}
\begin{align}
C(\zeta,X) &= \frac{c_{fric}}{\sqrt{\mu}}\left(1+\frac{c_{solid}}{\beta}+\frac{\varepsilon}{\tilde{M}\beta}\int_{\operatorname{supp}(\mathfrak{b})+
    X}\zeta\,dx\right),\\
\overline{C}(\zeta,X) &=
                         \frac{\varepsilon}{\tilde{M}}\int_{\mathbb{R}}\zeta\partial_x\mathfrak{b}(x-X)\,dx.
\end{align}
\end{subequations}
 
We wish to construct an appropriate numerical scheme. Let us suppose that we
are at time step $n$, so that quantities $\zeta^n$, $\overline{V}^n$,
and $X^n$ are known up until the index $n$. This implies that the constants
$C^n=C(\zeta^n,X^n)$ and $\overline{C}^n=\overline{C}(\zeta^n,X^n)$
are also known (since they do not involve any time differentiation).

We base our discretization on the reformulation
(\ref{eq_solidprediscr}) of the equation at hand. Let us apply a
second order accurate central finite difference scheme on the acceleration
$\ddot{X}$ and on the velocity $\dot{X}$, furthermore let us
apply a second order accurate backwards finite difference scheme for
its absolute value $|\dot{X}|$.

This yields
\begin{equation}\label{eq_soliddiscr}
\frac{X^{n+1} - 2X^n + X^{n-1}}{(\Delta t)^2} =
-C^n\frac{\frac{X^{n+1}-X^{n-1}}{2\Delta
    t}}{\left|\frac{3X^n-4X^{n-1}+X^{n-2}}{2\Delta
      t}\right|+\delta} + \overline{C}^n,
\end{equation}
so by rearranging the terms we get that
\begin{equation}
\begin{aligned}
\bigg(1+(\Delta t)^2&\frac{C^n}{|3X^n-4X^{n-1}+X^{n-2}|+2\Delta
    t\delta}\bigg) X^{n+1} \\
&= 2X^n - \left(1-(\Delta t)^2\frac{C^n}{|3X^n-4X^{n-1}+X^{n-2}|+2\Delta
    t\delta}\right)X^{n-1}+(\Delta t)^2\overline{C}^n.
\end{aligned}
\end{equation}

Notice that by the definition of $C$, the factor of $X^{n+1}$ is
strictly positive, thus we may multiply by its inverse to obtain an
explicit formula for $X^{n+1}$, namely
\begin{equation}
X^{n+1} = \overline{G}(\zeta^n,X^n,X^{n-1},X^{n-2}),
\end{equation}
where
\begin{equation*}
\overline{G}(\zeta^n,X^n,X^{n-1},X^{n-2}) = \frac{2X^n - \left(1-(\Delta t)^2\dfrac{C^n}{|3X^n-4X^{n-1}+X^{n-2}|+2\Delta
    t\delta}\right)X^{n-1}+(\Delta t)^2\overline{C}^n}{1+(\Delta t)^2\dfrac{C^n}{|3X^n-4X^{n-1}+X^{n-2}|+2\Delta
    t\delta}}.
\end{equation*}

\begin{rem}
Notice that we choose the same
order for the central finite difference schemes for the first and second order
derivatives. The accuracy of the backwards finite difference scheme
for the absolute value of the velocity estimate was chosen
accordingly, and may be adapted.
\end{rem}

The important property of this discretization is that the discrete
equation verifies a certain a priori estimate, just like the one
presented in Corollary \ref{cor_velocityestimate}.
\begin{lem}\label{lem_discretedissip}
Owing to  equation (\ref{eq_soliddiscr}), the velocity associated to
the displacement 
\begin{equation}
\dot{X}^n := \frac{X^{n+1}-X^{n-1}}{2\Delta t}
\end{equation}
verifies the following
\begin{equation}
  |\dot{X}^n|\leqslant \overline{C}^nn\Delta t + \mathcal{O}((\Delta t)^2).
\end{equation}
\end{lem}
\prove Let us multiply equation (\ref{eq_soliddiscr}) by
$\dot{X}^n$. Then, we obtain
\begin{equation*}
\frac{(X^{n+1}-X^{n})^2-(X^n-X^{n-1})^2}{2(\Delta t)^3} = -C^n\frac{(\dot{X}^n)^2}{\left|\frac{3X^n-4X^{n-1}+X^{n-2}}{2\Delta
      t}\right|+\delta} + \overline{C}^n\dot{X}^n.
\end{equation*}
Notice first of all that the first term on the right hand side is
non-positive. Moreover, on the left hand side we have a second order
approximation of the value $\partial_t((\dot{X})^2)(n\Delta t)$,
as such we may replace it with another second order approximation,
namely $\partial_t((\dot{X}^n)^2)$ modulo an error term of order
$\mathcal{O}((\Delta t)^2)$. Combining these two remarks, we get that
\begin{equation}
\partial_t((\dot{X}^n)^2)(t)^2\leqslant \overline{C}^n\dot{X}^n+\mathcal{O}((\Delta t)^2).
\end{equation}
A Gr\"onwall type inequality allows us to conclude in the desired
way. \quod

\subsection{The wave tank and its boundaries}\label{sec_tank}

Notice that the weakly nonlinear Boussinesq system
(\ref{eq_boussinesq}) has its domain of validity for $x\in\mathbb{R}$
which is clearly not the case for numerical models. Hence we
consider the discretized model in a wave tank of sufficiently large
size (it shall be detailed for each experiment in the next
section). The idea is to place the solid in the middle of the tank and
numerically generate the waves (soliton or wave train) near the solid,
therefore allowing us to focus on the middle section of the tank, the
effects of boundary conditions imposed on the horizontal limits would
be negligible.

In all the cases, the width of the wave tank is taken to be at least
$100L$ with $L$ being the wavelength of a wave, and the analysis is
focused on the central $20L\sim 30L$ wide region of the domain.

As for the boundary conditions, for the sake of clarity, solid walls
are implemented, implying reflective boundaries for the fluid
variables. This means Neumann-type boundary condition for the ``scalar''
parameters $\zeta$ and $h$ (the derivative equals to $0$) and
Dirichlet-type boundary condition on the fluid velocity $\overline{V}$
(equaling to $0$).

Variables for the solid are not concerned by these limiting conditions
since they are independent of the spatial variable. The only
technical effect is that the simulation is to be stopped when the
object touches the boundary. As evoked before, due to the size of the
wave tank, it is not a likely scenario.

\subsection{Further remarks}\label{sec_remarks}

The first remark concerns a reference solution for the Boussinesq
system with completely flat bottom topography. It is known that a
solitary wave solution (soliton) exists for this equation (see for
example \cite{stevansoliton}). Naturally, by introducing a solid
object on the bottom of the fluid domain, this referential solution
will not stay a solitary wave propagating at a constant speed, it will
be transformed, deformed according to the governing equations and the
change in bottom topography. Nevertheless it serves as a basic tool to
analyze the effects of the object on a single wave, as it will be done
during the first half of the next section.

Some words should be mentioned about the effect of the solid
displacement $X$ on the fluid equations in
(\ref{eq_systembous1}). Due to the horizontal motion of the
object, it is present in the variable $b$ as a translation. Since
there are multiple instances when derivatives of the bottom topography
function are taken, the order of the operations has to be established
with respect to the discretization.

In our algorithm, at each time step, the actual bottom surface shall
be calculated via the translation by $X^n$ of the initial state
$\mathfrak{b}$ and then it is discretized on the grid points and the
mid-points as well. If we were to discretize the initial bottom
topography and then translate it, it is clear that the fitting of the
translated discrete bottom to the actual grid would create additional
error terms which could potentially decrease the overall accuracy.

However we shall mention that our approach works mainly because the initial
bottom surface $\mathfrak{b}$ as a function defined on $\mathbb{R}$ is
known, so that we can discretize any translated instance of it. It
would not be possible without accurate interpolations, if the solid height was initially given only on
grid points.

An important remark concerns the integral terms present in the solid
equation (\ref{eq_newtonorder1}). Except for the volume of the solid,
these integrals involve integration over the support of the object (at
its current state) of fluid variables that are initially only defined
on mesh points (or mid-points). As such the
applicable accurate methods are somewhat limited. In our situation we
chose a third order Simpson method (resulting in a global error of
order $4$) and it writes as follows
\begin{equation}
\int_{j\Delta x}^{k\Delta x}\zeta(x)\,dx\approx\Delta
x\frac{1}{6}\left(\zeta_j+\zeta_k+2\sum_{l=j+1}^{k-1}\zeta_l+4\sum_{l=j}^{k-1}\zeta_{l+1/2}\right).
\end{equation}

\section{Numerical results}

In this part we present a multitude of numerical experiments and
simulation results for different one dimensional wave propagation and transformation
scenarios. Effects of a solid allowed to move will be compared to a
fixed solid case to highlight the main features of this new approach.

As explained in Remark \ref{rem_coefficient}, the appropriate base
water depth for our case is (at least) $10$ meters (due to the
difference in the order of magnitude of the coefficients). As such,
comparable physical experiments are not exactly available. 

For our
study in general, the wave tank is taken with the following physical
parameters: its width is exactly $1000$ (meters), its height is $H_0=20$. The
attributed shallowness parameter is then determined by the choice of
the wavelength $L$, with the two main cases being $L=40$ for
$\mu=0.25$, and $L=20\sqrt{10}$ for $\mu=0.1$. This implies that the
wave-tank is $25$ wavelengths long in the former case and approximately
$16$ wavelengths long in the latter case. The principal observational
area is in the vertical section $[400, 600]$ of the wave tank.

The solid will be considered to be given by a truncated Gaussian function, the discrete truncation determined by an error term $\epsilon = 10^{-4}$, that is
\begin{equation*}
\mathfrak{b}(x) = b_0(x)\mathbf{1}(b_0>\epsilon),\quad\textnormal{ where }b_0(x) = a_{bott}\exp\left(-10\left(\frac{x}{L}\right)^2\right).
\end{equation*}

As for its physical parameters, the default choice for the mass
parameter is chosen to be $\tilde{M}=\frac{2}{3}$, corresponding to an
approximate solid density of $\varrho_S=2\frac{g}{cm^3}$. The vertical
size as well as the coefficient of the friction will be varied during
most of the simulations.

\subsection{Order of the numerical scheme}\label{sec_verif}

The first set of tests concerns the verification of the convergence of
the scheme. In \cite{linman} the overall algorithm for
the fluid equations was
second order accurate both in time and in space. Since we made some
adaptations on the original scheme, it is reasonable to verify how
this has affected the order of the convergence. It turns out that
considering all situations and parameter regimes relevant to the model
we always have a
convergence of order at least $\mathcal{O}((\Delta x)^3+(\Delta
t)^2)$.

\begin{figure}[t]
\centering
\begin{subfigure}{.5\textwidth}
  \centering
  \includegraphics[width=\linewidth]{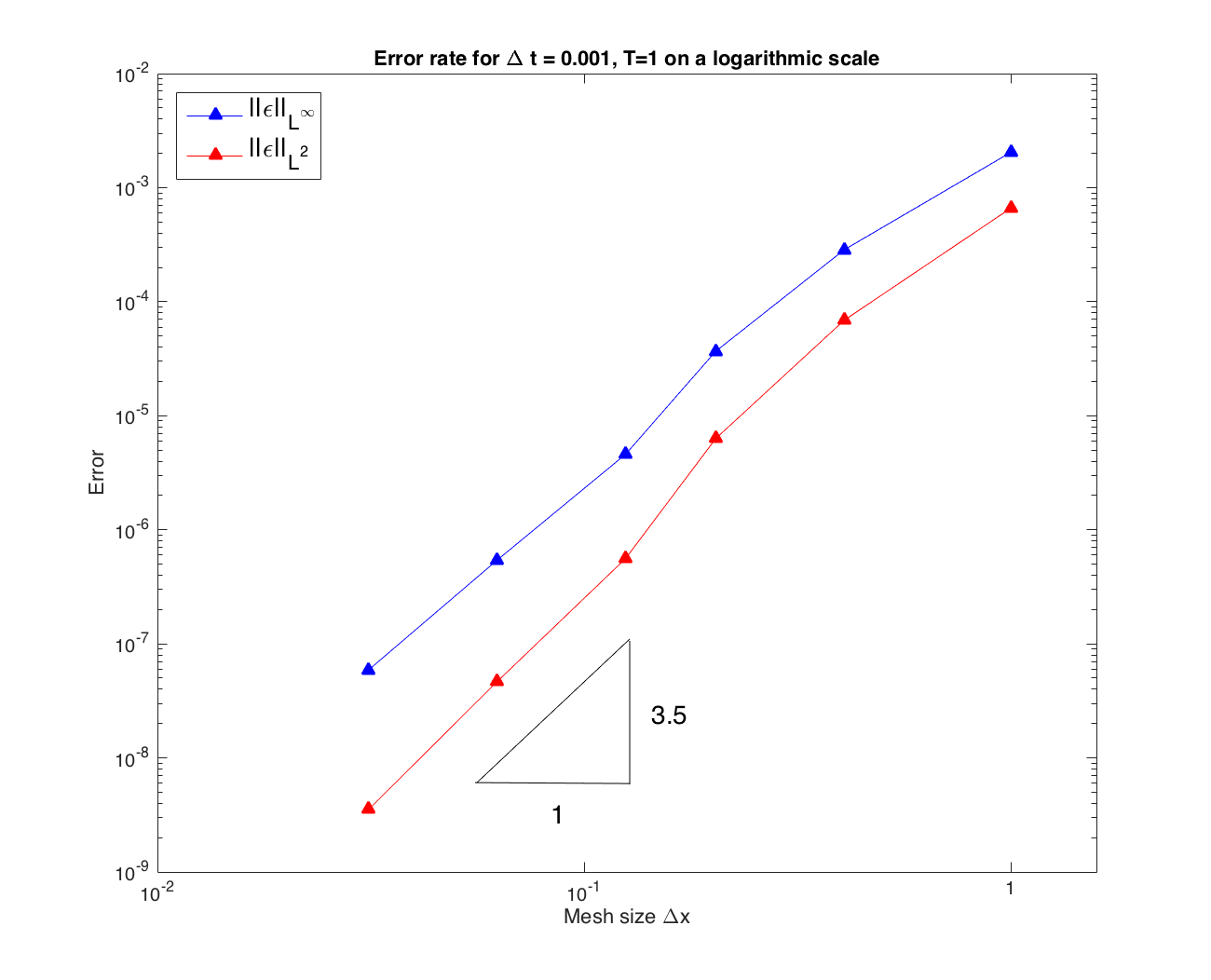}
  \caption{$\mu=\varepsilon=0.1$ solitary wave evolution, spatial error}
  \label{fig_spatialerr}
\end{subfigure}%
\begin{subfigure}{.5\textwidth}
  \centering
  \includegraphics[width=\linewidth]{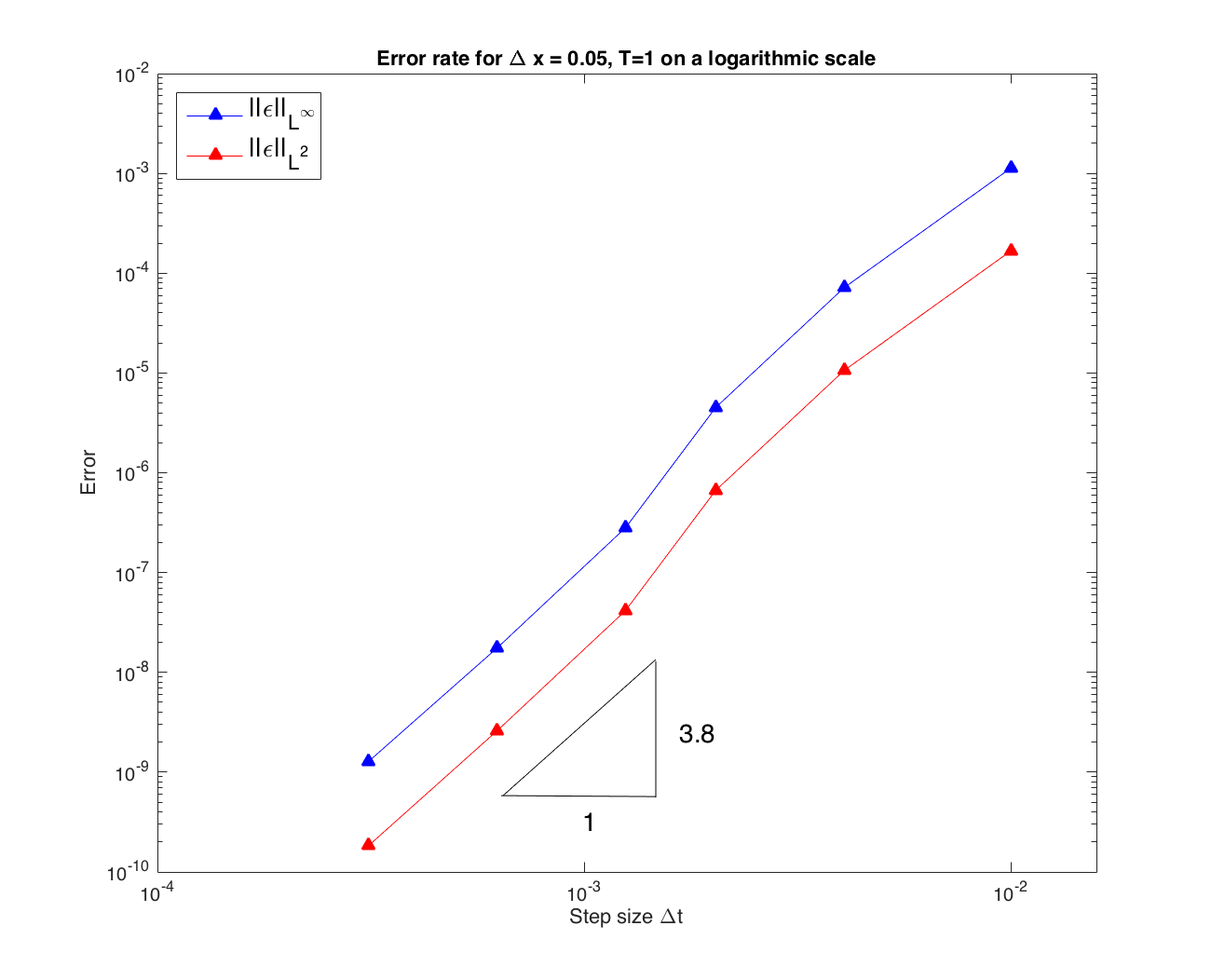}
  \caption{$\mu=\varepsilon=0.1$ solitary wave evolution, temporal error}
  \label{fig_temporalerr}
\end{subfigure}
\caption{Discretization error for solitary wave evolution over flat bottom}
\label{fig_simpleerr}
\end{figure}

\subsubsection{Convergence of the scheme over a flat bottom topography}

First of all, we consider the wave tank without the presence of the solid
and its effects. This means that we are considering the simplified
algorithm for a flat bottom case. The convergence is then checked numerically by the exact traveling wave solution for (\ref{eq_boussinesq}) with $b\equiv 0$. The existence of a solitary wave solution is a well known fact, for explicit computations, please refer to \cite{chensoliton}, for a more general approach, one may see \cite{stevansoliton} for example. By searching the solution for (\ref{eq_boussinesq}) as a solitary traveling wave with constant speed $c>1$ ($\zeta = \zeta_c=\zeta_c(x-ct)$, $\overline{V}=\overline{V}_c=\overline{V}_c(x-ct)$) it is easy to check that the velocity profile has to satisfy
\begin{equation}
-(\overline{V}_c)'' = -\frac{3}{\mu c}\overline{V}_c\left(c-\frac{1}{c-\varepsilon\overline{V}_c}-\frac{\varepsilon}{2}\overline{V}_c\right),
\end{equation}
with $\overline{V}_c<c/\varepsilon$. From the velocity profile, the surface elevation then can be recovered simply by
\begin{equation}
\zeta_c=\frac{\overline{V}_c}{c-\varepsilon\overline{V}_c}.
\end{equation}

With this at our disposal, the solitary wave solution serving as a reference can be reconstructed by a fourth order Runge-Kutta method (in order to avoid any influence on the accuracy of the overall scheme).

Two basic verification tests have been designed, each one testing for
the wave-tank and wave parameters $\mu=\varepsilon=0.1$. The simulation is run for a time $T=1$.

The first test measures the error of the spatial discretization,
compared to the solitary wave solution, by
fixing the time step $\Delta t$ sufficiently small ($\Delta t=0.001$) and varying the spatial discretization's step
size, giving us $\Delta x=2^{-k}$ for $k\in\{0,1,2,3,4,5\}$. We
measure the $L^2$ and $L^\infty$ norms of the error for the surface elevation
$\zeta$. Figure \ref{fig_spatialerr} shows the results, on a
logarithmic scale. The scheme is behaving as an algorithm of order $3.5$
in the spatial discretization.

The second test measures the error of the time discretization,
compared to the solitary wave solution, by
fixing the grid size $\Delta x$ sufficiently small ($\Delta x=0.05$) and varying the time discretization's step size,
giving us $\Delta t=0.01\cdot2^{-k}$ for
$k\in\{0,1,2,3,4,5\}$. Once again we measure the $L^2$ and $L^\infty$ norms of the
error for the surface elevation $\zeta$. Figure \ref{fig_temporalerr}
shows the results, on a logarithmic scale. The scheme is behaving as
an algorithm of order $3.8$ for the time discretization, an almost
order $4$ convergence which would be the ideal scenario for the
applied Adams predictor-corrector method.

\begin{figure}
\centering
\begin{subfigure}{.5\textwidth}
  \centering
  \includegraphics[width=\linewidth]{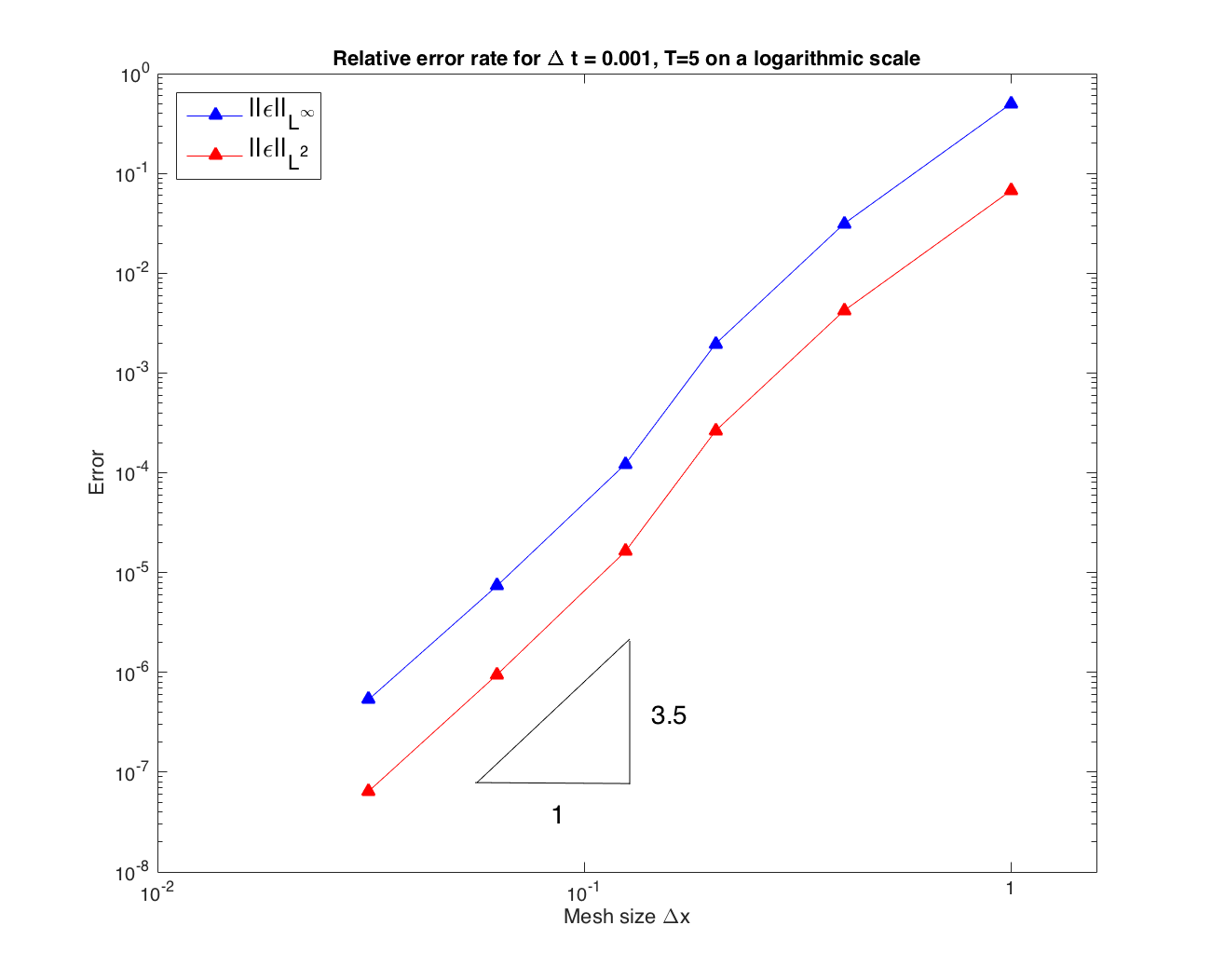}
  \caption{$\mu=\varepsilon=0.1$ solitary wave evolution, spatial error}
  \label{fig_smallerrspace}
\end{subfigure}%
\begin{subfigure}{.5\textwidth}
  \centering
  \includegraphics[width=\linewidth]{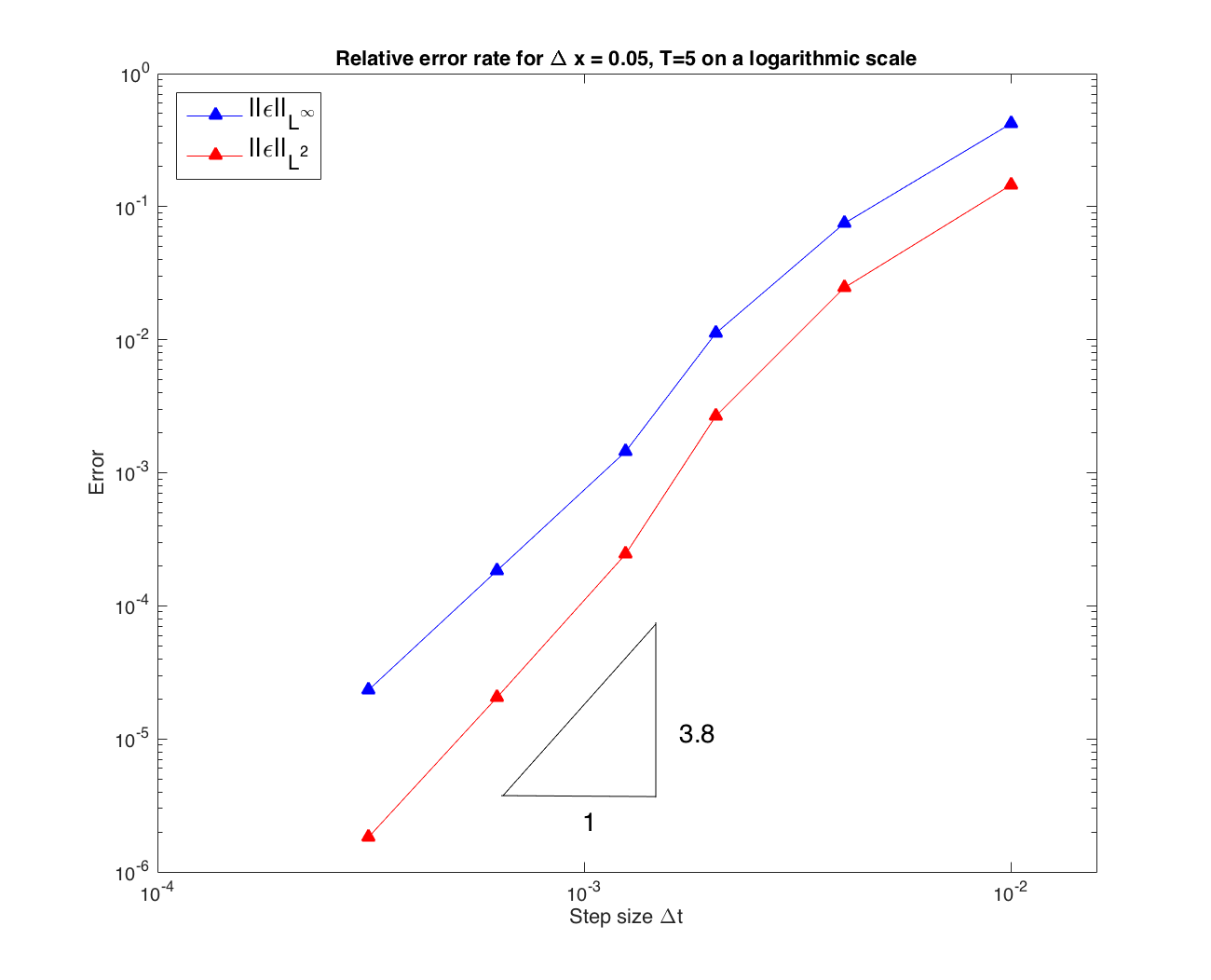}
  \caption{$\mu=\varepsilon=0.1$ solitary wave evolution, temporal error}
  \label{fig_smallerrtime}
\end{subfigure}
\caption{Discretization error for solitary wave evolution over non-flat
bottom ($\beta = 0.3$, $c_{fric}=0.01$)}
\label{fig_smallerr}
\end{figure}

Therefore we have established a major improvement over \cite{linman},
by an accurate interpolation, and more coherent accuracy in the finite
difference schemes, one can indeed obtain an almost fourth order convergence.

As elaborated in Section 3.4 of \cite{linman}, a von Neumann stability
analysis can be carried out for the linearized system. Since no
changes have been made in the time discretization of the fluid
equations (system (\ref{eq_boussinesq})), their analysis can be
adapted in a straightforward way to our
case too, leading to a CFL condition of the form:
\begin{equation}
\sqrt{gH_0}\frac{\Delta t}{\Delta x}\leqslant 0.5.
\end{equation}

In order to respect this stability condition, in
what follows we set $\Delta x=0.05$ and $\Delta t = 0.001$.

\subsubsection{Convergence of the scheme for a non-flat bottom topography}

Lacking an explicit solution for the non-flat bottom case, we
performed a relative error analysis to test the global convergence of
the full coupled system as well, meaning that as a reference solution
we calculated the surface elevation for $\Delta x = 0.01$, $\Delta t =
10^{-4}$ and we compared it to the calculated surface elevations for
less refined mesh sizes and time steps. Two physically different testing
parameters were chosen, the first one corresponding to an
immobile solid at the bottom, the second one representing the case of
the fully coupled problem where solid movement is observed. For the
first case, the physical parameters of the system
were chosen to be $\mu = \varepsilon =0.1$, $\beta = 0.3$, with a frictional
coefficient of $c_{fric}=0.01$. The simulation is run for a time $T=5$.

Figure \ref{fig_smallerr} shows the relative error of the surface
elevation compared to our choice of reference solution (in $L^2$ and
$L^\infty$ norms), carried out for
either a fixed grid size ($\Delta x=0.05$) or for a fixed time step
($\Delta t = 0.001$), dividing by $2$ the other step parameter for each
consecutive measurement, just as before. One can observe that the overall spatial
discretization stays of order $3.5$, and the temporal error stays in
the previously observed order of $3.8$ too.

\begin{figure}
\centering
\begin{subfigure}{.5\textwidth}
  \centering
  \includegraphics[width=\linewidth]{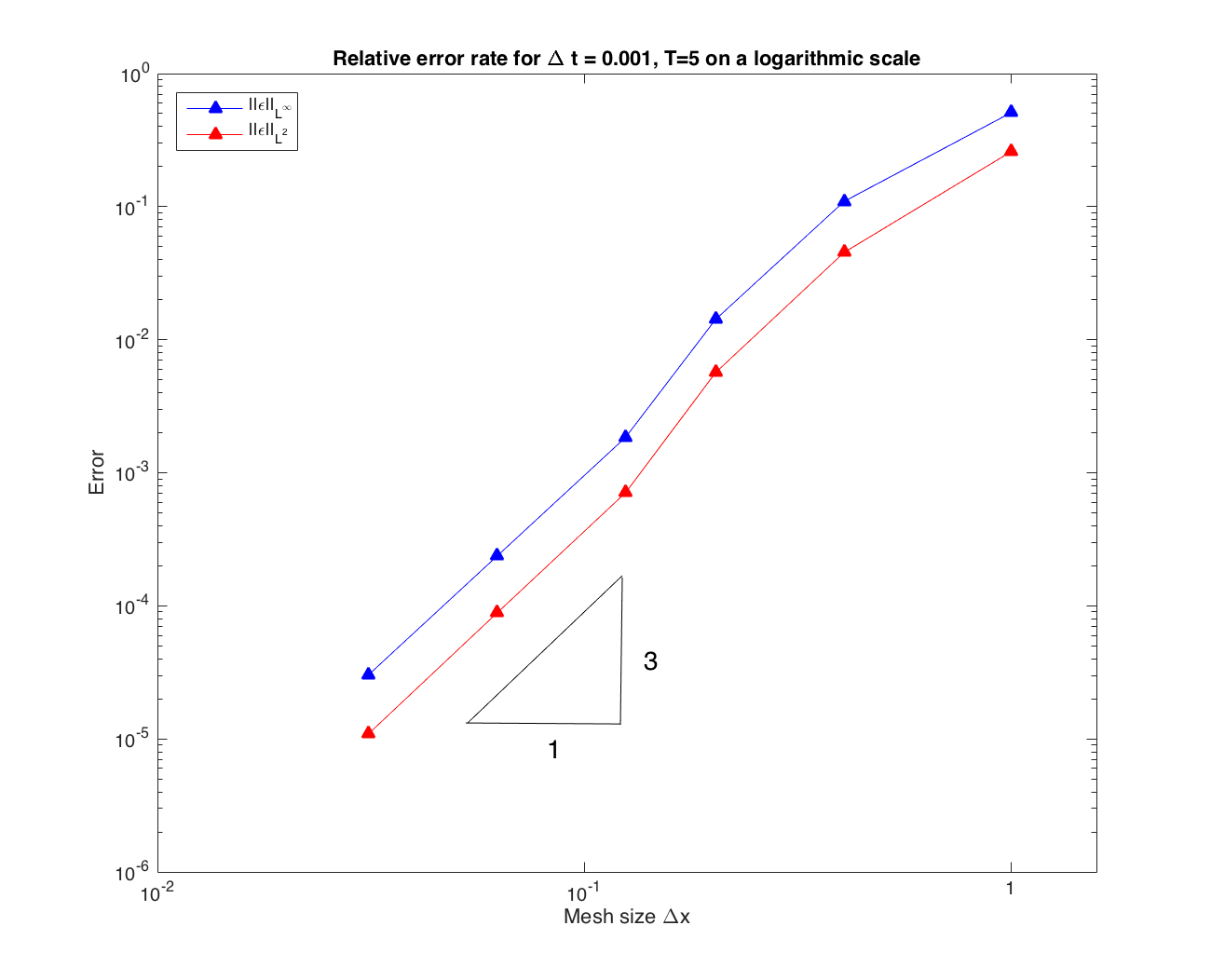}
  \caption{Relative spatial error for $\Delta t=0.001$}
  \label{fig_coupleerr1}
\end{subfigure}%
\begin{subfigure}{.5\textwidth}
  \centering
  \includegraphics[width=\linewidth]{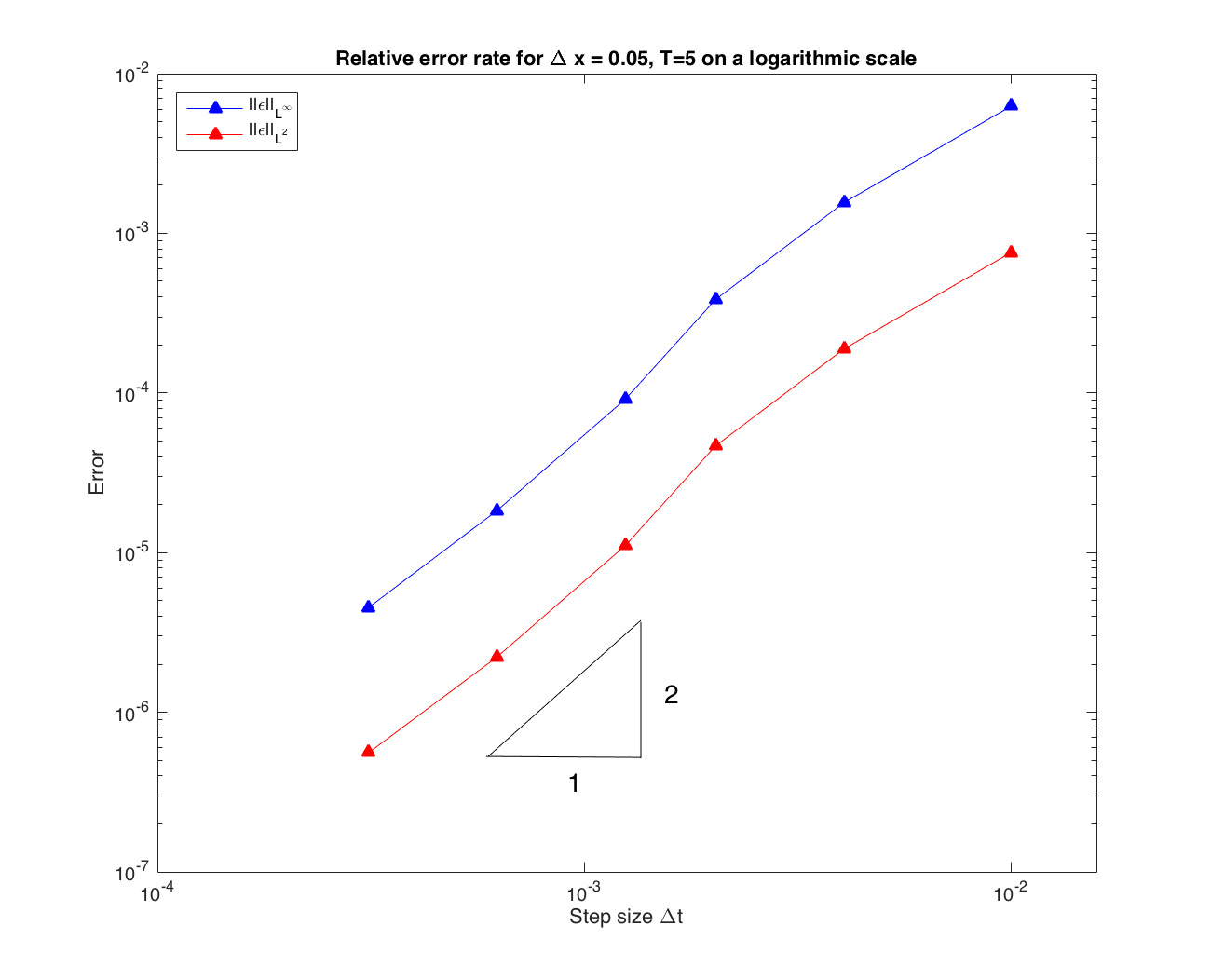}
  \caption{Relative temporal error for $\Delta x=0.05$}
  \label{fig_coupleerr2}
\end{subfigure}
\caption[Relative error measurements for the coupled problem]{$L^2$
  and $L^\infty$ error for a wave evolution over time $T=5$ ($\mu =
  \varepsilon =0.2$, $\beta = 0.4$, $c_{fric}=0.001$)}
\label{fig_coupleerr}
\end{figure}

For the
second case, the physical parameters of the system
were chosen to be $\mu = \varepsilon =0.2$, $\beta = 0.4$, with a frictional
coefficient of $c_{fric}=0.001$. The simulation is run until a time
$T=5$ allowing for sufficiently long interaction between the solid and
the incoming wave.

Figure \ref{fig_coupleerr} shows the relative errors (in $L^2$ and
$L^\infty$ norms) for the surface elevation, carried out once again for
either a fixed grid size ($\Delta x=0.05$) or for a fixed time step
($\Delta t = 0.001$). One can observe that the overall spatial
discretization has decreased to an order $3$, attributed to the loss
in accuracy represented by the observable solid motion. Notice also
that we have a temporal convergence rate of $2$, attributed to the fact
that the time discretization scheme for the solid was chosen to be
only of order $2$.

\subsection{Evolution of the surface of a single passing wave}\label{sec_singli}

 In this part we present the two main characteristic situations that have
 been observed during the ensemble of the simulations, each with a
 passing wave and a breaking wave case. Two
 representative examples were chosen, showing step by step the
 ``transformation'' of a single passing wave over the solid, first for
 a highly frictional case, then for the almost perfect sliding case. As
 initial conditions, the approaching wave is taken to be the solitary
 wave solution for the flat bottom case.

\subsubsection{The regime of a large coefficient of friction}

The first example presents what most of the test cases looked liked in
the numerical experiments to follow; the wave passing over the solid,
getting slightly perturbed by the bottom topography irregularity (the
presence of the solid) and continuing its trajectory with a modified
amplitude and an altered form. With a shallowness parameter of $\mu=0.1$, the
initial wave amplitude is taken as $4$ (meaning $\varepsilon =
0.2$). The solid has a maximal vertical size of $6$ and is subjected
to a frictional sliding on the bottom, with a coefficient of
$c_{fric}=0.5$.

\begin{figure}
\centering
\begin{subfigure}{.5\textwidth}
  \centering
  \includegraphics[width=0.92\linewidth]{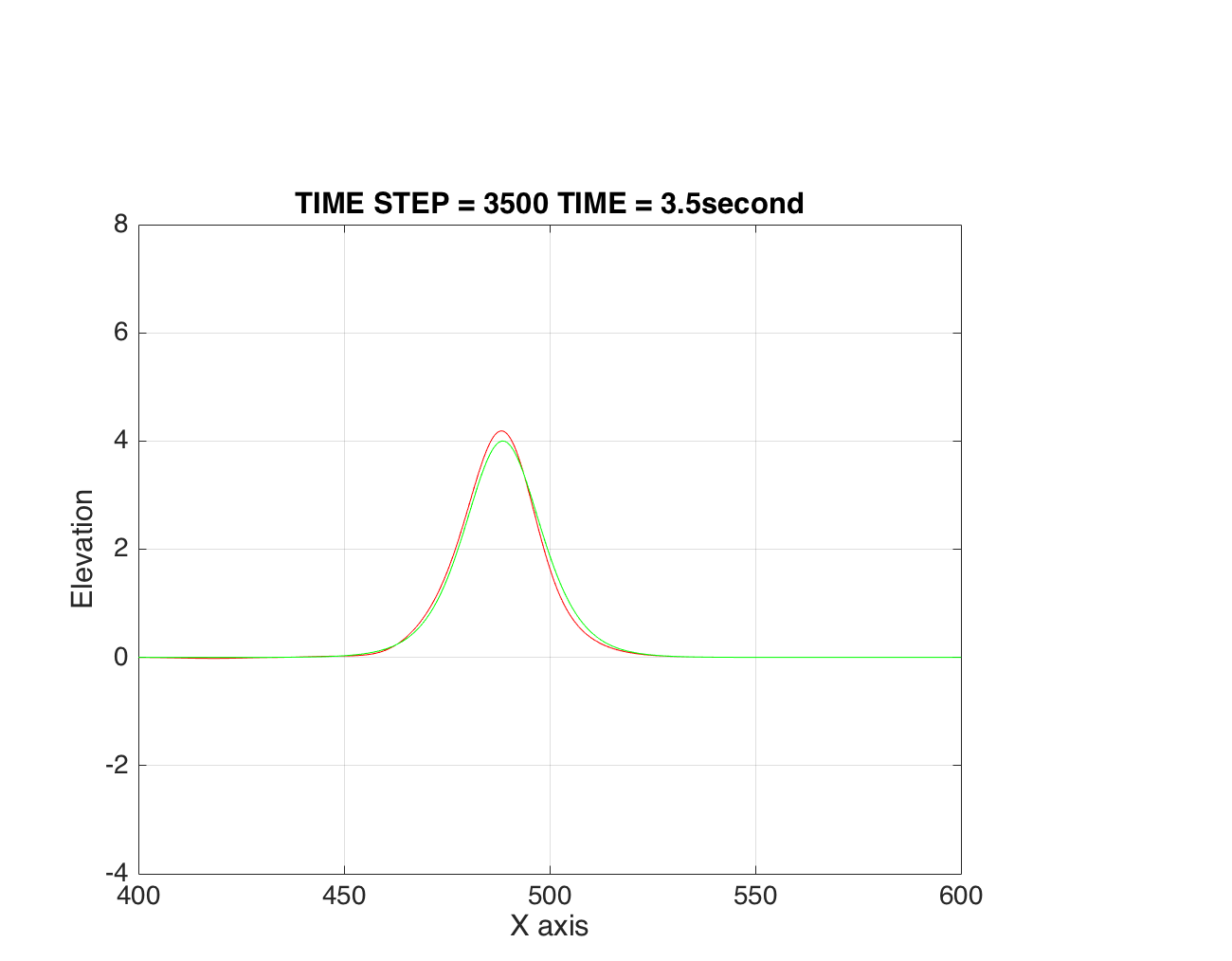}
\end{subfigure}%\hspace{-2em}
\begin{subfigure}{.5\textwidth}
  \centering
  \includegraphics[width=0.92\linewidth]{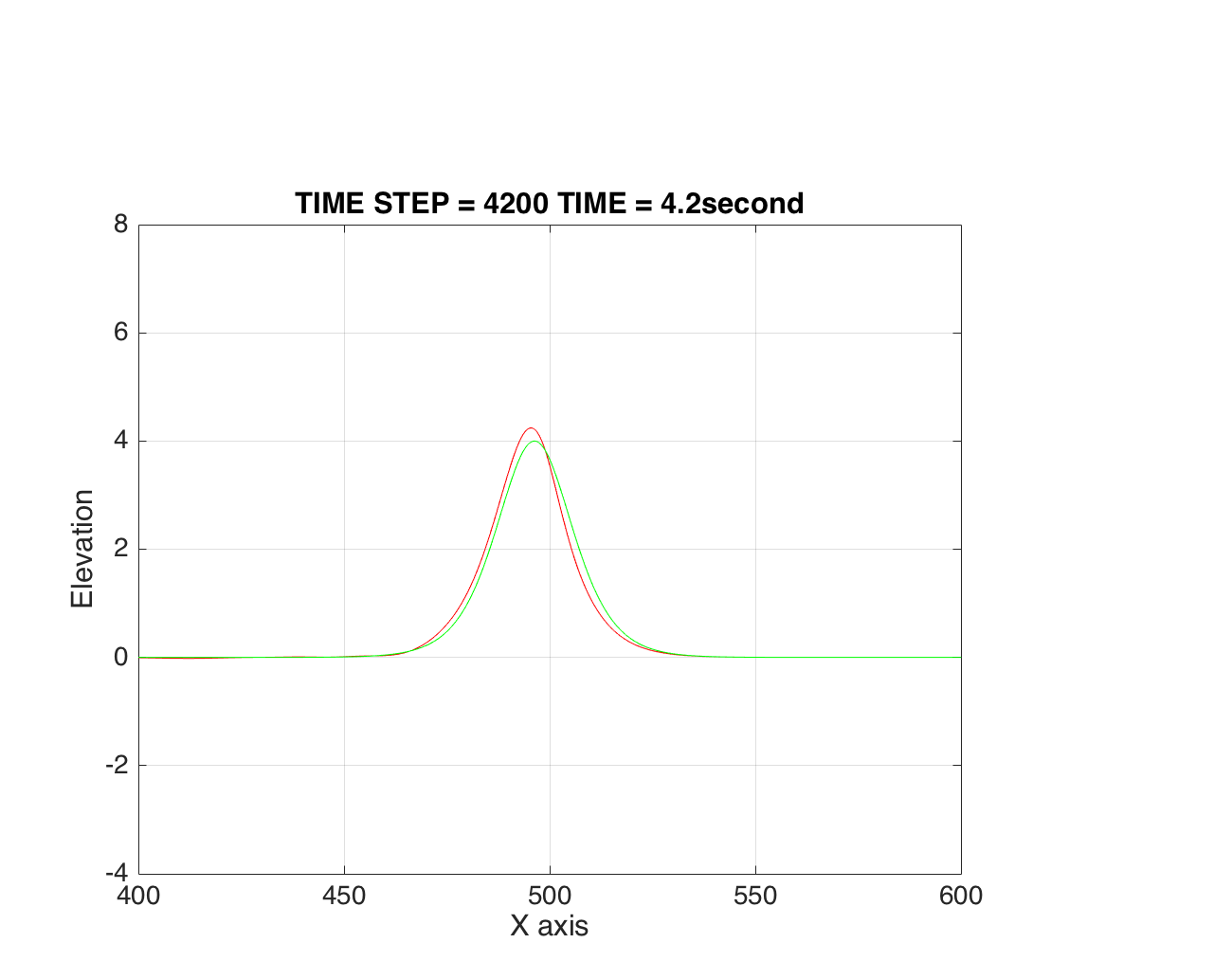}
\end{subfigure}
\newline
\begin{subfigure}{.5\textwidth}
  \centering
  \includegraphics[width=0.92\linewidth]{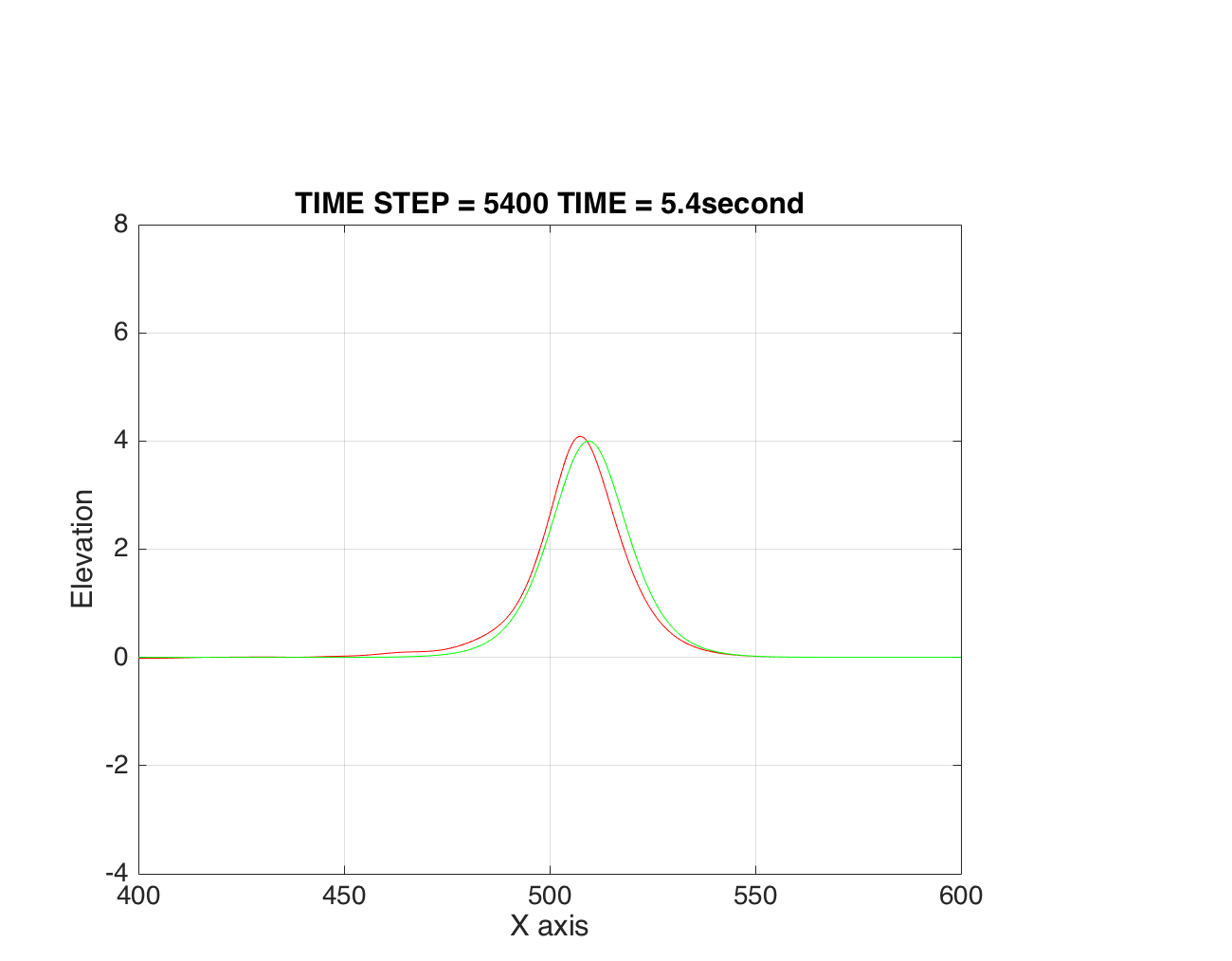}
\end{subfigure}%
\begin{subfigure}{.5\textwidth}
  \centering
  \includegraphics[width=0.92\linewidth]{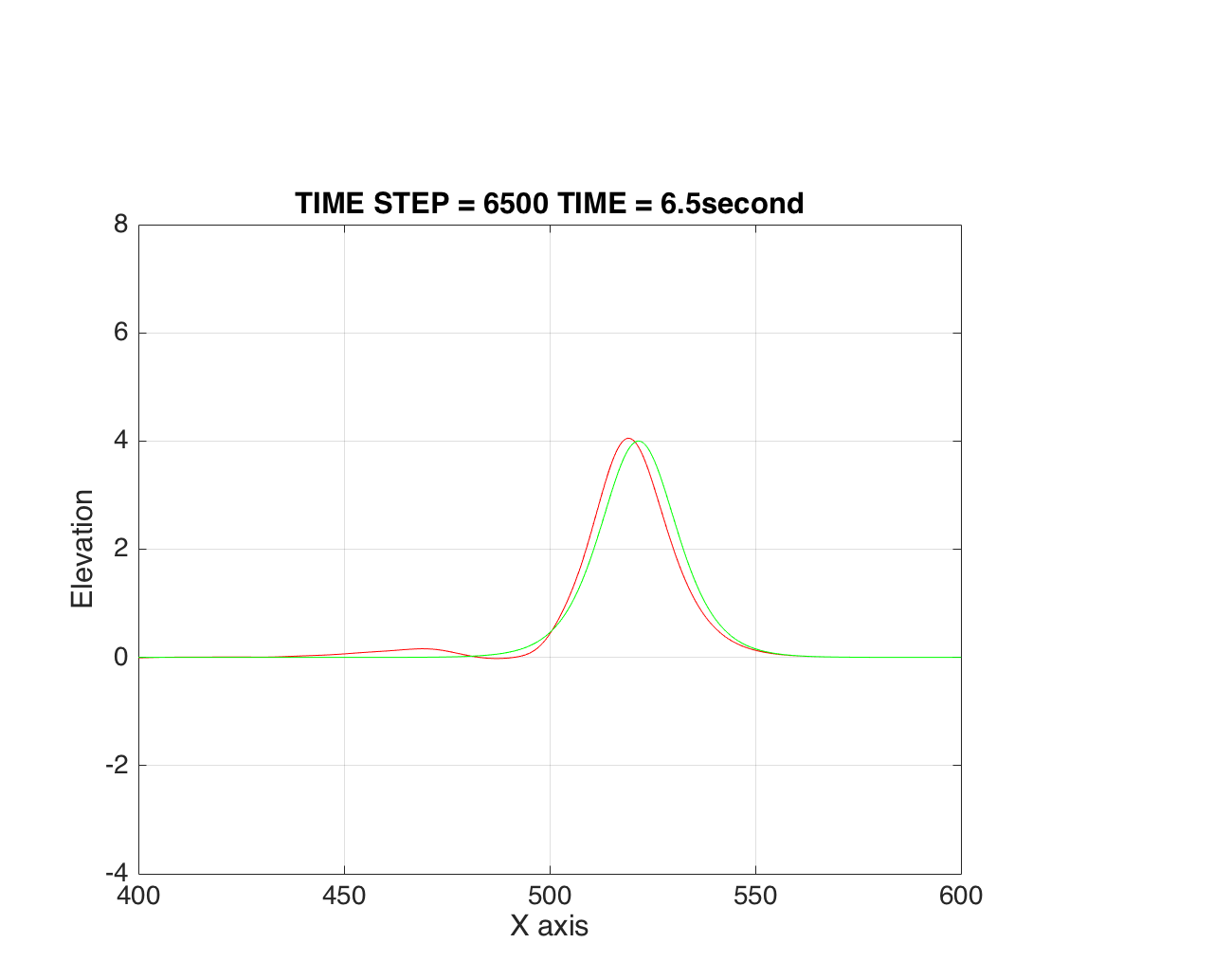}
\end{subfigure}
\newline
\begin{subfigure}{.5\textwidth}
  \centering
  \includegraphics[width=0.92\linewidth]{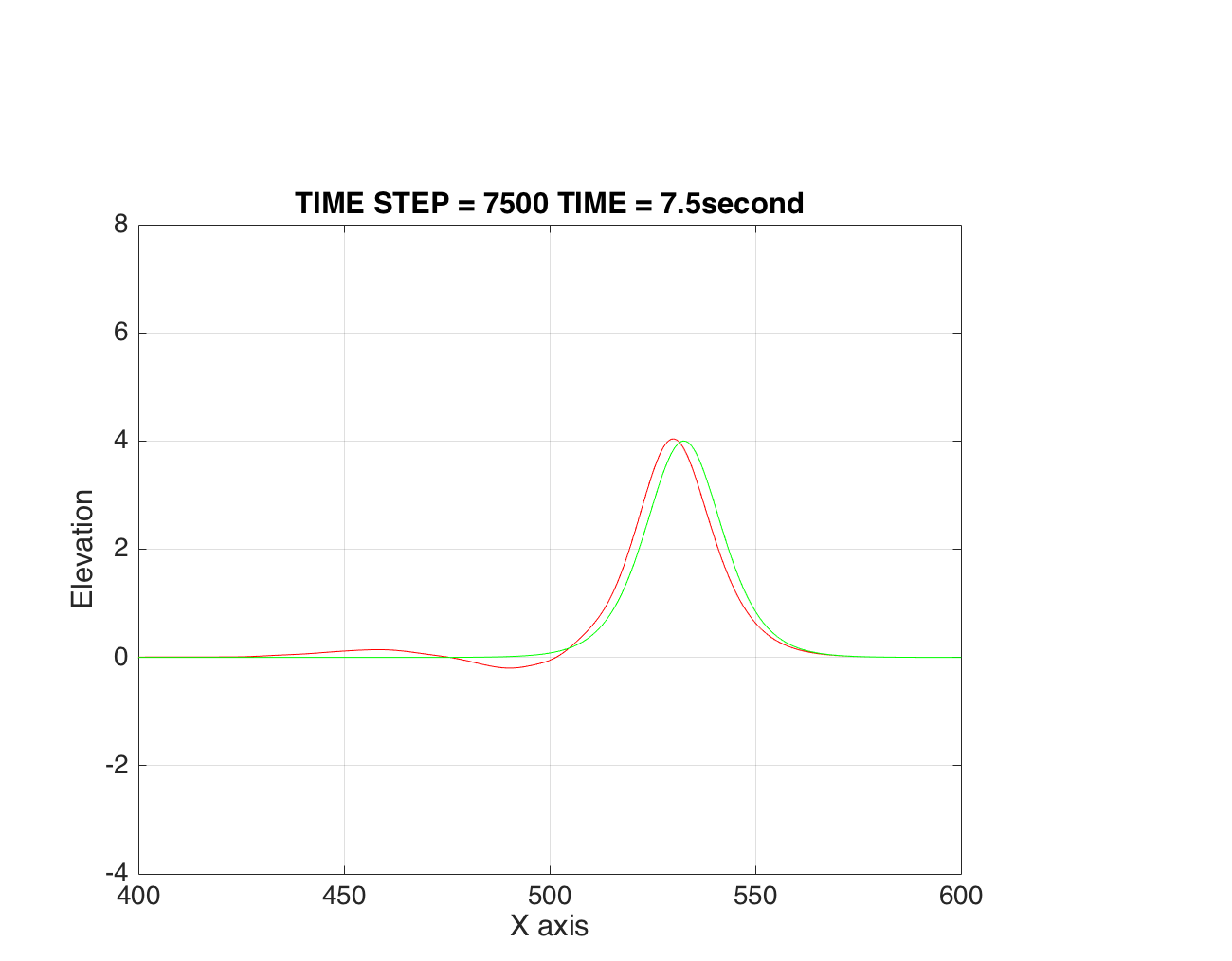}
\end{subfigure}%\hspace{-2em}
\begin{subfigure}{.5\textwidth}
  \centering
  \includegraphics[width=0.92\linewidth]{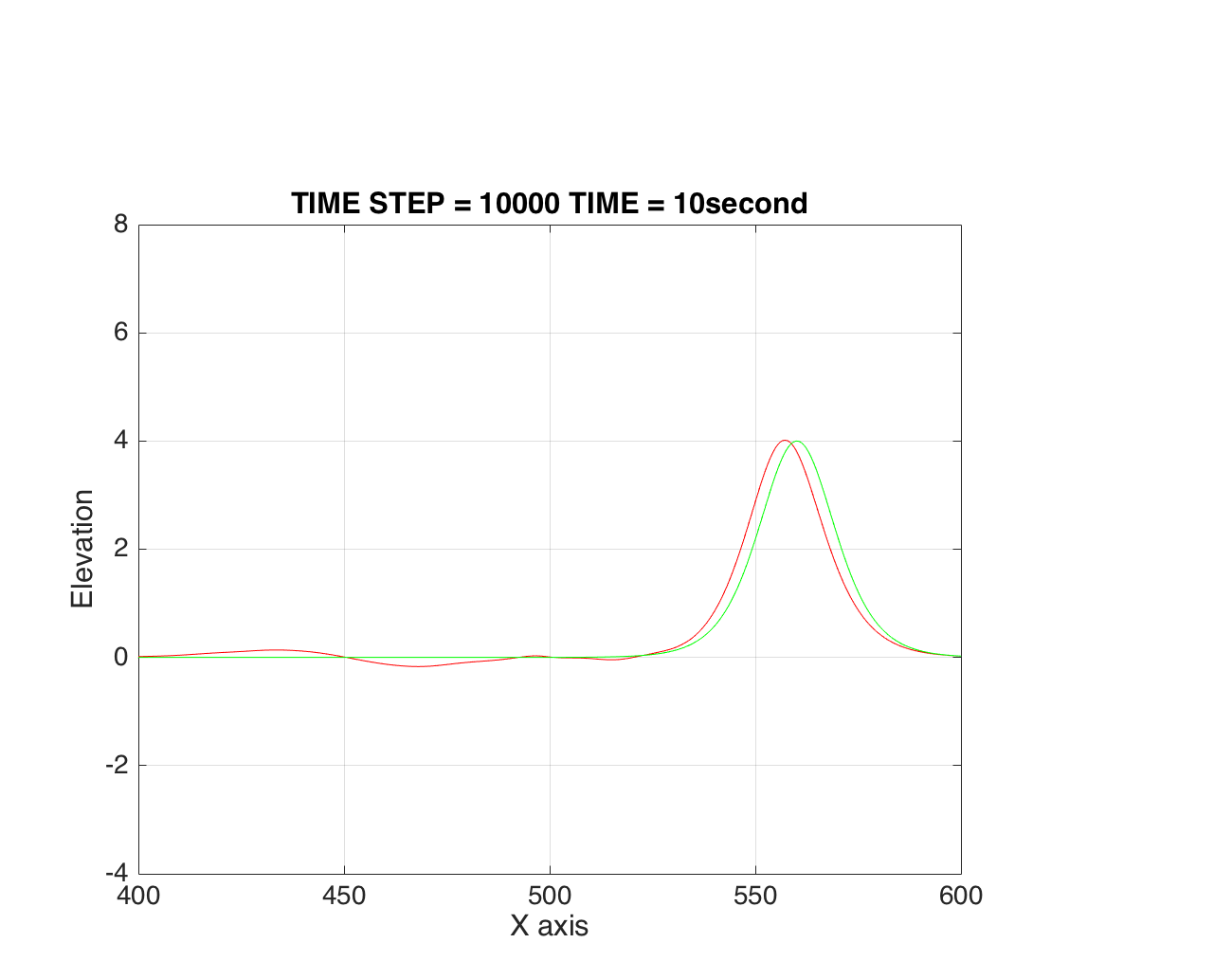}
\end{subfigure}
\caption[Evolution of a passing wave]{Evolution of a passing wave
  ($\mu=0.1$, $\varepsilon=0.2$) over
  a small obstacle ($\beta=0.3$, $c_{fric}=0.5$); {\color{red}red}
    curve is the passing wave,
    {\color{green}green} curve is the
    reference soliton for flat bottom}
\label{fig_wavepass}
\end{figure}

In Figure \ref{fig_wavepass} this passing wave is plotted ({\color{red}red})
at different time steps. As a reference, the flat bottom solitary wave
is also visualized ({\color{green}green}), propagating at a constant
speed, and allowing for a better qualitative comparison for the
changes the approaching wave undergoes. The solid is centered around the
horizontal coordinate $x=500$ with a numerical support spanning through
the interval $[480,520]$, for a better visibility, it has been omitted
  from the figure.

In the first figures the wave approaches the solid, and thus wave
shoaling is observed (its amplitude increases) until its peak at around time step $4200$,
after which the wave crest passes over the peak of the solid and
drops (step $5400$), due to the drop in the bottom
topography. After this drop the initial wave, now slightly asymmetric,
continues onward, with an amplitude (slightly) less than its initial
value. 

Moreover the wave became out
of phase due to the ``bump'' in its motion. One can also observe a
backwards going small amplitude long wavelength wave created by the
drop after passing over the solid (starting from the back trough at around
step $6500$). It is important to remark that
due to the high frictional term, essentially no solid displacement
was observed.

\begin{figure}
\centering
\begin{subfigure}{.4\textwidth}
  \centering
  \includegraphics[width=\linewidth]{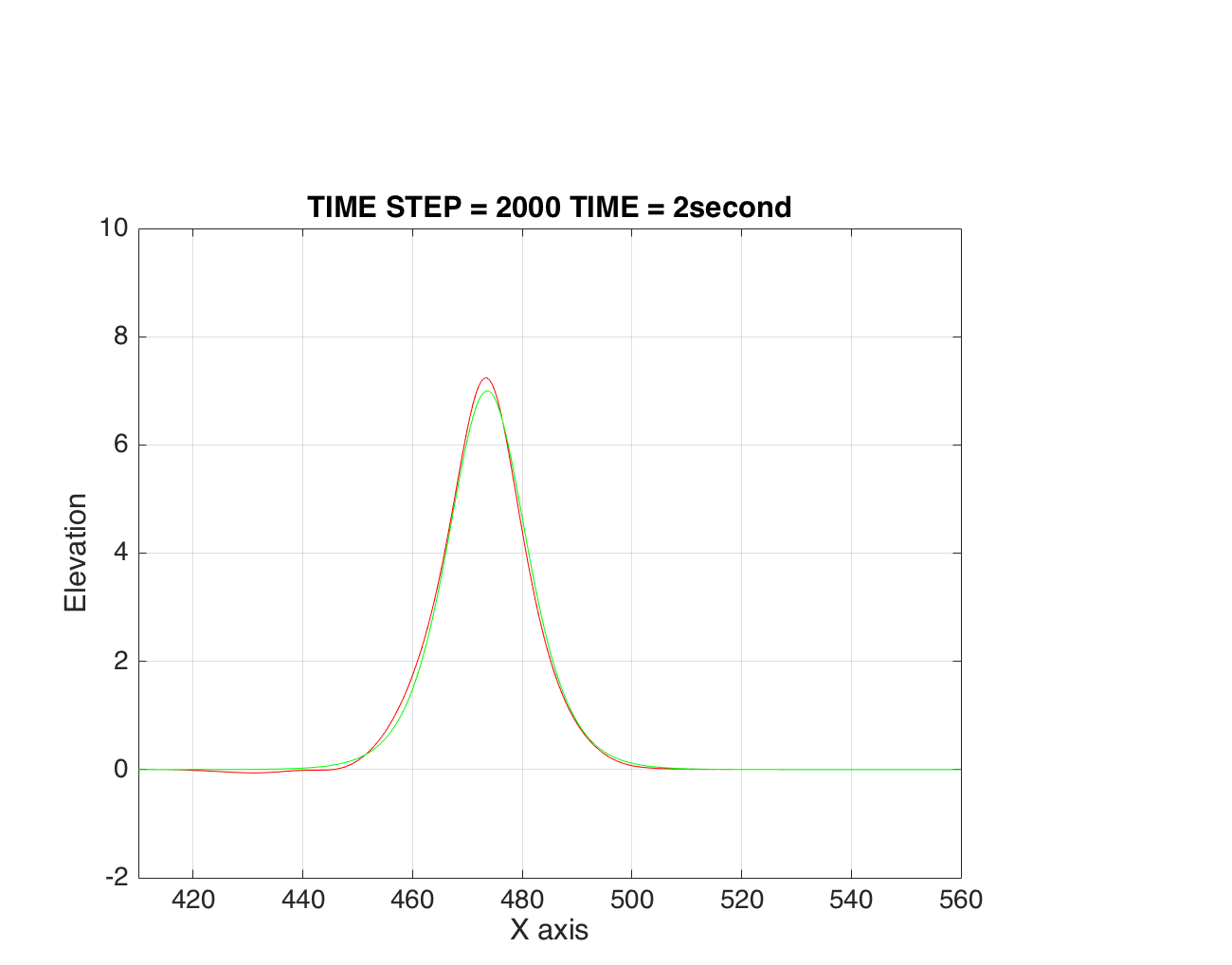}
\end{subfigure}%
\begin{subfigure}{.4\textwidth}
  \centering
  \includegraphics[width=\linewidth]{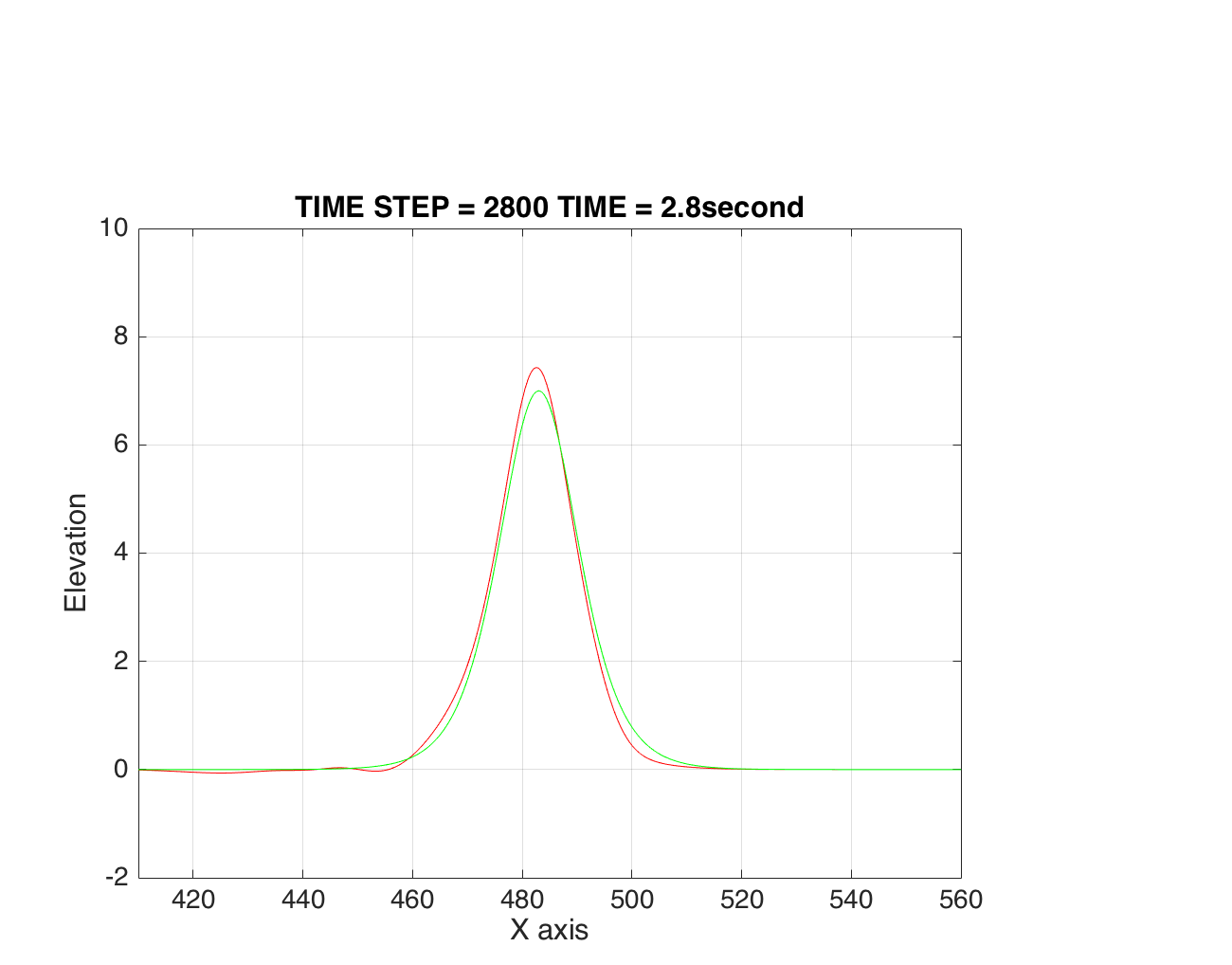}
\end{subfigure}\newline
\begin{subfigure}{.4\textwidth}
  \centering
  \includegraphics[width=\linewidth]{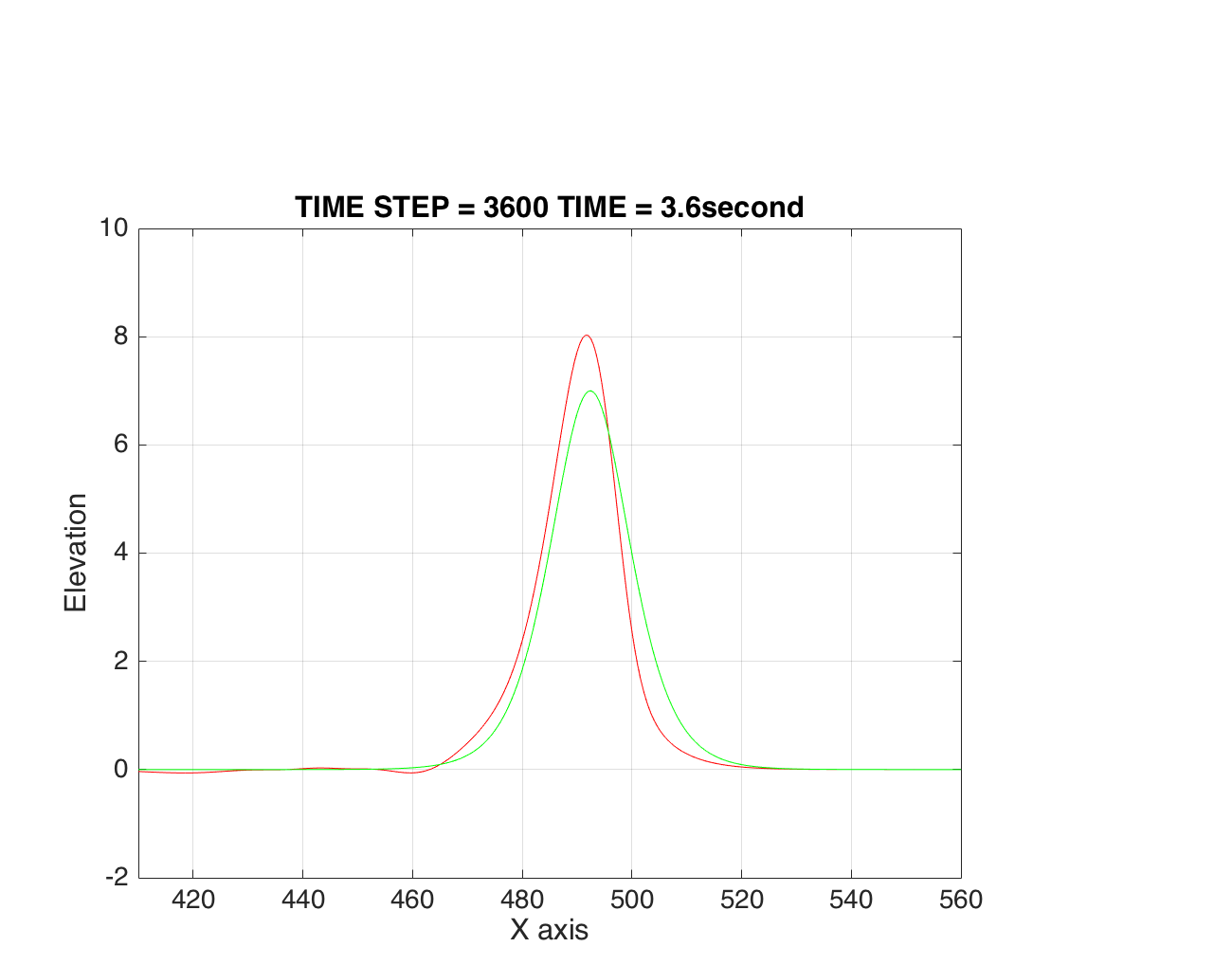}
\end{subfigure}%
\begin{subfigure}{.4\textwidth}
  \centering
  \includegraphics[width=\linewidth]{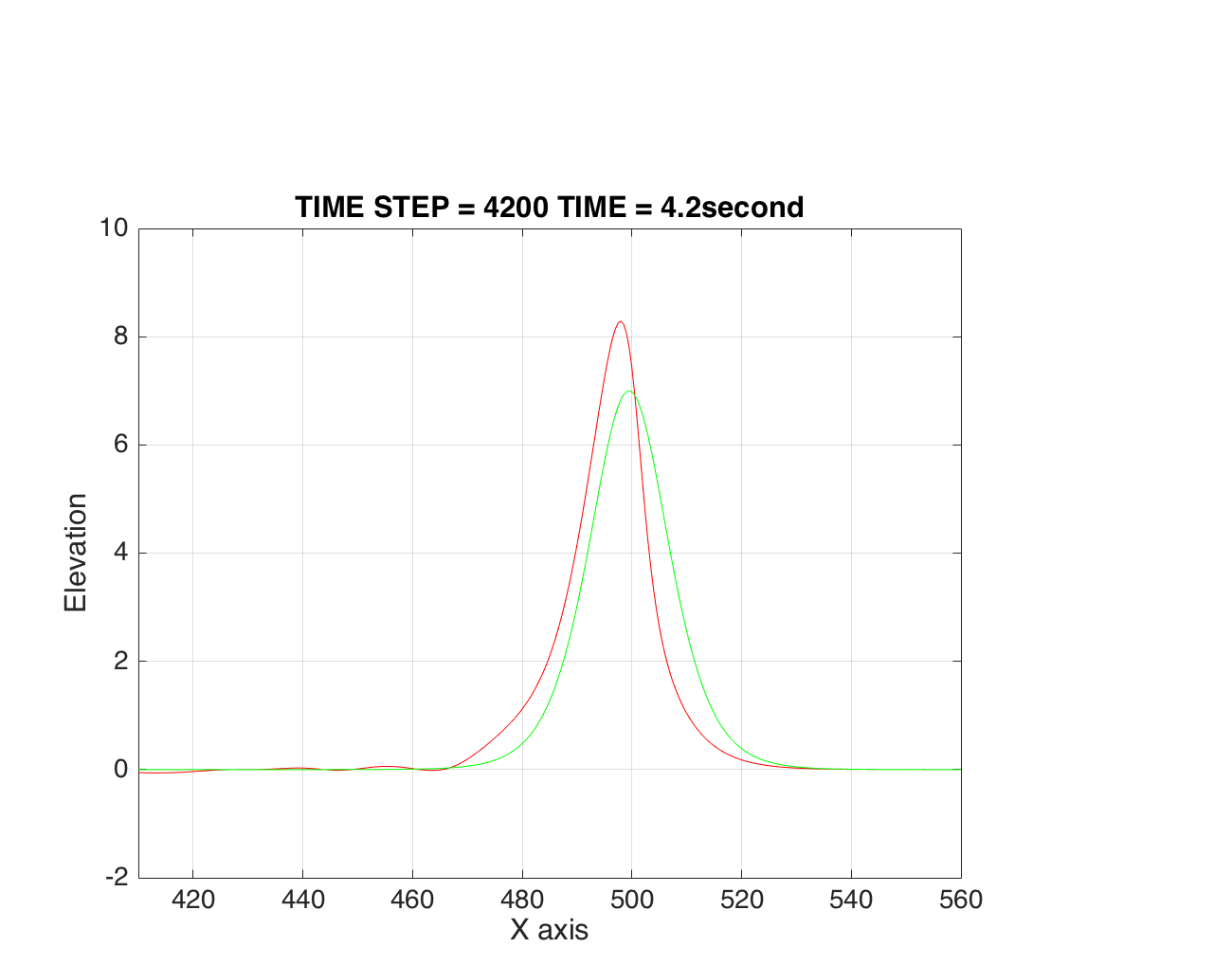}
\end{subfigure}\newline
\caption[Evolution of a breaking wave]{Evolution of an approaching
  large wave ($\mu=0.25$, $\varepsilon=0.35$ and breaking when
  reaching the obstacle ($\beta = 0.5$, $c_{fric}=0.5$); {\color{red}red}
    curve is the approaching wave,
    {\color{green}green} curve is the
    reference soliton for flat bottom}
\label{fig_wavebreak}
\end{figure}

The set of images in Figure \ref{fig_wavebreak} depict the wave breaking
encountered during the simulations for regimes with a large
coefficient of friction. In
this case a numerical
condition was detected in the experiments that signals a possible
wave-breaking due to a steepening wave slope (for more details, see
Section \ref{sec_wavebreaking}).

A plunging (or spilling) type wave-breaking is most
commonly indicated by a critical increase in steepness in the middle
section of the front wave slope, and was observed for relatively large wave
amplitudes with a large object at the bottom. For a representative test case, a
wave of wavelength $40$ is chosen over a base water depth of $20$, with
initial amplitude $7$. The solid has a maximal height of $10$ and
is sliding on the bottom with a dynamical friction coefficient of
$c_{fric}=0.5$. As represented in Figure \ref{fig_wavebreak}, the wave approaching the
solid increases in amplitude, just like before, however this increase
becomes critical as the wave crest approaches the solid peak. The flat
bottom solitary wave is represented only as a reference.

\subsubsection{The regime of a small coefficient of friction}

As for a second set of examples, we present a characteristic situation
in the regime of small coefficient of friction, that is, when solid displacement is observed. This
was observed for all examined vertical solid dimensions, appearing to
be more pertinent in intermediate to high wave amplitude regimes. As a
test case, we chose $c_{fric}=0.001$ for a solid height of $8$. The
wave amplitude is chosen to be $4$ with a shallowness parameter of the
system equaling $\mu=0.1$.

\begin{figure}
\centering
\begin{subfigure}{.5\textwidth}
  \centering
  \includegraphics[width=0.95\linewidth]{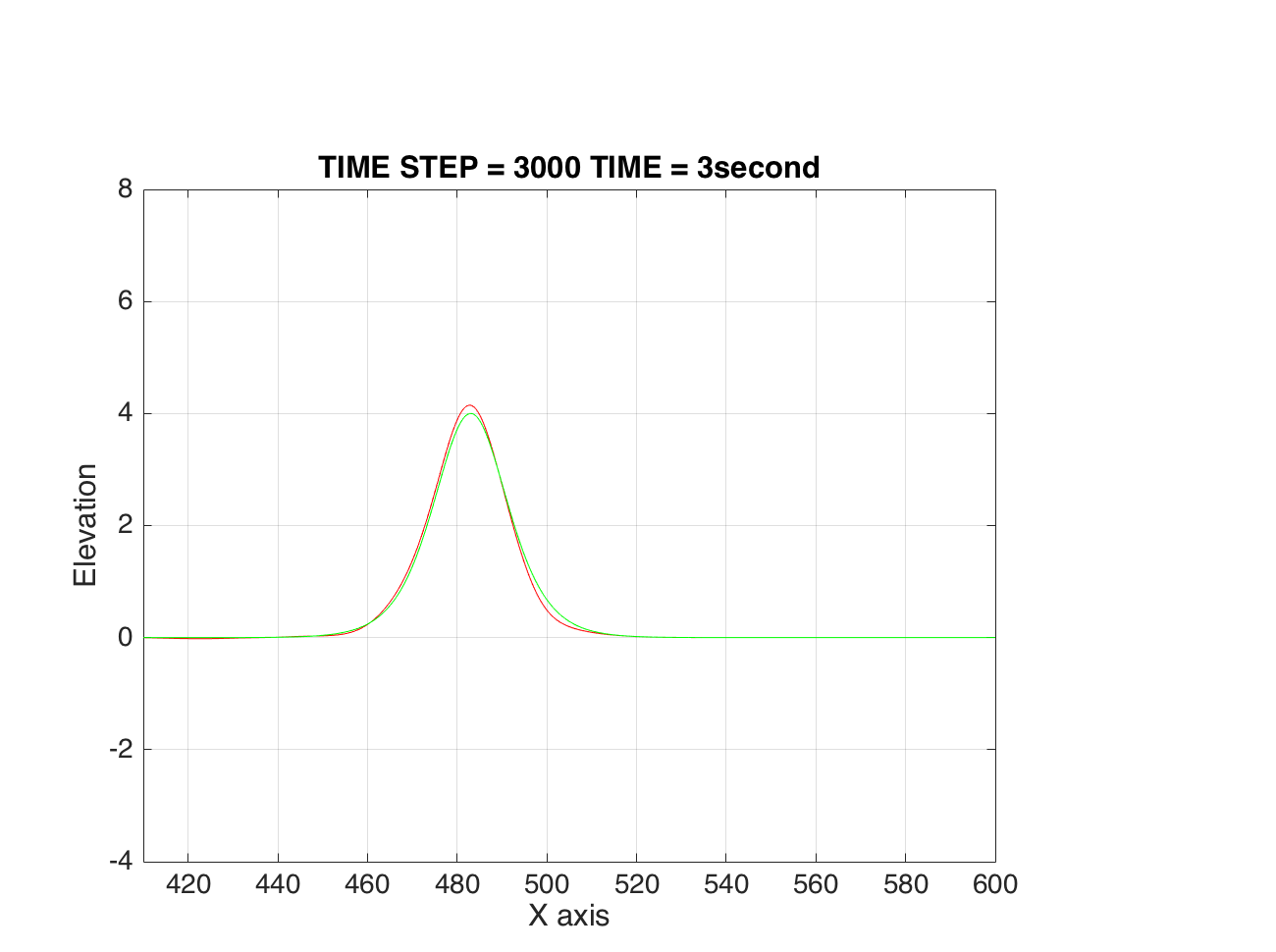}
\end{subfigure}%
\begin{subfigure}{.5\textwidth}
  \centering
  \includegraphics[width=0.95\linewidth]{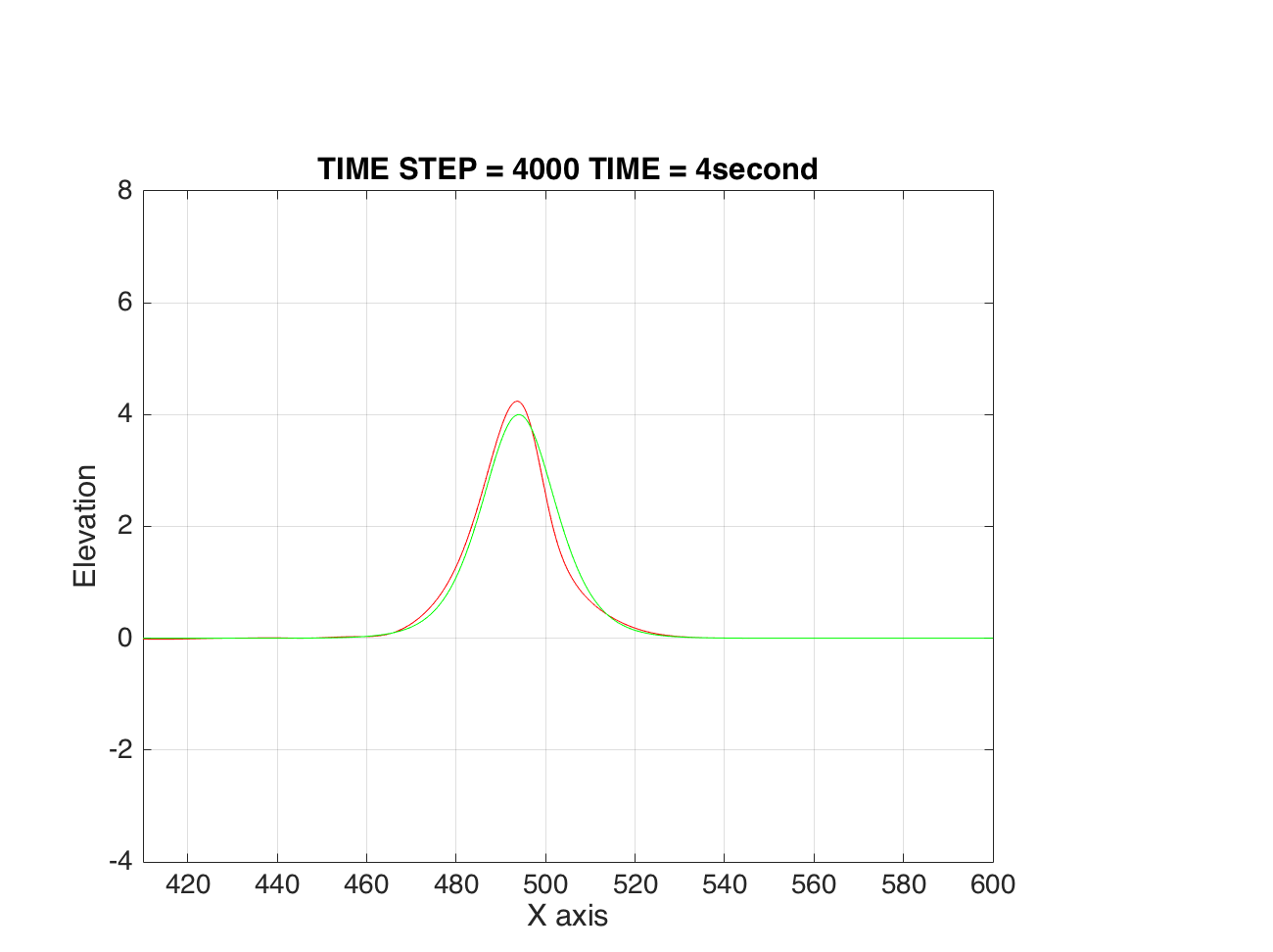}
\end{subfigure}
\newline
\begin{subfigure}{.5\textwidth}
  \centering
  \includegraphics[width=0.95\linewidth]{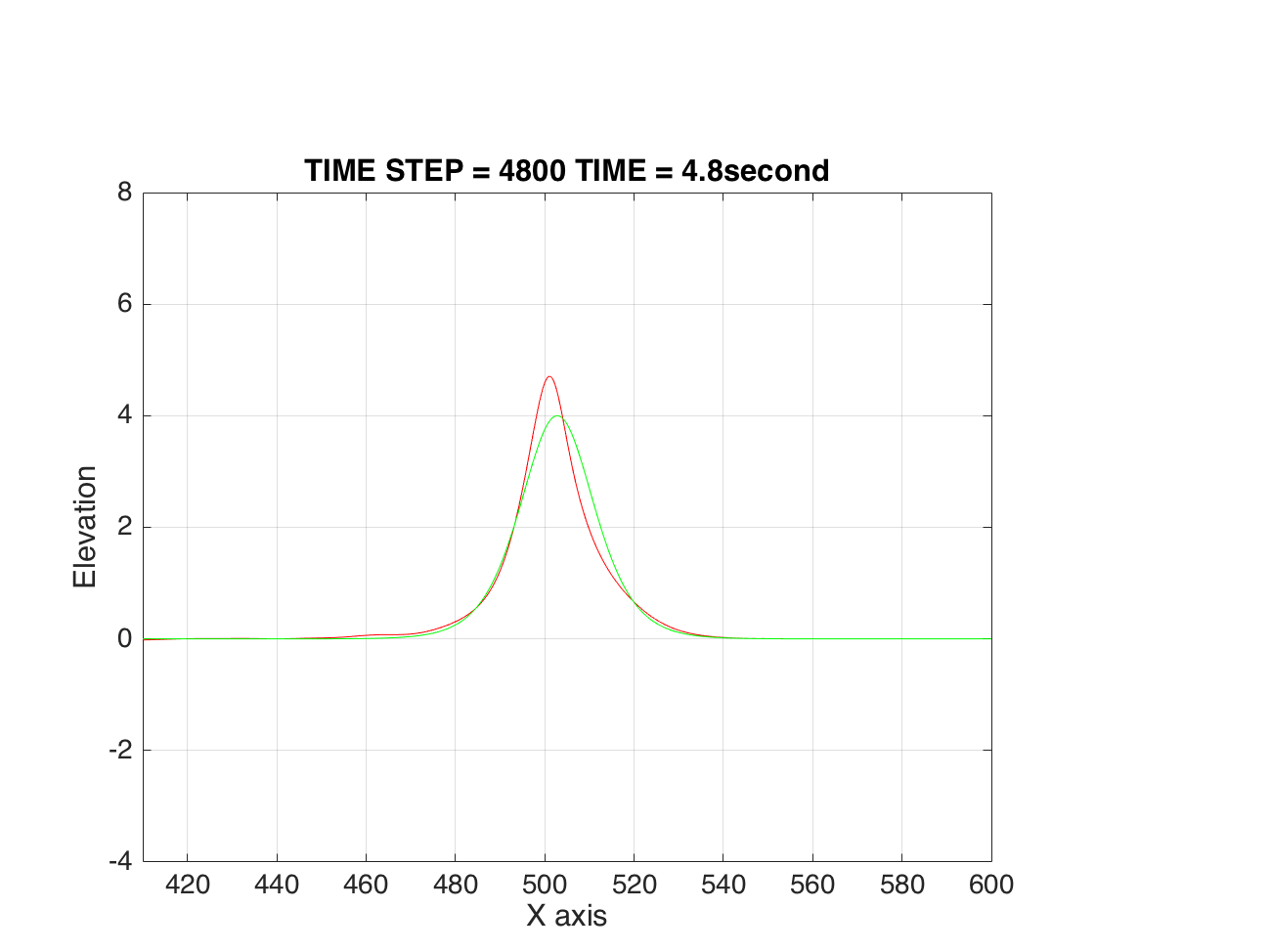}
\end{subfigure}%
\begin{subfigure}{.5\textwidth}
  \centering
  \includegraphics[width=0.95\linewidth]{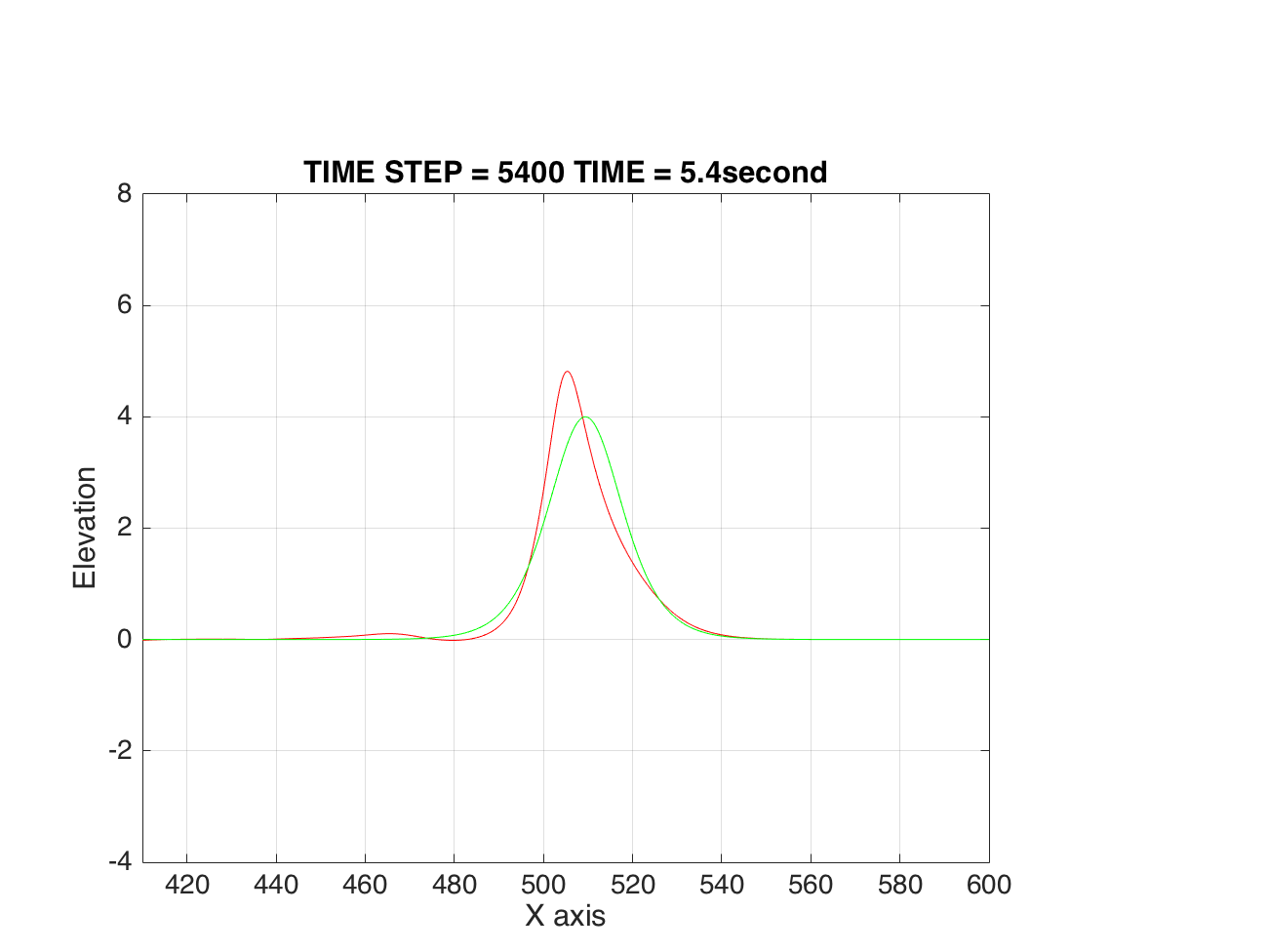}
\end{subfigure}
\newline
\begin{subfigure}{.5\textwidth}
  \centering
  \includegraphics[width=0.95\linewidth]{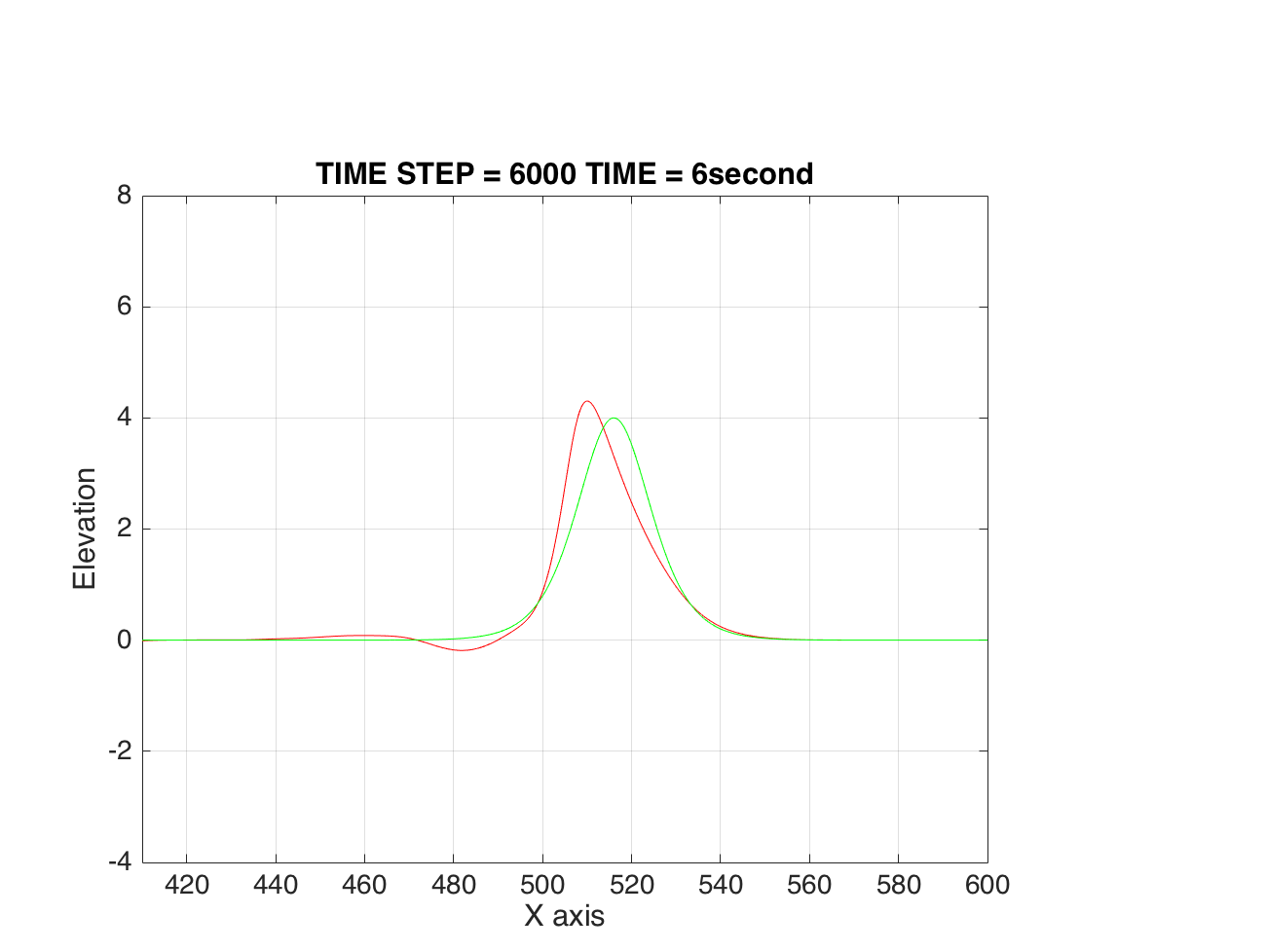}
\end{subfigure}%
\begin{subfigure}{.5\textwidth}
  \centering
  \includegraphics[width=0.95\linewidth]{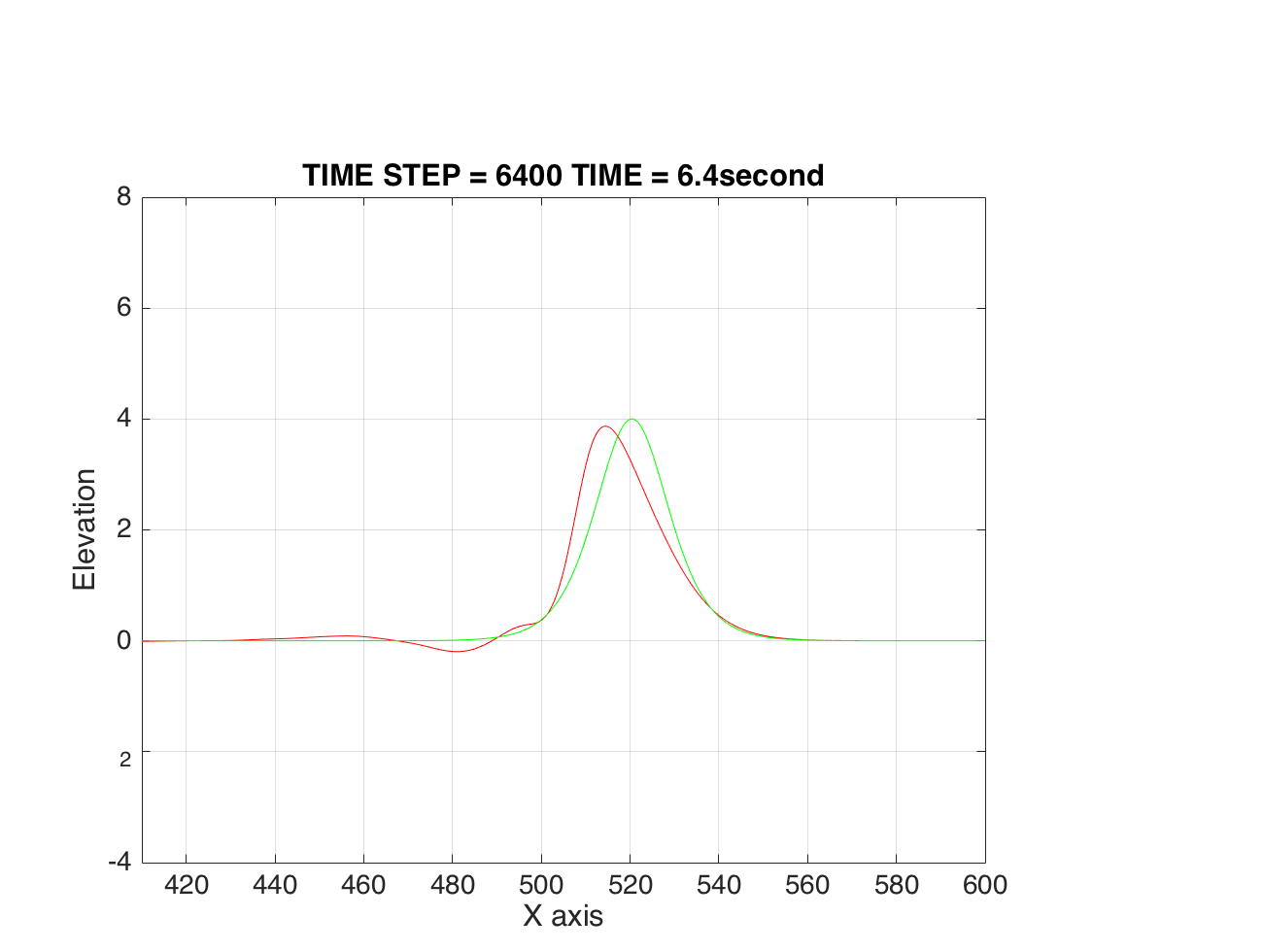}
\end{subfigure}
\newline
\begin{subfigure}{.5\textwidth}
  \centering
  \includegraphics[width=0.95\linewidth]{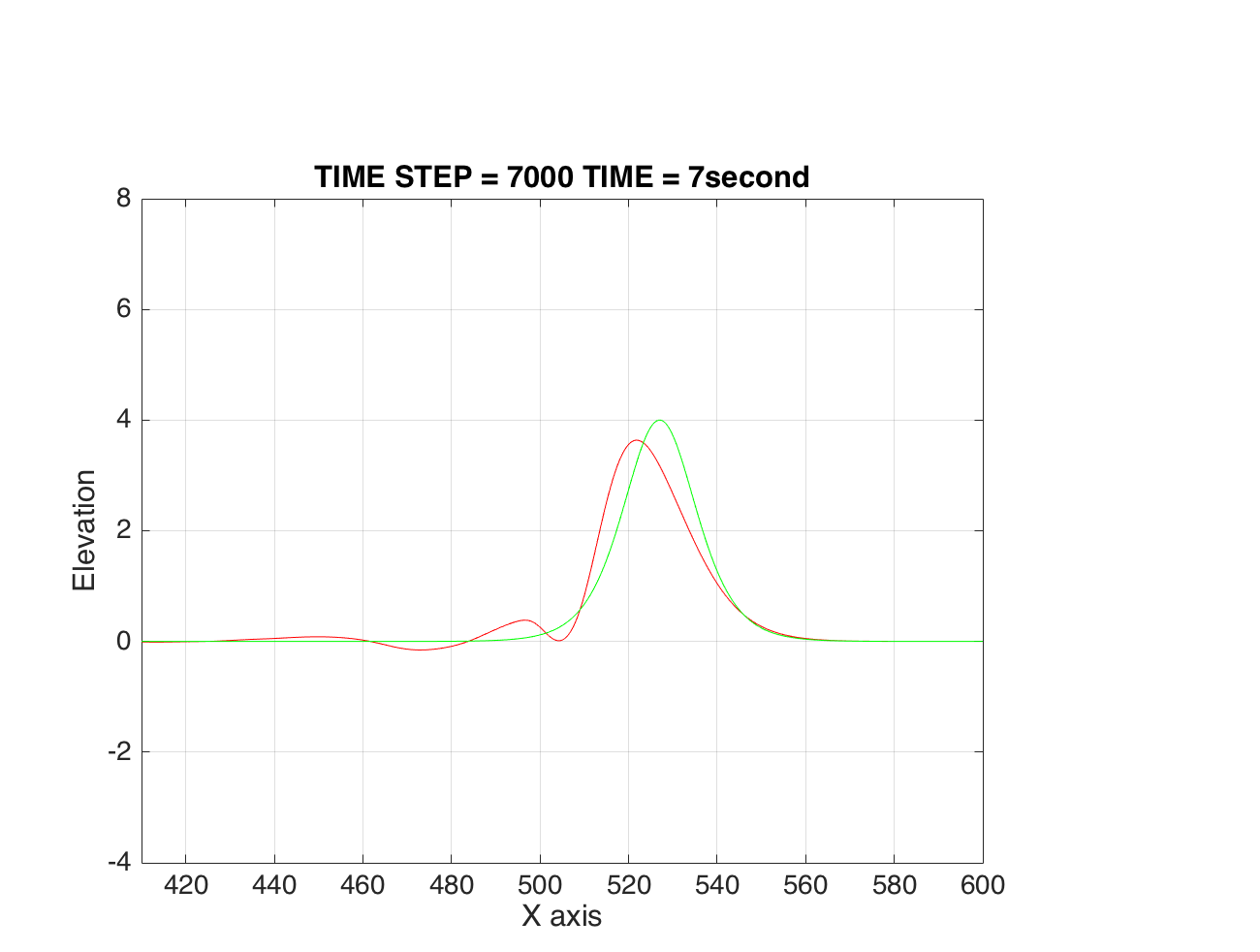}
\end{subfigure}%
\begin{subfigure}{.5\textwidth}
  \centering
  \includegraphics[width=0.95\linewidth]{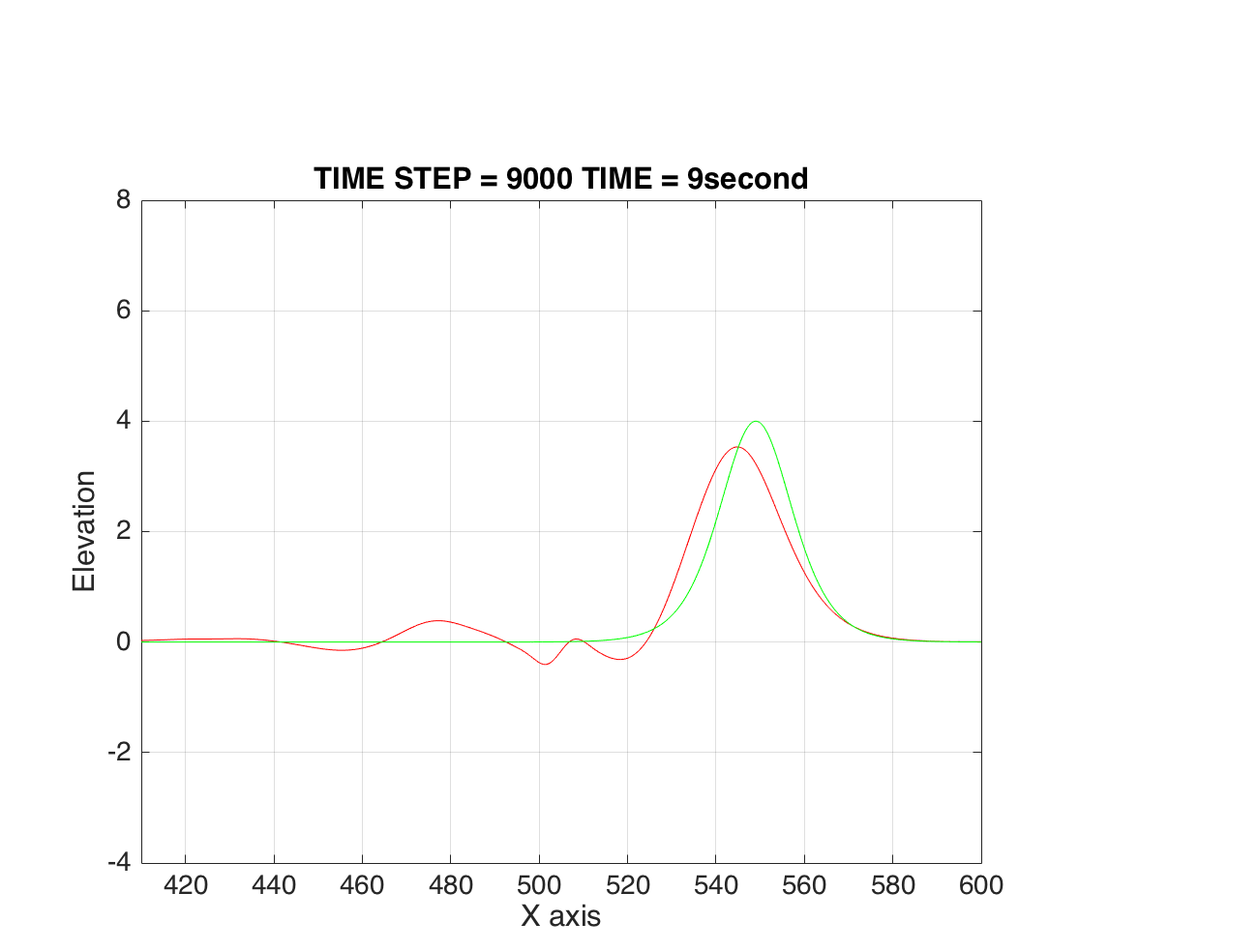}
\end{subfigure}
\caption[Evolution of a passing wave over a sliding object]{Evolution of a passing wave
  ($\mu=0.1$, $\varepsilon=0.2$) over
  a larger sliding obstacle ($\beta=0.4$, $c_{fric}=0.001$); {\color{red}red}
    curve is the passing wave,
    {\color{green}green} curve is the
    reference soliton for flat bottom}
\label{fig_waveslip}
\end{figure}

As it can be seen from Figure \ref{fig_waveslip} the approaching wave gains amplitude as it starts running up on the
slope created by the solid object on the bottom. In this
physical scenario the friction term characterizing the solid movement
becomes less important than the pressure term, resulting in the solid
sliding for a long distance based on the pressure generated by the wave. This solid motion
results in a perturbation that creates a new wavefront in front of the
solid, thus in front of the initial wave (time step $4000$),
transmitting its energy to the solid and the newly forming wavefront. 

In the meantime the solid starts
propagating with increased velocity, thus further amplifying its
amplitude (time step $5400$). When the wave peak finally passes over the top
of the solid, it drops due to the downwards slope in the bottom
topography (time step $6000$, generating secondary waves traveling
backwards (time step $6400$). Notice
that the passing wave has a significant loss in amplitude (almost $10\%$) and
an attenuated shape (time step $9000$).

\begin{figure}
\centering
\begin{subfigure}{.4\textwidth}
  \centering
  \includegraphics[width=\linewidth]{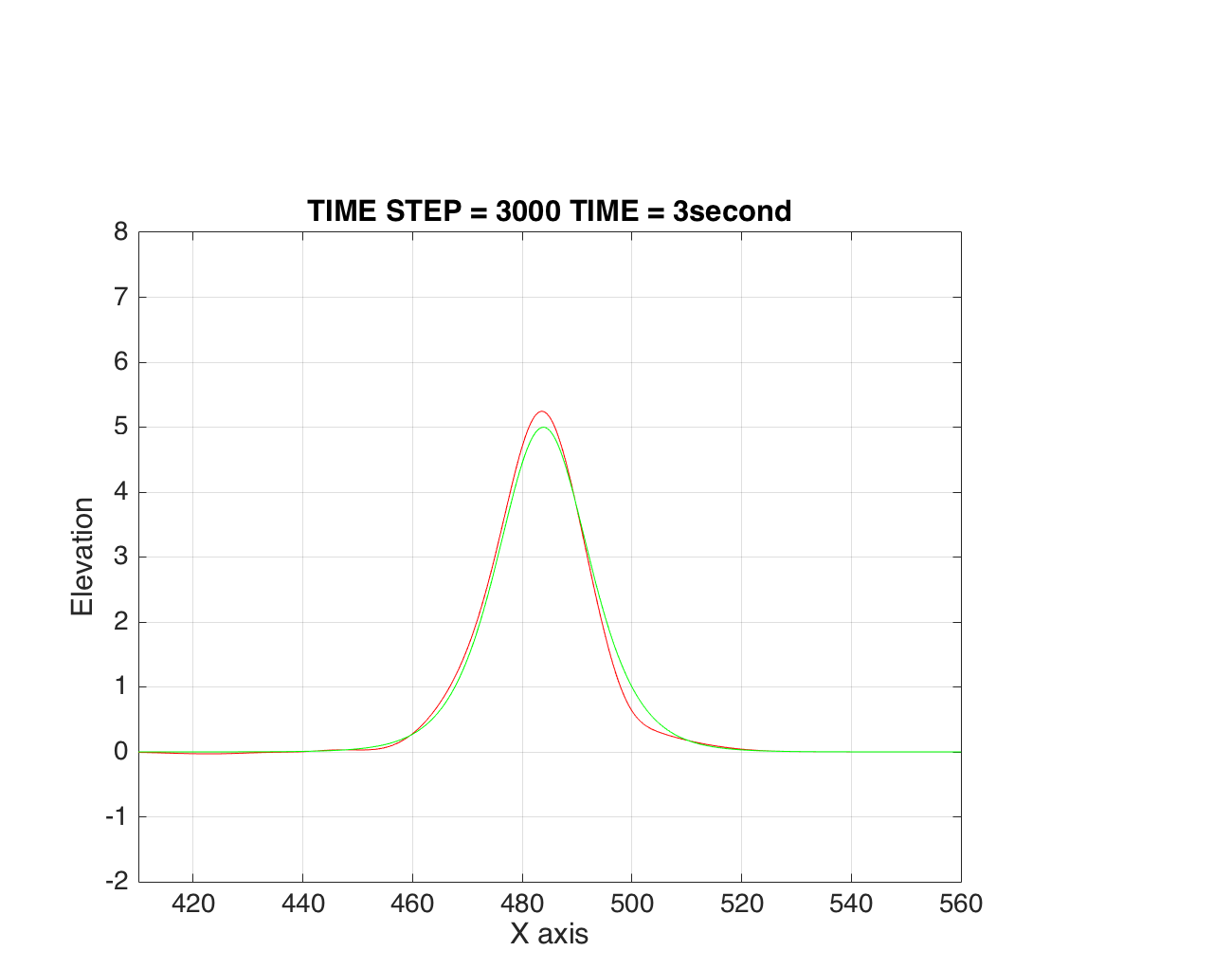}
\end{subfigure}%
\begin{subfigure}{.4\textwidth}
  \centering
  \includegraphics[width=\linewidth]{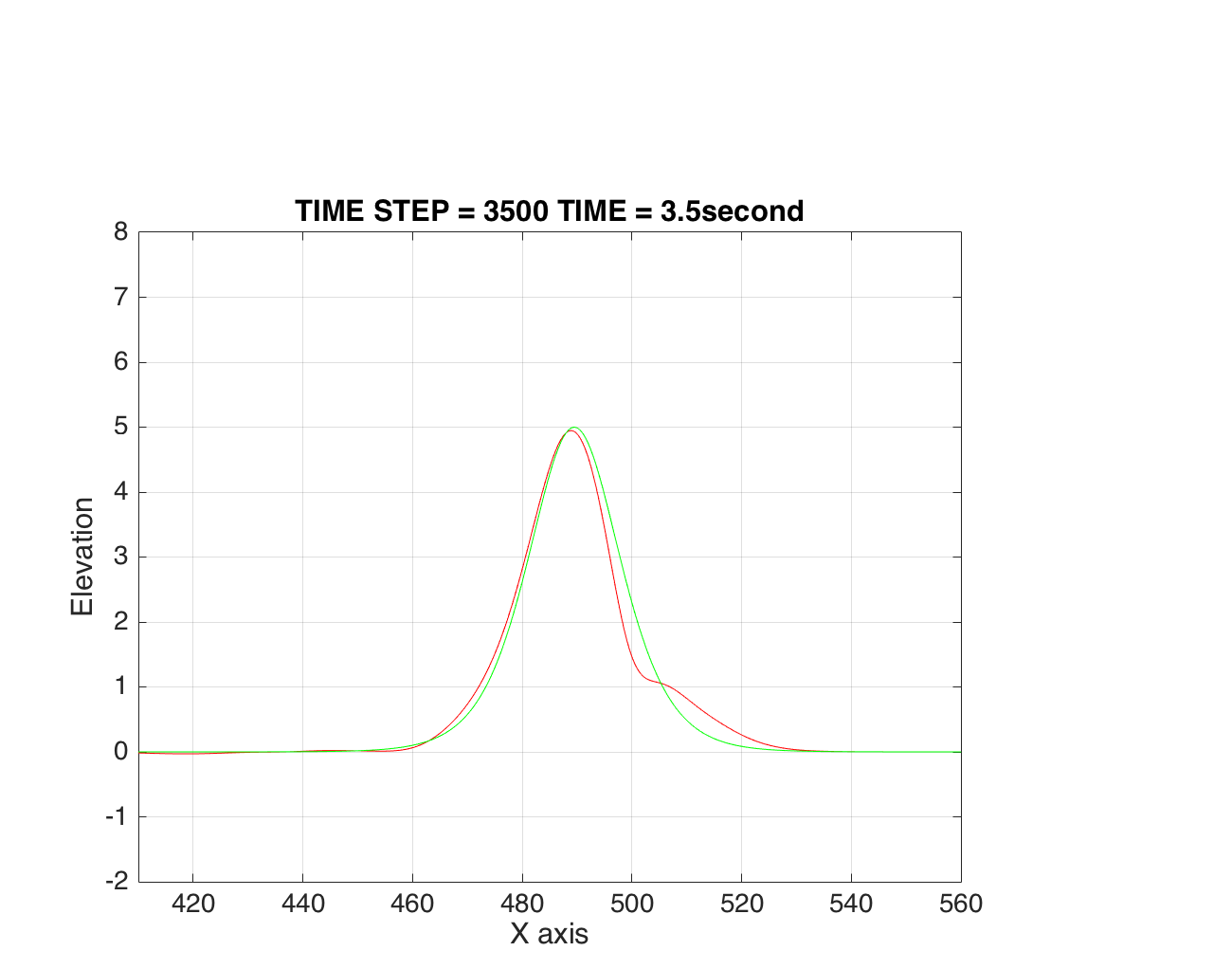}
\end{subfigure}\newline
\begin{subfigure}{.4\textwidth}
  \centering
  \includegraphics[width=\linewidth]{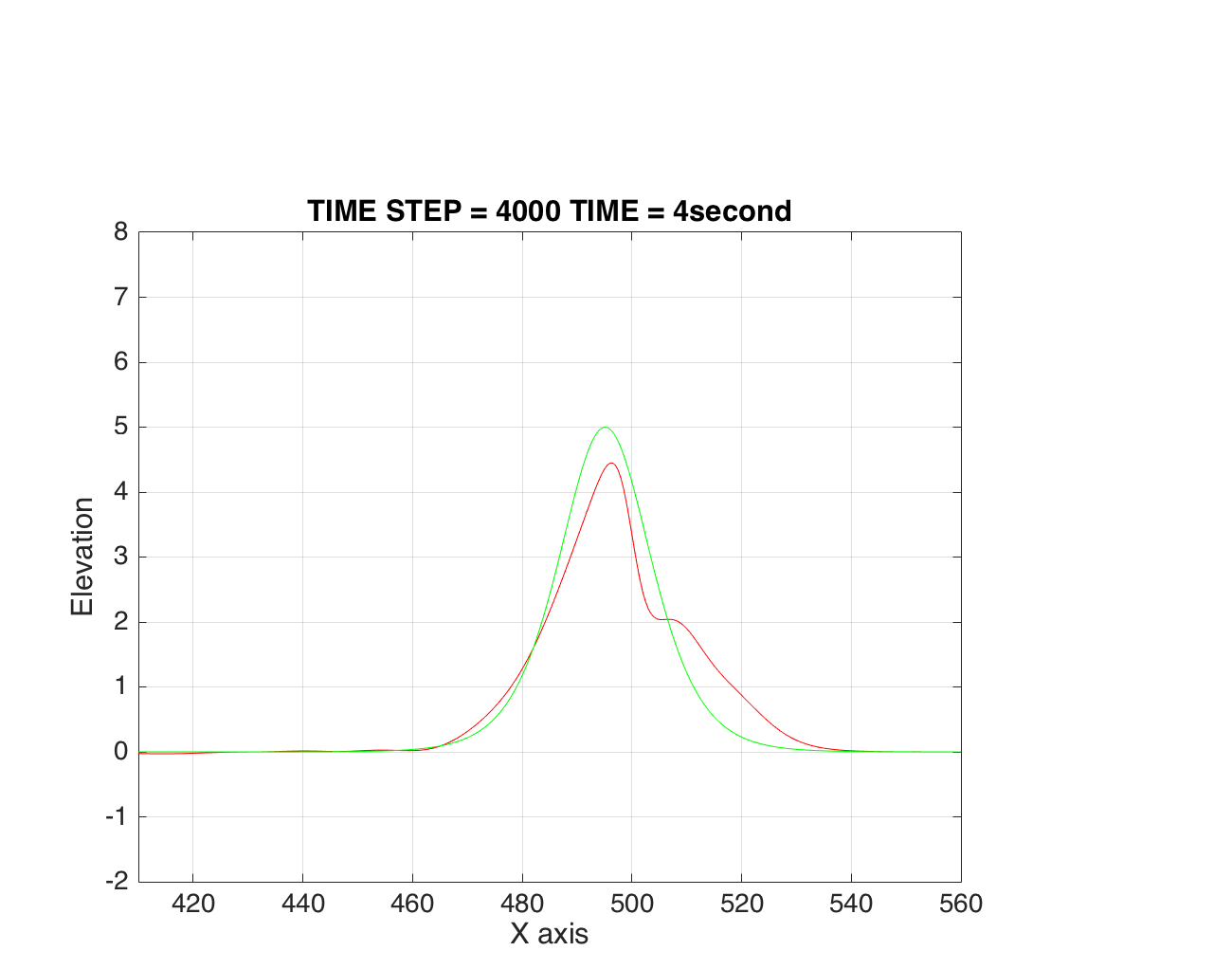}
\end{subfigure}%
\begin{subfigure}{.4\textwidth}
  \centering
  \includegraphics[width=\linewidth]{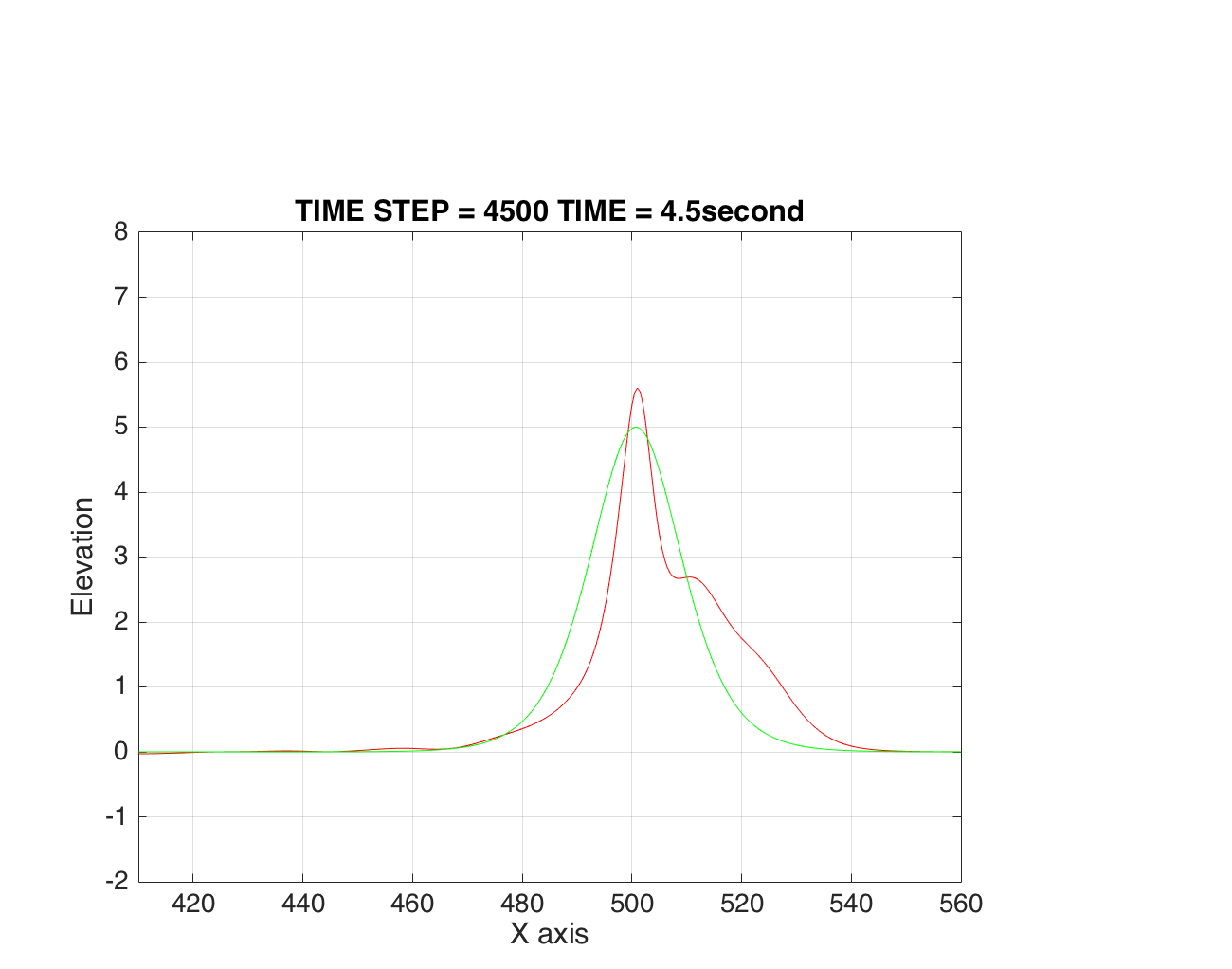}
\end{subfigure}\newline
\caption[Evolution of a surging wave]{Evolution of an approaching wave
  ($\mu=0.25$, $\varepsilon=0.25$) over an almost perfectly
  sliding large solid ($\beta=0.5$, $c_{fric}=0.001$); {\color{red}red}
    curve is the approaching wave,
    {\color{green}green} curve is the
    reference soliton for flat bottom}
\label{fig_wavesurge}
\end{figure}

Another type of wave breaking was observed in the nearly perfect,
almost frictionless ($c_{fric}\ll 0.01$)
sliding cases for large solid objects and intermediate to
large wave amplitudes (see Figure \ref{fig_wavesurge}). The approaching wave gains in amplitude due to
the slope presented by the object, further amplified by the fact that
the solid starts propagating as well. In these critical cases however,
the solid velocity becomes comparable to the passing wave's, thus
maintaining the critical position of the wave (time step $4000-4500$, the
wave peak propagates with the solid peak). A solid displacement of $5$
to $20$ meters can be detected. The critical increase in steepness is detected closer to
the front trough, implying a characteristically different wave
breaking (surging waves). 

\subsection{Amplitude variation for a passing wave}\label{sec_amplvari}

The following simulations measure the effect of a bottom topography
deformation (a solid object) on a single passing wave, most notably
the variation of the amplitude. We will compare the amplitude of a
single wave approaching the solid with the amplitude of the same
wave after having traversed the interaction zone (that is the section
of the water above the object). The main interest is to observe the
difference between the cases when the solid is allowed to move and the
case where the solid is fixed to the bottom.

A solitary wave for the flat bottom case will be taken as an initial
condition, situated two wavelengths from the solid object, traveling
towards it. Throughout the simulations a characteristic base wavelength of $L=40$
is taken for a uniform shallowness parameter of $\mu=0.25$. 

We will consider $3$ different (vertical) sizes for the solid, a small
object corresponding to $\beta = 0.1$, a medium sized object with
$\beta =0.3$ and a relatively large object for $\beta = 0.5$. For each
of these three cases we examine the qualitative effects of the
solid. We are testing three different frictional domains as well, an
almost frictionless, perfect sliding, with $c_{fric} = 0.001$, a
relatively smooth sliding for $c_{fric} = 0.01$, and a hard frictional sliding for $c_{fric}=0.5$. The results
will be compared to the two reference cases: one being the simulation
running for the same time for a flat bottom, giving a $1$ to $1$ ratio
between the entering and the leaving wave amplitude, the other one being the case when the same initial object is fixed to the ground, not being allowed to move throughout the simulation, meaning that $b(t,x)=\mathfrak{b}(x)$ independent of the time.

\begin{figure}
\centering
\begin{subfigure}{.5\textwidth}
  \centering
  \includegraphics[width=\linewidth]{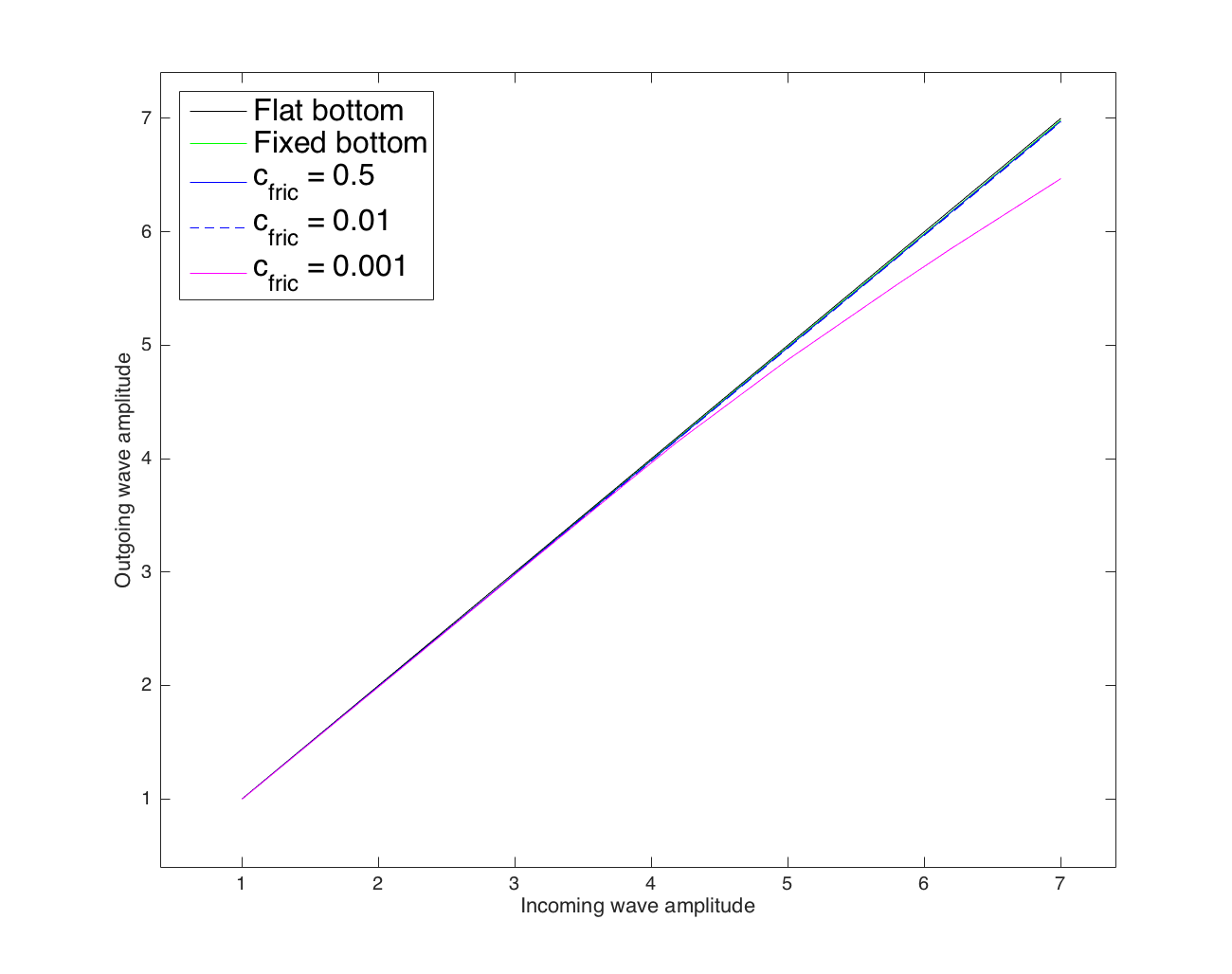}
  \caption{Change in wave amplitude over a
  small solid}
  \label{fig_avsmall1}
\end{subfigure}%
\begin{subfigure}{.5\textwidth}
  \centering
  \includegraphics[width=\linewidth]{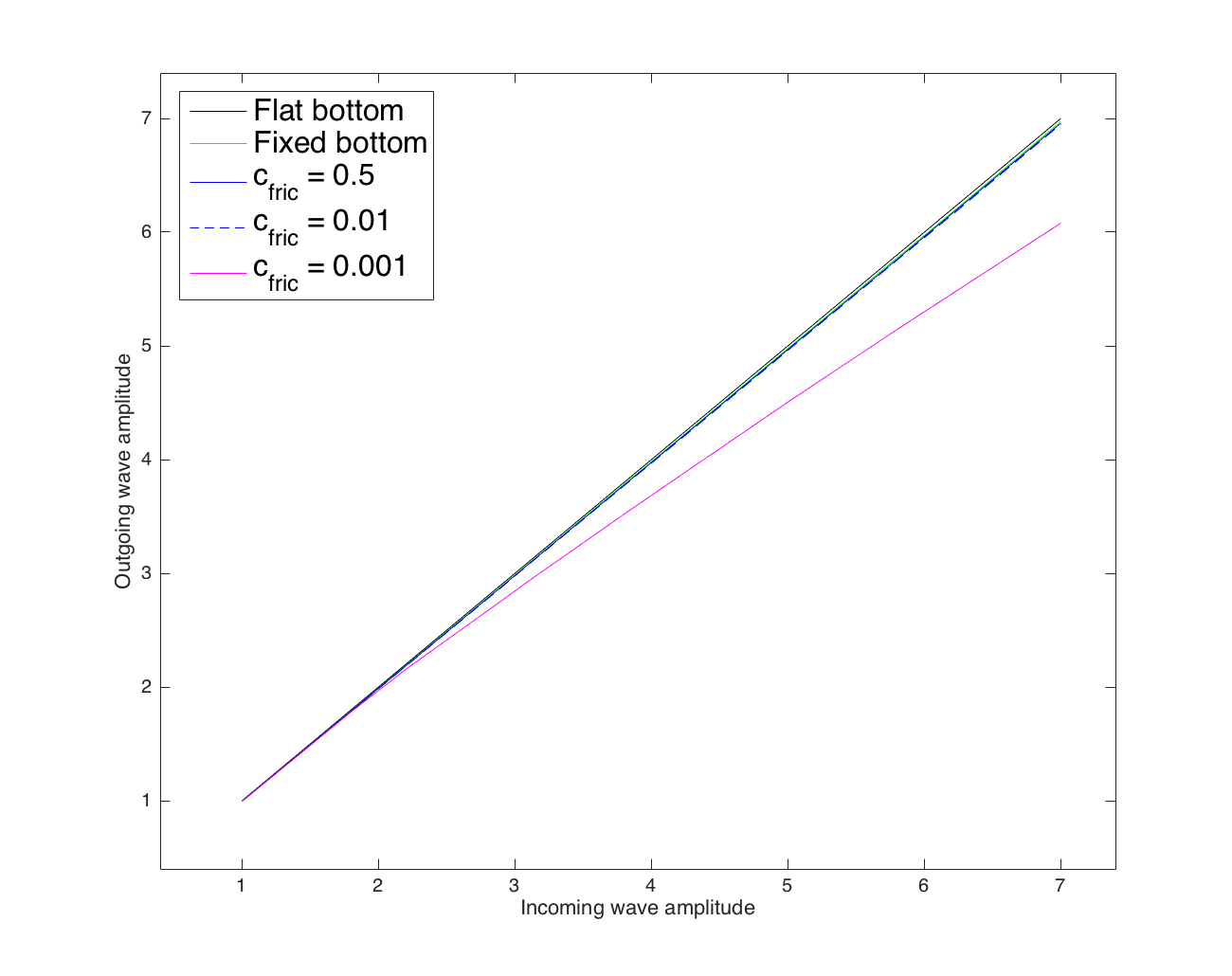}
  \caption{Change in wave amplitude over a
  medium sized solid}
  \label{fig_avmed1}
\end{subfigure}
\caption[Wave amplitude variation for a smaller object]{Wave amplitude variation for small ($\beta = 0.1$) to
  medium ($\beta = 0.3$) sized solid ($\mu=0.25$)}
\label{fig_avsmallmed}
\end{figure}

For a small object ($\beta = 0.1$, see Figure \ref{fig_avsmall1}) we
observe that the $1$ to $1$ ratio for amplitude variation for a flat
bottom is essentially preserved for both the fixed bottom and the
cases of a solid with frictional movement with relatively large frictional
coefficient. Only an almost negligible ($10^{-3}$) drop in amplitude is
detected in all these situations.

The decrease in amplitude is measurable however in the
$c_{fric}=0.001$ regime for wave amplitudes above $4$. This relative
drop in amplitude is
increased with the increase of the incoming wave height.

For the $\beta=0.3$ intermediate vertical solid scale, the simulation results are
summarized in Figure \ref{fig_avmed1}. One can notice similar effects
as the ones concluded in the previous case for a small sized
solid. Except for the almost perfect sliding, the other cases respect
the $1\colon 1$ ratio of wave amplitudes.

The impact of the
solid is more visible in the almost perfect sliding case, since the
considerable drop in amplitude is observed for smaller waves as
well. The drop in amplitude amounts for up to $12\%$ amplitude loss
for the higher waves. Here we note that for intermediate sized waves, the
initial wave is completely absorbed by the new wave produced by the
solid motion, which is a longer and much flatter wave (thus the remarkable
drop in amplitude). For small wave amplitudes, this is not observed,
since the generated pressure force on the bottom is not large enough to
create significant solid motion.

\begin{figure}
\centering
\begin{subfigure}{.5\textwidth}
  \centering
  \includegraphics[width=\linewidth]{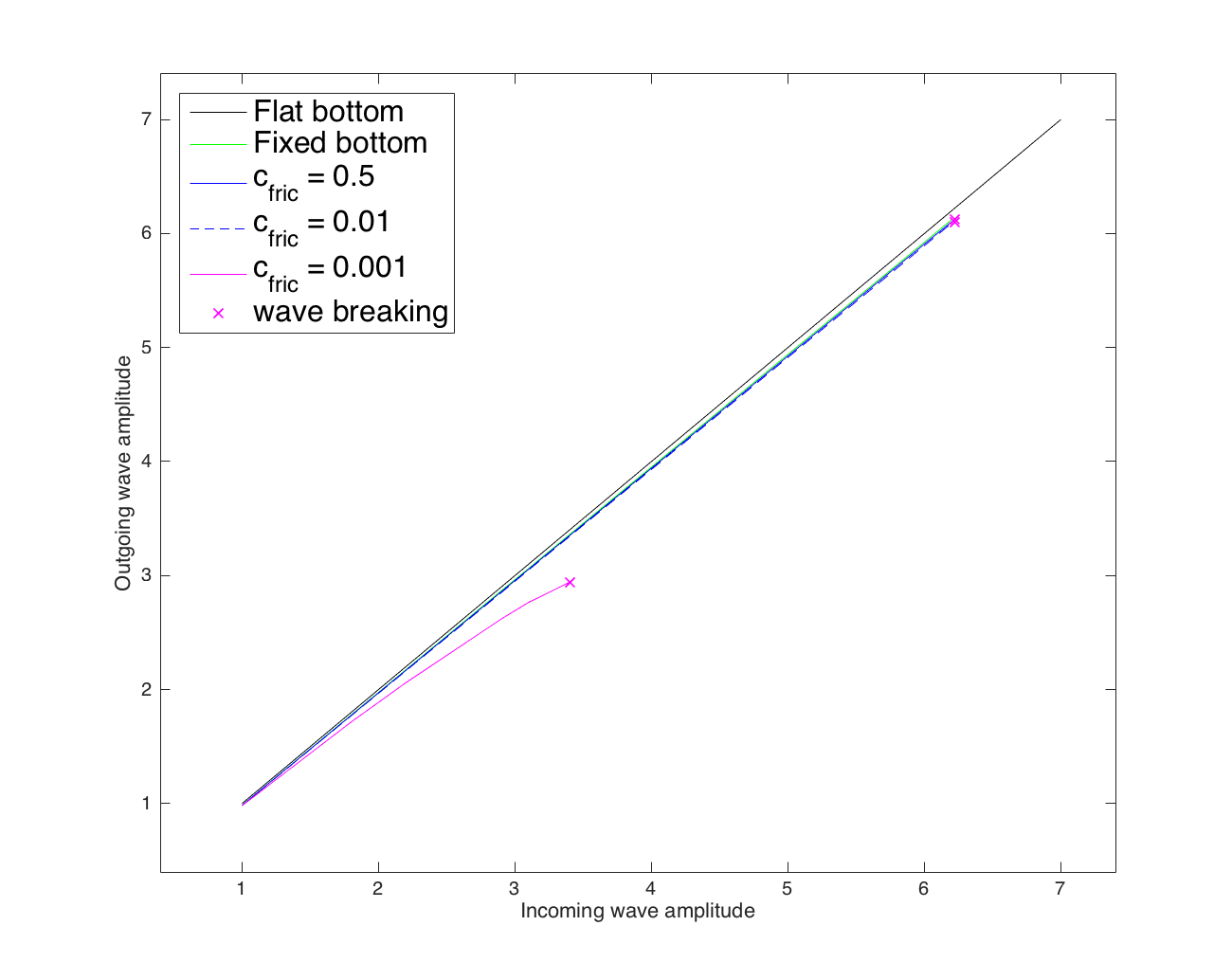}
  \caption{Change in wave amplitude for a passing solitary wave}
  \label{fig_avlarge1}
\end{subfigure}%
\begin{subfigure}{.5\textwidth}
  \centering
  \includegraphics[width=\linewidth]{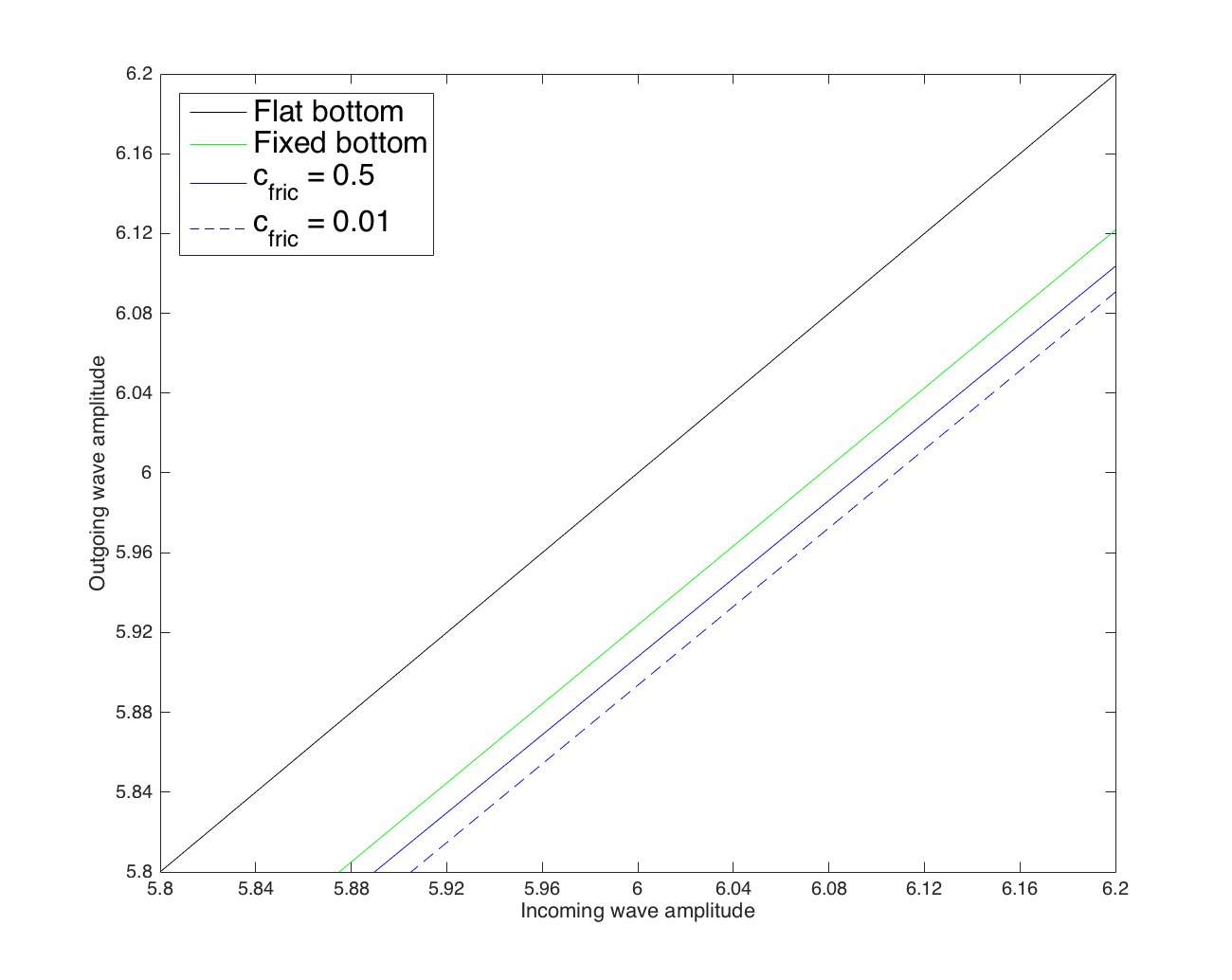}
  \caption{Zoom to the subcritical section of the amplitudes}
  \label{fig_avlarge2}
\end{subfigure}
\caption[Wave amplitude variation for a large object]{Wave amplitude variation for a large ($\beta = 0.5$) solid ($\mu=0.25$)}
\label{fig_avlarge}
\end{figure}

Finally, the results for the $\beta=0.5$ case (Figure
\ref{fig_avlarge}) indicate an even more complicated behavior. Amplitude
decrease is observed for all the non flat bottom cases, with the drop
of the amplitude being more and more important the freelier the object
can move (the less friction is imposed). The heavily frictional case
($c_{fric}=0.5$) still matches quite well the fixed bottom
scenario but the other two test cases show a more important decrease in
wave amplitude.

Another feature, made clearer in the zoomed image (Figure
\ref{fig_avlarge2}) is a slight layering between the cases in
the presence of a solid, with the fixed bottom case being the closest
to preserving the amplitude, closely followed by the hard frictional
case, signifying that these two regimes are physically close to each
other. Notice that the less important the friction is in the
system, the more the amplitude is dropping, especially for middle to
high wave amplitudes.

The main difference however is the newfound presence of wave breaking, meaning
the simulation was stopped because the numerical condition for the
initial wave breaking phase was observed (for more details, please
refer to the next section). Not only can we observe wave breaking for
all non flat bottom cases, but also there is a slight difference between the
fixed and the moving bottom cases too (further examined in the next section).

For the nearly negligible friction ($c_{fric}=0.001$) surging waves are observed, highlighting the fact
that for a larger solid in almost perfect sliding, solid motion can
also lead to critical wave transformation. The phenomenon was already described in the last part of the
previous section.

To sum it up, we observed that allowing bottom movements, especially
for an almost frictionless sliding, accounts for a measurable decrease
in wave amplitude. This is due to the fact that part of the energy of
the wave is transferred to the solid as a kinetic energy, propelling
its motion on the bottom. However, in actuality the situation is more
complicated, we will address it in Section \ref{sec_solidmov}. As for the plunging
type wave breaking and its correlation with the solid
movement, we will address it in what follows.

\subsection{An effect on the wave breaking}\label{sec_wavebreaking}

In this section we summarize the numerical test cases related to the
last part of the previous section dealing with the amplitude variation
for large solid objects where we observed (plunging type)
wave breaking for larger wave amplitudes. It manifested in the fact
that the wave became too steep, indicating that physically the wave
entered in a wave surging phase, closely followed by the breaking or
plunging of the wave. As such the simulation was stopped, even before the wave could pass over the solid.

A sufficient, but rather lenient numerical condition for wave breaking in our case (based
on \cite{southgatebraking}) is the following
\begin{equation}
\max_{i}\frac{\zeta_{i+1}-\zeta_i}{\Delta x}>1.
\end{equation}

For a more detailed analysis on wave breaking conditions for
Boussinesq type models, we refer to \cite{davidbreaking} and \cite{mariobreak}.

With a series of experiments we now examine the position of the wave
breaking point according to this criteria for different parameter
choices. We chose a maximal solid height of $H_0/2$ giving us
a topography parameter of $\beta=0.5$. For this case wave breaking
could be observed for large wave amplitudes. We implement two different frictional situations ($c_{fric} =
0.001$,  and $c_{fric}=0.5$), as well as the
reference case for a fixed bottom.

\begin{figure}
\centering
  \includegraphics[width=0.7\linewidth]{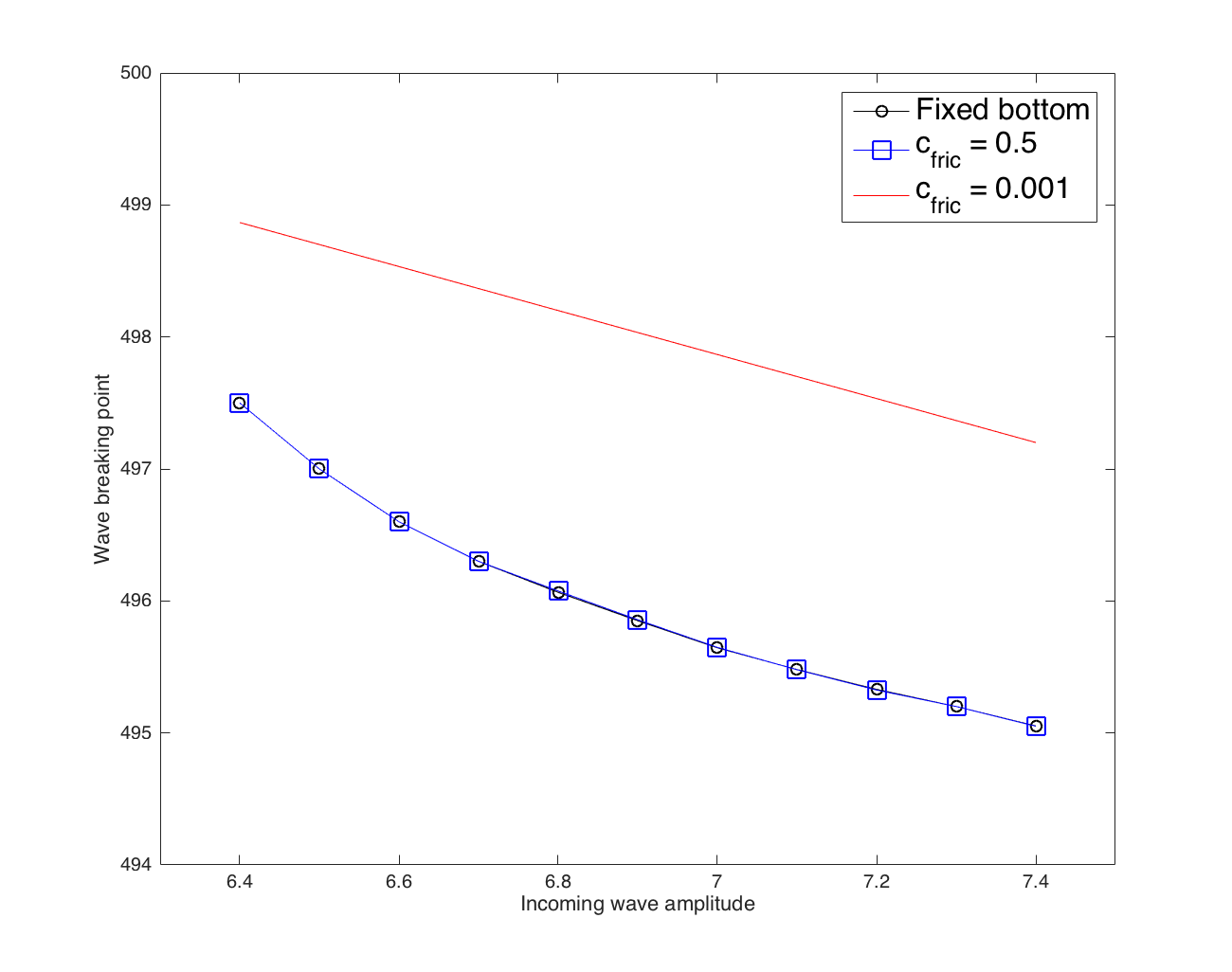}
  \caption[Wave breaking point]{Wave breaking point for large amplitude incoming waves
    ($\mu=0.25$, $\beta=0.5$)}
  \label{fig_breaker}
\end{figure}

A single traveling wave of wavelength $L=40$ is sent $2L$ distance
away from the solid and with a wave amplitude as the principal
parameter for the simulations. As a reference, the position of the
wave at the moment of the numerical wave breaking point is given by
the position of the wave crest. The simulations have also been carried out for the case of a fixed bottom, having the topography of the initial solid state.

The main interest of this numerical experiment was to see the effect of
a solid that is allowed to move on the wave breaking point. As it can be
seen from the results (Figure \ref{fig_breaker}) the $c_{fric}=0.5$ frictional case
and the fixed bottom barely differ. Actually the same was observed for
other $c_{fric}$ values as well (in the range of $0.01$ to $1$). 

In the case when the object is sliding almost frictionlessly
($c_{fric}=0.001$) however, the position of the wave crest is much
further away from the initial position. The
qualitative difference is due to the change of the nature of the wave
breaking (as it was remarked in Section
\ref{sec_singli}). Nevertheless a measurable delay is observed for the
wave breaking point, owing to the fact that the initial wave loses
some of its energy while the new, frontal wave is created by the solid.

\subsection{Observations on the hydrodynamical damping}

The previous simulations gave us some insight on the effects of a
moving solid on the wave motion. Now we reverse our point of view, so
to speak, in order to examine the inverse, that is how the waves
affect the solid motion.

\begin{figure}
\centering
\begin{subfigure}{.5\textwidth}
  \centering
  \includegraphics[width=\linewidth]{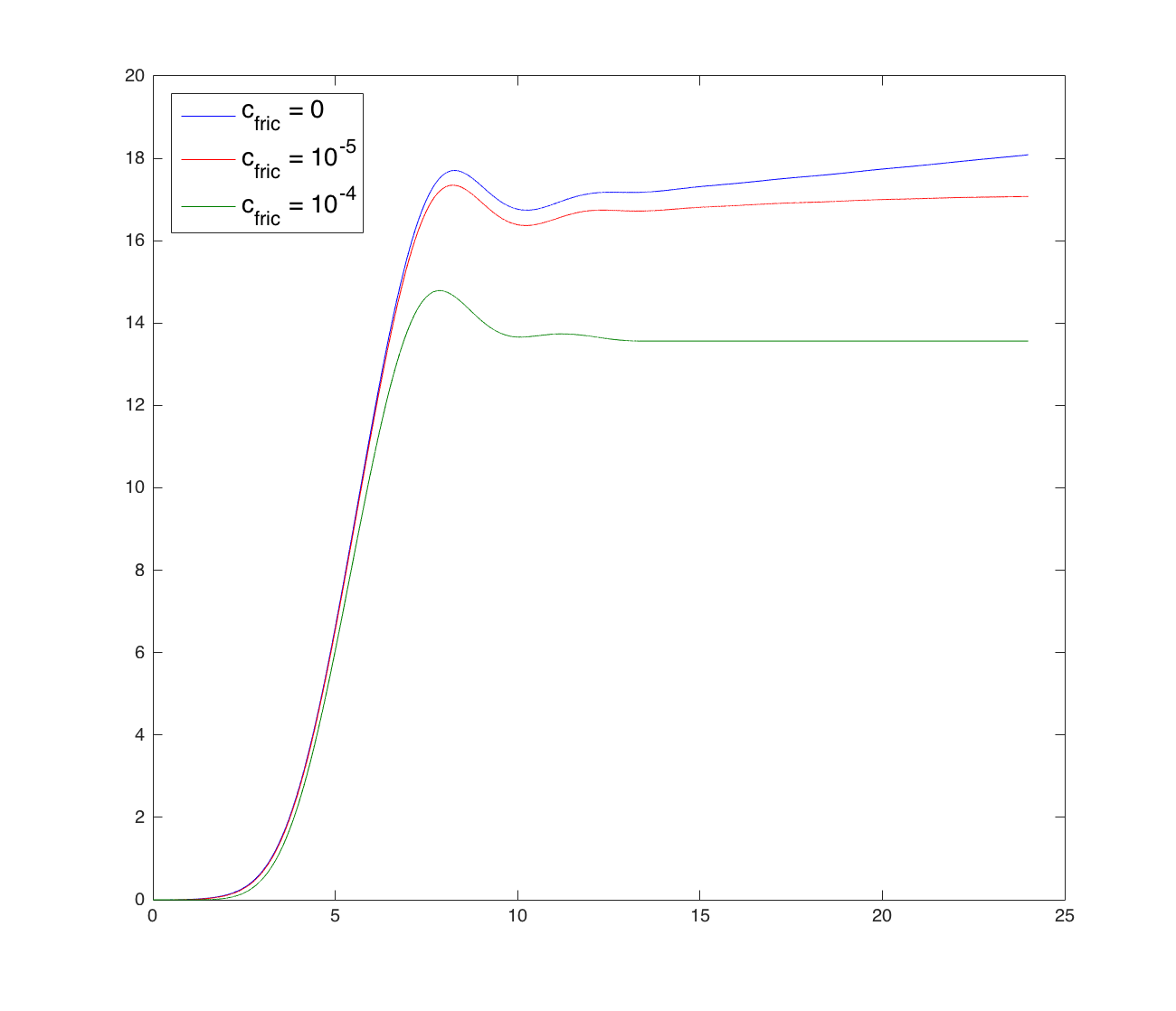}
  \caption{Solid position}
  \label{fig_posilowfric}
\end{subfigure}%
\begin{subfigure}{.5\textwidth}
  \centering
  \includegraphics[width=\linewidth]{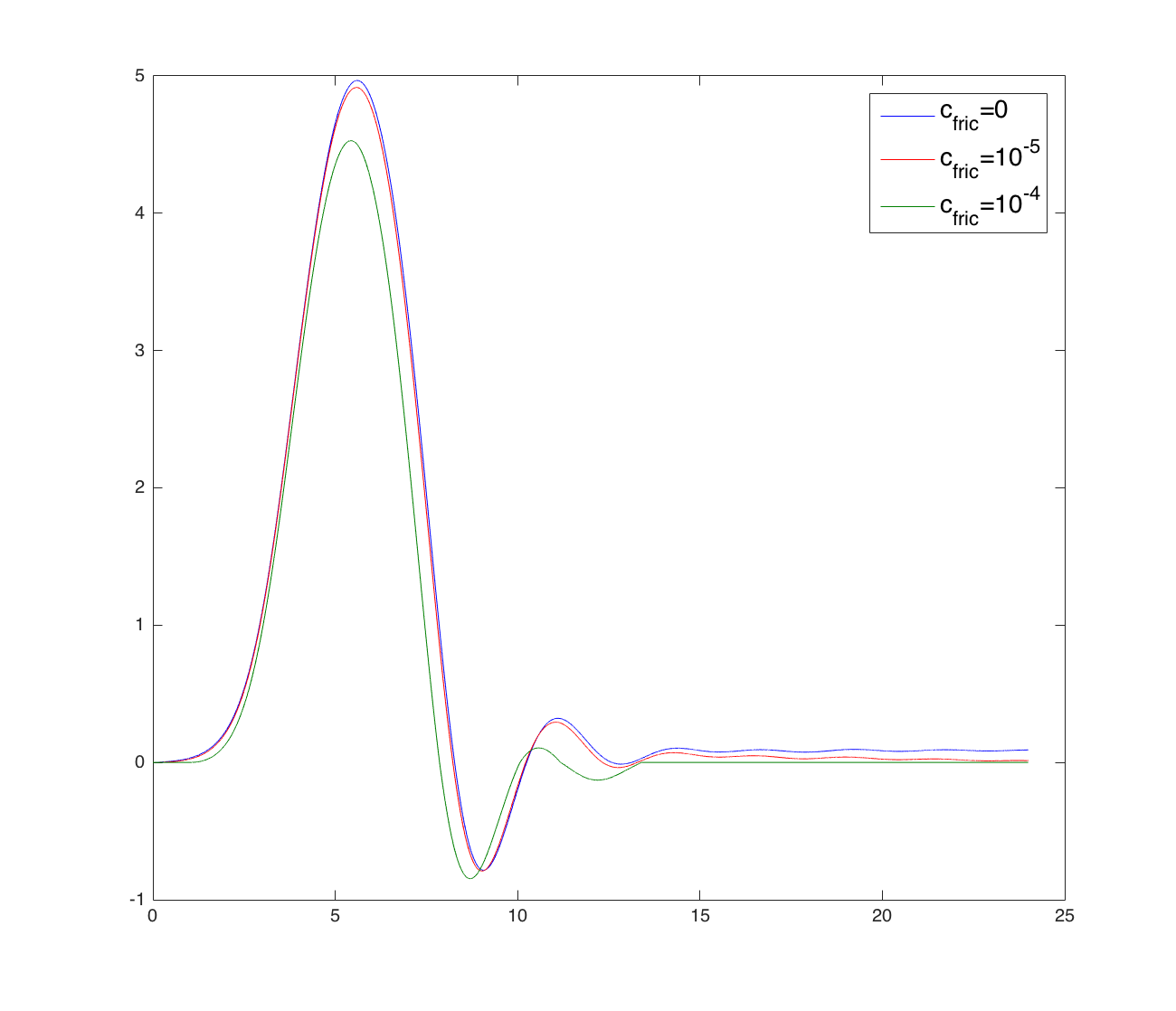}
  \caption{Solid velocity}
  \label{fig_velolowfric}
\end{subfigure}
\caption[Solid motion for low coefficients of friction]{Solid motion
  for low coefficients of friction ($\mu=\varepsilon=0.2$, $\beta=0.4$)}
\label{fig_lowfric}
\end{figure}

Our numerical simulations show
that by choosing $c_{fric}$ of order $10^{-4}\sim 10^{-5}$ the solid
takes and extended amount of time before finally coming to a halt
(Figure \ref{fig_lowfric}), creating small amplitude, large wavelength waves in
front of and behind it. By creating such waves, its motion is damped
by a phenomenon similar to the dead-water phenomenon, described in detail for Boussinesq type models for example
in \cite{vincentintern}.

In the limiting situation of $c_{fric}=0$ it continues its movement
without stopping. Notice the rapidly stabilized oscillatory behaviour
in the velocity profile of the solid (Figure \ref{fig_velolowfric}), further highlighting
the small amplitude wave-generation around the peak of the solid.

\begin{figure}
\centering
\begin{subfigure}{.5\textwidth}
  \centering
  \includegraphics[width=\linewidth]{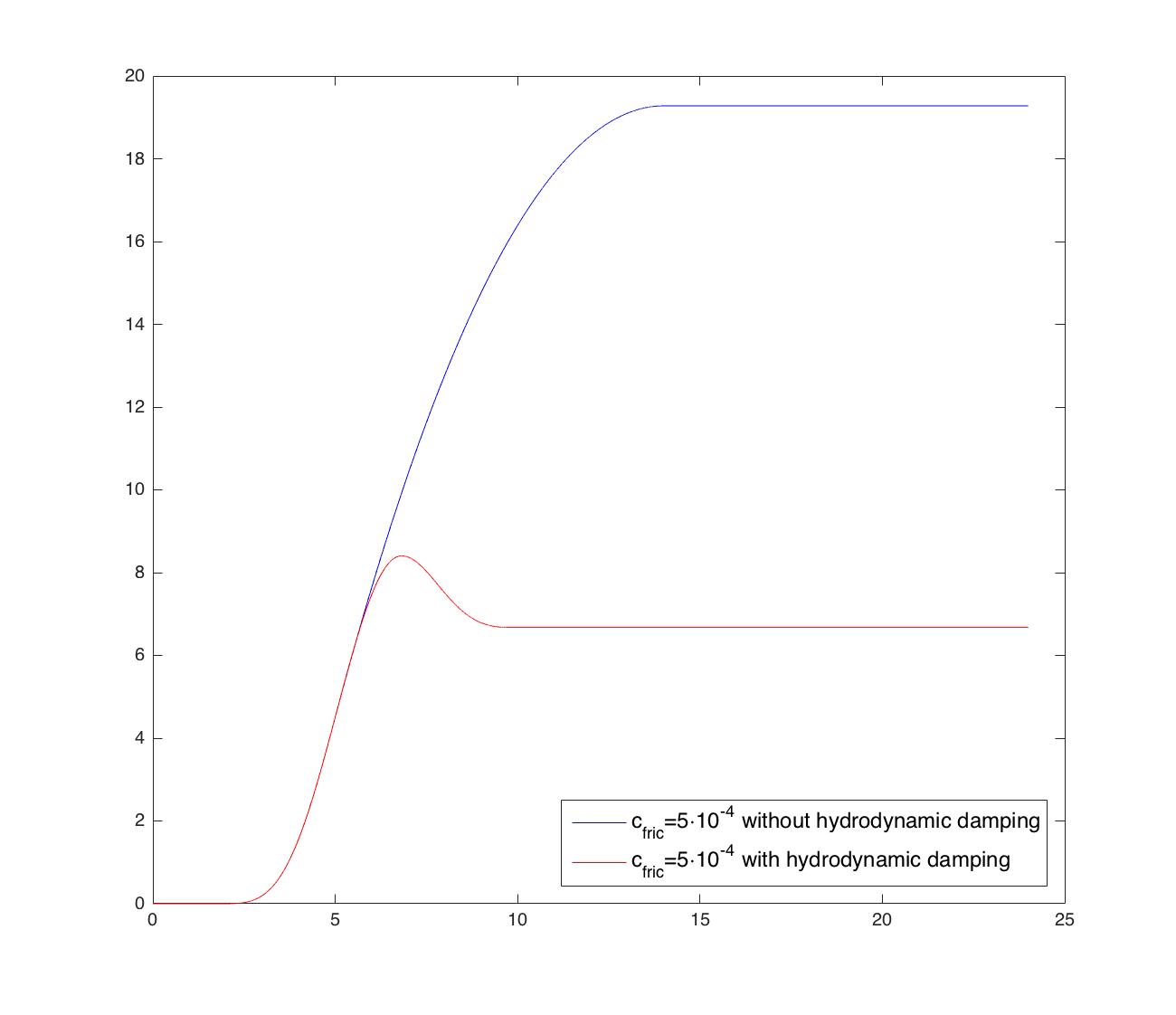}
  \caption{Solid position}
  \label{fig_posihydro}
\end{subfigure}%
\begin{subfigure}{.5\textwidth}
  \centering
  \includegraphics[width=\linewidth]{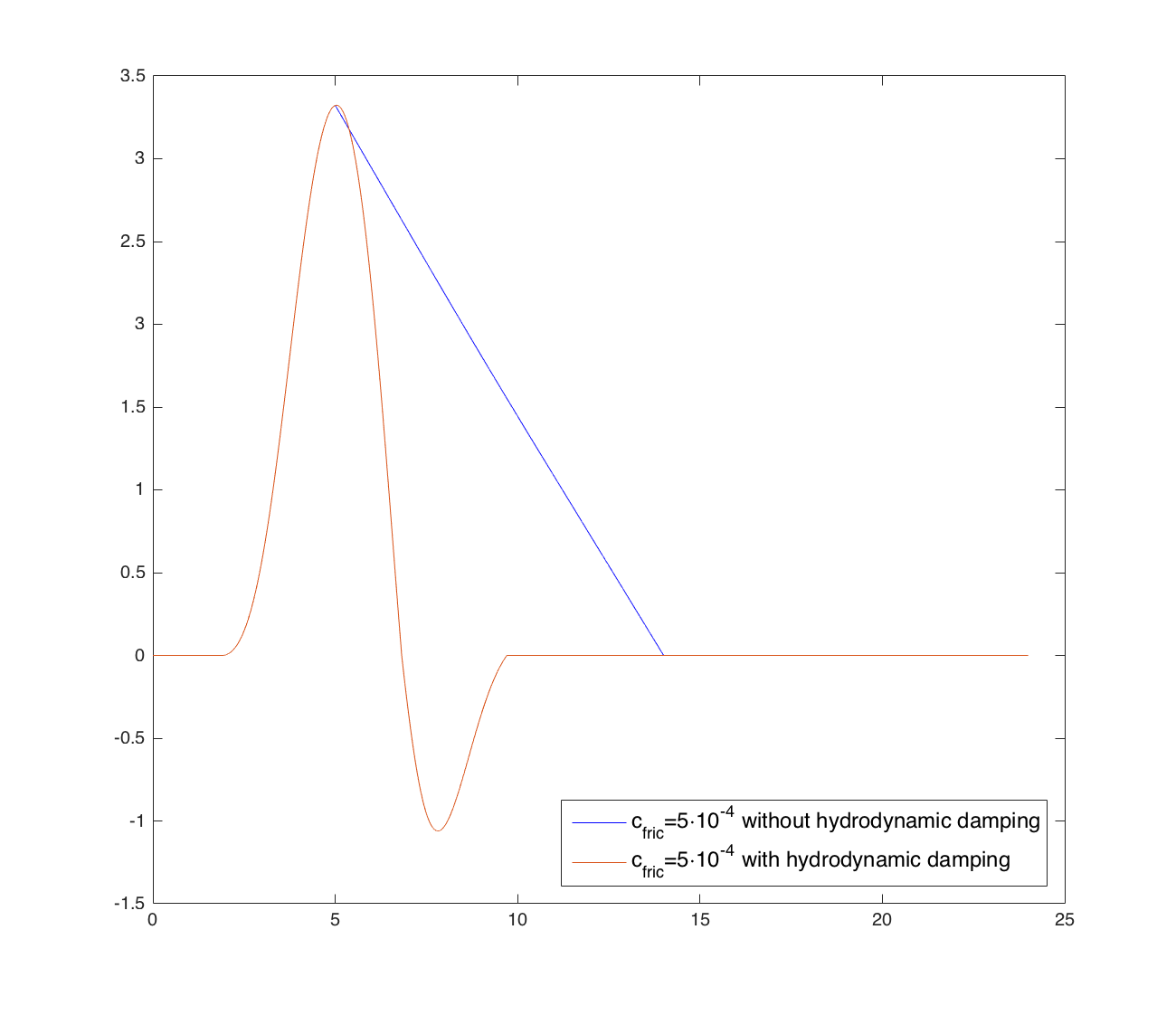}
  \caption{Solid velocity}
  \label{fig_velohydro}
\end{subfigure}
\caption[Solid motion with and without hydrodynamic damping]{Solid
  motion with and without hydrodynamic damping ($\mu=\varepsilon=0.2$,
  $\beta=0.4$, $c_{fric}=0.0005$)}
\label{fig_hydro}
\end{figure}

The hydrodynamic effects are shown not only by the change in direction
for the solid motion as well as the increased changes in velocity. In
Figure \ref{fig_hydro} we compared two situations: as a reference the standard model
was left running for the whole time ($T=24$), for the second simulation the
pressure term was removed at the moment when the velocity hit its
maximal value, taking into consideration only the frictional damping
of the system.

As it can be clearly seen, without the hydrodynamic effects the solid
velocity decreases essentially linearly, as it is dictated by the
corresponding equation of motion, the solid slowly comes to a
halt. The hydrodynamic damping not only increase the deceleration
process but it also keeps the solid relatively close to its initial position.

\subsection{Measurements of the solid displacement}\label{sec_solidmov}

Three additional sets of simulations were conducted in order to
measure the effects of friction on the system, as well as to
highlight long term effects by means of simulating an approaching
wave train.

\subsubsection{Solid displacement for varying coefficients of friction}

\begin{figure}
\centering
  \includegraphics[width=0.7\linewidth]{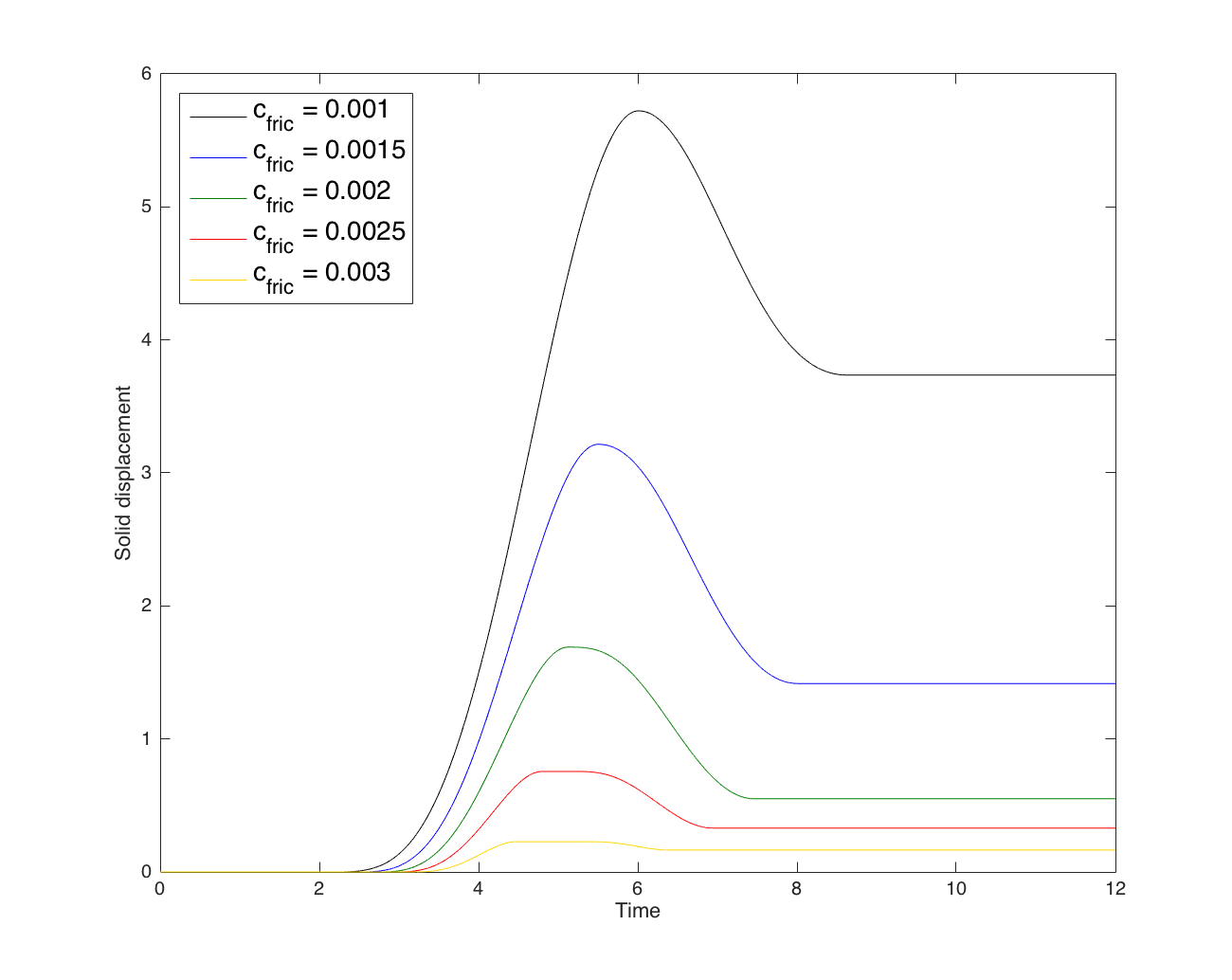}
\caption[Evolution of the solid position for varying
coefficients of friction]{Solid displacement for a varying coefficients of friction ($\mu=\varepsilon=0.25$, $\beta=0.3$)}
\label{fig_friction}
\end{figure}

The first set of tests was conducted in tandem with the observations
on the wave surging for almost perfect sliding. For the physical
parameters $\beta = 0.3$ and $\varepsilon = 0.25$ the friction
coefficient is varied from $0.001$ to $0.003$. Figure
\ref{fig_friction} shows the change in the solid motion for this case.

We remark a significant drop in the maximal solid displacement, attributed
to the qualitative change of the system from a frictional to a
frictionless sliding. The critical values for $c_{fric}$ for
this transformation depend on the shape of the object as well as its
physical parameters.

Furthermore we remark that the solid does in fact stop in its motion, a property
due to not only the damping effect of the frictional forces but to the
hydrodynamic damping as well. The latter one is clearly visible by the
fact that the solid motion changes direction after the wave peak
passes over it and rapidly looses velocity by the time the wave leaves
the interaction zone.

\subsubsection{Solid displacement for varying wave amplitudes}

The second
one consists of measuring the effect of a single
traveling wave on the solid motion, for a relatively frictionless sliding
($c_{fric} = 0.001$). The maximal vertical size of the solid is $6$, with
waves having an amplitude in the range $[3, 5]$. It also serves as a reference case for the wave train simulations
presented later on.

\begin{figure}\centering
\includegraphics[width=0.7\linewidth]{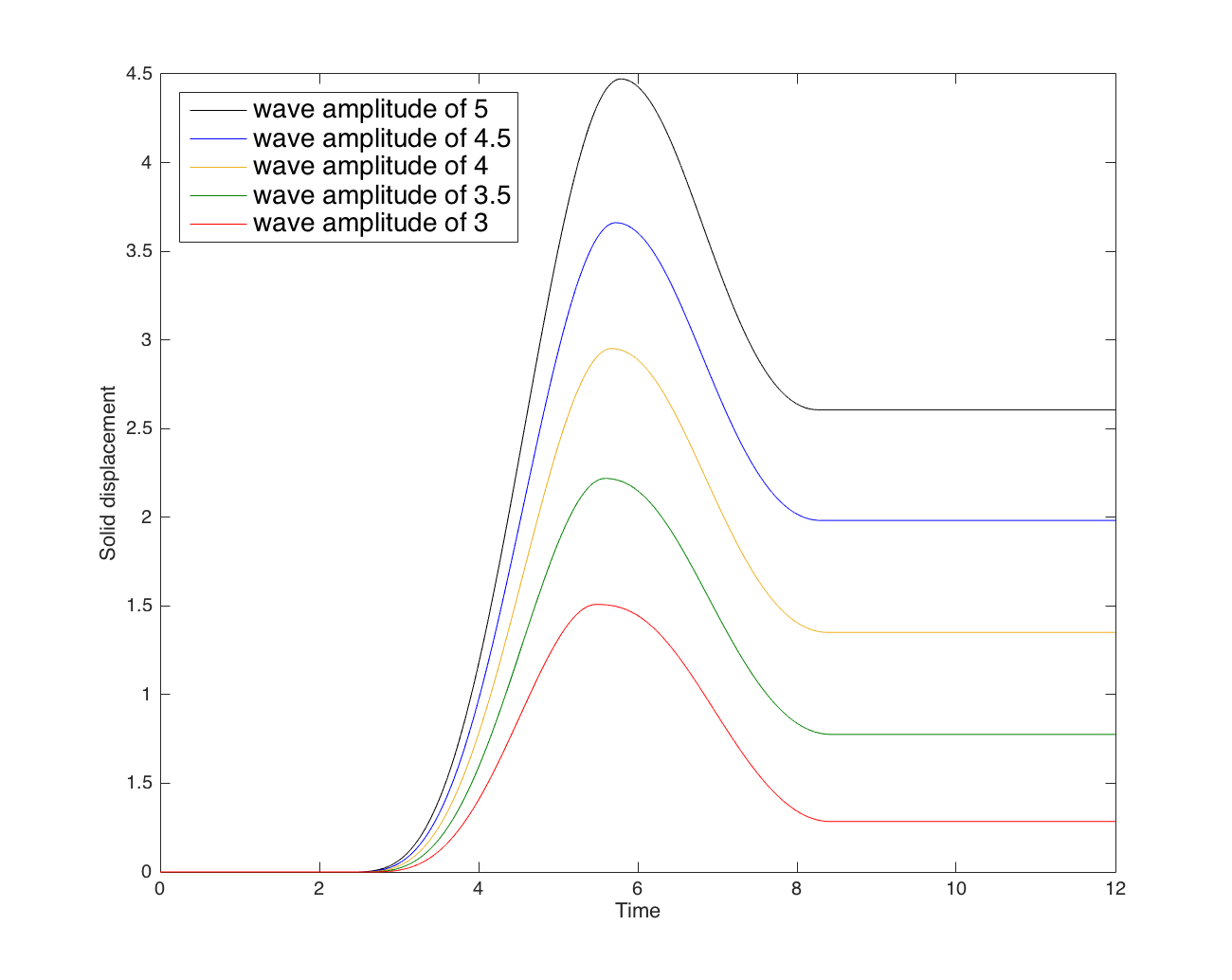}
\caption[Evolution of the solid position for a single passing wave]{Evolution of the solid position due to a single
  solitary wave, with frictional sliding ($\mu=0.25$, $\beta=0.3$, $c_{fric}=0.001$)}
\label{fig_trainref}
\end{figure}

On Figure \ref{fig_trainref} we can see the
almost perfect ``oscillation circle'' for the solid position, meaning
that after the approaching wave pushes the solid forward, it passes
over the object, and then it pushes the solid backwards, making it
return nearly to its initial state. Notice the slight asymmetry of the
curves as well as the rather extended calming phase, attributed mainly
to the slowing effects of frictional forces.

The final set of numerical experiments concerns the long term effects
of wave motion on a solid that is allowed to move freely, subjected to
a frictional sliding on the bottom of the wave tank. This was carried out by simulations on a long time scale, it involves sending a wave train consisting of $10$ consecutive solitary waves in the direction of the solid and measuring the evolution of the solid displacement.

The wave tank is now taken to be twice as large as before, with a length of $2000$ to properly
accommodate the wave train. A medium sized solid is chosen ($\beta = 0.3$) for intermediate wave
amplitudes ($\varepsilon$ ranging from $0.15$ to $0.25$) with a
wavelength of $L=40$. The simulation is run for a time of $T=54$,
allowing for $7$ waves to pass over the solid.

\begin{figure}\centering
\includegraphics[width=0.75\linewidth]{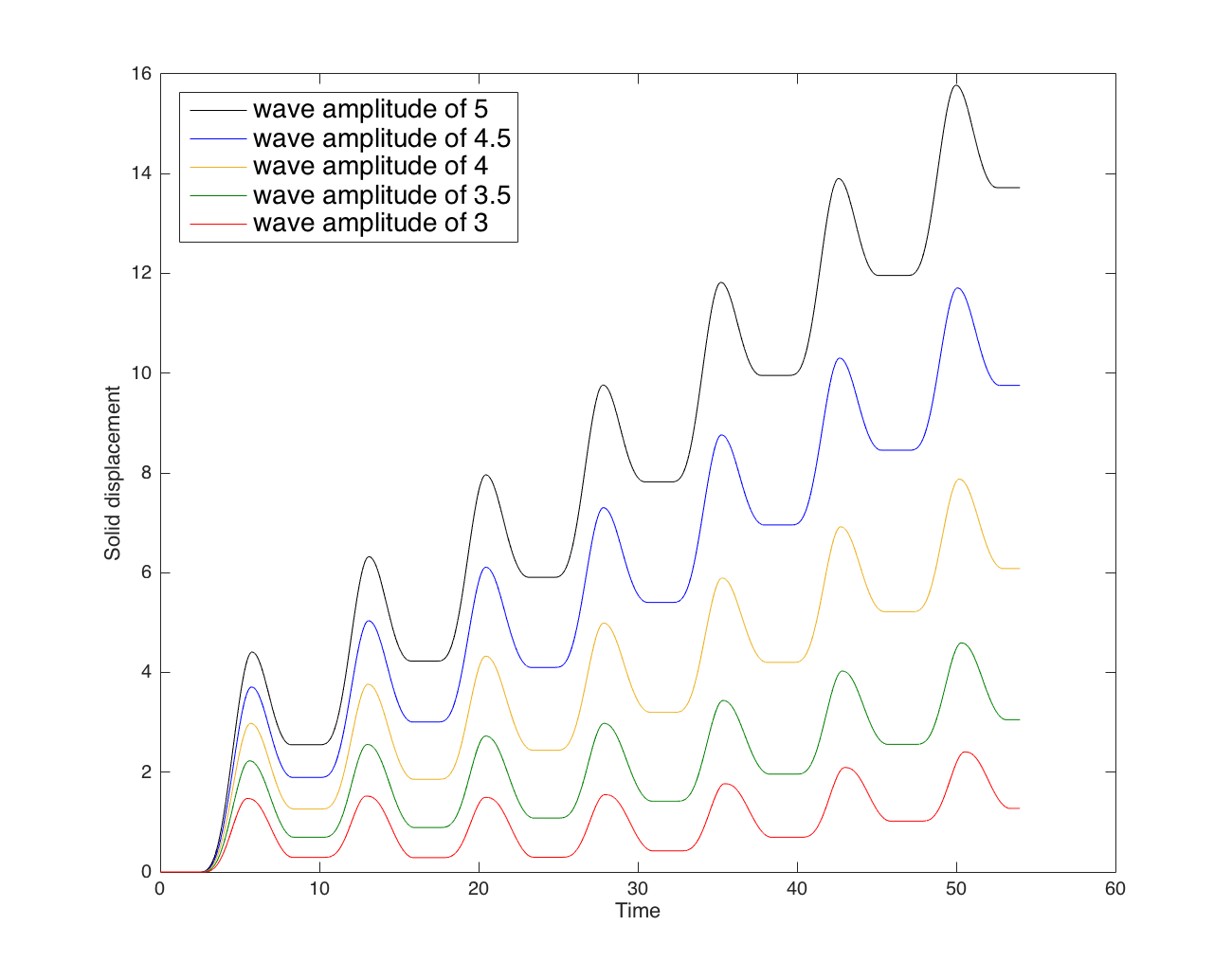}
\caption[Evolution of the solid position for multiple passing waves]{Evolution of the solid position due to a wave
  train, with almost perfect sliding ($\mu=0.25$, $\beta=0.3$, $c_{fric}=0.001$)}\label{fig_trainsmall}
\end{figure}

The evolution of the position of the center of the mass is plotted in
Figure \ref{fig_trainsmall}, corresponding to an almost perfect
sliding ($c_{fric} = 0.001$). We observe solid displacement of order
$1$ for each passing wave, as well as some qualitative differences in the wave
cycles.

There is a positive net solid displacement per incoming wave, resulting in the solid
getting further and further away from the initial position after each
passing wave in the train. We remark that the displacement per passing
wave is not constant, attributing to an overall nonlinear increase in
distance. This is due to the fact that a passing wave is perturbed by
the solid, resulting in backwards traveling small amplitude long
waves (as pointed out in Section \ref{sec_singli}), meaning that later
members of the wave train get perturbed even before reaching the
interaction zone. This results in a non-constant incoming wave
amplitude even though initially it was constant.

Also, comparing the effects of a single wave (Figure \ref{fig_trainref}) with a
wave train (Figure \ref{fig_trainsmall}), we can
notice that, even though the elements of the wave train were launched
sufficiently apart to avoid undesired interactions with each other, due to the extended period of time it takes for the
frictional and hydrostatic damping to slow down the solid motion,
nonlinear effects are observable in the solid displacement.

\section{Conclusion}

In the present paper, we complemented the author's previous, mostly
theoretical work (\cite{egoboost}) by a numerical analysis of relevant
shallow water situations. A staggered grid based finite difference
method incorporating an Adams predictor-corrector algorithm (\cite{linman}) was
adapted to the coupled Boussinesq equations. A numerical verification
indicates a convergence of order three in space and an order two
in time variable. We not only showed that the motion of the solid
reduces the wave amplitude in general, but also that it affects the wave breaking. We observed the presence of
surging type wave breaking for low friction regimes as well as
nonlinear interactions between the wave propagation and the solid displacement.

\bibliographystyle{alpha}
\bibliography{biblio}

\end{document}